\documentclass[AMS,Times1COL]{Template/WileyNJDv5} %STIX1COL,STIX2COL,STIXSMALL
\articletype{Article}%

\received{Date Month Year}
\revised{Date Month Year}
\accepted{Date Month Year}
\journal{}
\volume{00}
\copyyear{2023}
\startpage{1}

\usepackage{orcidlink}
\usepackage{nicematrix}
\usepackage[capitalise]{cleveref}
\usepackage[colorinlistoftodos, textwidth=3.5cm]{todonotes} 
\setuptodonotes{inline}
\usepackage{comment}
\newtheorem{Theorem}{Theorem}

\newtheorem{example}{Example}
\newtheorem{Remark}{Remark}
\crefname{construction}{Construction}{Constructions}
\usepackage{subcaption}
\usepackage{lineno}
\begin{document}
\large

\title{Combinatorial solutions to the Social Golfer Problem and Social Golfer Problem with Adjacent Group Sizes}

% % Author Orchid ID: enter ID or remove command
% \newcommand{\orcidA}{\orcidlink{0000-0002-0941-1717}} % Add \orcidA{} behind the author's name - ALICE
% \newcommand{\orcidB}{\orcidlink{0000-0003-1116-875X}} % Add \orcidB{} behind the author's name - Ivaylo
% \newcommand{\orcidD}{\orcidlink{0000-0002-3632-9612}} % Add \orcidD{} behind the author's name - Julian

% \author[1]{Alice Miller \orcidA{} }
\author[1]{Alice Miller}

% \author[1]{Ivaylo Valkov \orcidB{} }
\author[1]{Ivaylo Valkov}

% \author[2]{{R. Julian R.}  Abel \orcidD{} }
\author[2]{{R. Julian R.}  Abel}

\authormark{MILLER \textsc{et al.}}
\titlemark{Combinatorial solutions to the Social Golfer Problem and Social Golfer Problem with Adjacent Group Sizes}

\address[1]{\orgdiv{School of Computing Science}, \orgname{University of Glasgow}, \orgaddress{\state{Glasgow}, \country{United Kingdom}}}

\address[2]{\orgdiv{School of Mathematics and Statistics}, \orgname{UNSW Sydney}, \orgaddress{\state{Sydney}, \country{Australia}}}

% \address[3]{\orgdiv{Department Name}, \orgname{Institution Name}, \orgaddress{\state{State Name}, \country{Country Name}}}

\corres{Corresponding author Alice Miller. \email{alice.miller@glasgow.ac.uk}}

\presentaddress{This is sample for present address text this is sample for present address text.}

%\fundingInfo{Text}
%\JELinfo{ejlje}

\abstract[Abstract]{Resolvable combinatorial designs including Resolvable Balanced Incomplete Block Designs, Resolvable Group Divisible Designs, Uniformly Resolvable Designs and Mutually Orthogonal Latin Squares and Rectangles are used to construct optimal solutions to the Social Golfer problem (SGP) and the Social Golfer problem with adjacent group sizes (SGA). An algorithm is presented to find an optimal solution in general, and a complete set of solutions is provided for up to 150 players.} 

% Keywords % This is limited to 7

 % AMS subject classification: 05B05

\maketitle

\renewcommand\thefootnote{\fnsymbol{footnote}}
\setcounter{footnote}{1}

%%%%%%%%%%%%%%%%%%%%%%%%%%%%%%%%%%%%%%%%%%
\section{Introduction}\label{sect:intro}
The study of combinatorial designs \cite{codihandbook, hughesPiper_book_1985} is a well-established area of mathematical research.
%, with many dedicated journals (e.g. Journal of Designs and Codes \cite{jDesignsCodes} and the Journal of  Combinatorial Designs \cite{jCombDes}). 
A combinatorial design is an arrangement of a set of elements into defined substructures ({\it blocks}) in such a way that desired combinatorial properties are satisfied. They have many applications in, for example, communications, cryptography and networking \cite{colbournDinitzStinson_applications_1999}, optical orthogonal codes \cite{djordjevic2003}, cancer trial design \cite{baileyCameron2019} and genetic screening algorithms \cite{huber_genetic_2013}. 

In this paper we are mainly interested in designs for which every pair of elements appears in exactly one block together and one of the desirable properties is {\it resolvability}, i.e., where the blocks can be arranged into rounds in which every element appears once. 
Such designs are called resolvable designs \cite{furinoMiaoYin_book_1996}.

In recent work \cite{mibakavapu2021}, we highlighted the link between resolvable designs and another mathematical problem: the Social Golfer problem SGP. The Social Golfer problem is to allocate $v$ players into groups of equal size over as many rounds as possible. No two players should play together in more than one round, and the goal is to allocate as many rounds as possible. We note that the terms used in SGP (such as {\it group}) conflict with similar terms used in combinatorial designs. In \cref{sect:SGP} we restate the problem using combinatorial design terminology.

We were initially drawn to the problem when trying to allocate participants in large online calls (i.e. on Zoom or Teams) to successive breakout rooms during the COVID-19 pandemic. We soon realised that the allocation problem had more general application in a teaching context (when allocating students to successive study groups for example). Allocating as many rounds of groups as possible, whilst ensuring that no two students appear in the same group together more than once directly maps to the Social Golfer problem. Indeed these allocations are examples of a class of combinatorial designs known as resolvable maximum packings \cite{geLamLingShen_ResolvableMaximalPackingsQuads_2005,stinsonWeiYin_Packings_2007}.

In \cite{mibakavapu2021} we also considered the case where there are two group sizes which differ by one, and in each round the number of groups of each size is constant. We called this the Social Golfer problem with adjacent group sizes (SGA). This is very helpful in situations where the total number of participants does not have any suitable divisors.

We provided an initial list of the best solutions we could find for (student) class sizes of up to $50$ and group sizes either a divisor of the class size or a pair of consecutive values in $\{(4,5),(5,6)\}$. In most cases we found our solutions by adapting resolvable designs, in particular {\it resolvable balanced incomplete block designs} and a particular type of {\it resolvable group divisible design}, namely  {\it resolvable transversal designs}. In some cases we  used (or adapted) existing published best solutions for the SGP. We also made all our solutions available on our web application BoRAT, the link for which can be found at \cite{miller_research}.  We continue to add solutions as they are generated.

In this paper we re-introduce the SGP and the SGA, and present examples of the combinatorial designs that we use to find good solutions to them. We both improve on solutions given in \cite{mibakavapu2021} and extend our results to up to $150$ participants. We also give complete solutions for the SGP with group sizes $3$ or $4$ for any (suitable) number of participants. We explain the methods used in \cite{mibakavapu2021} and introduce additional methods based on other resolvable, or partially resolvable designs such as  {\it resolvable group divisible designs}, {\it uniformally resolvable designs} and {\it incomplete resolvable transversal designs}. We describe some of the constructions used and highlight some difficult cases - presenting the best current solution to solve them. 

We present an algorithm to determine which method to use to find an optimal solution when groups are of equal size. If groups have size, $3$ or $4$ the algorithm returns a solution for any number of participants. Otherwise, the algorithm suggests a solution. We discuss how to construct solutions for unequal group sizes ($4$ and $5$, or $5$ and $6$) by adding participants to smaller solutions, or removing from them from larger solutions. 

Tables containing SGP and SGA solutions for up to $150$ participants are provided in \cref{appendix:tables}. 
%when groups have equal size and a suitable resolvable design is known, otherwise an optimal solution for either equal group size or group sizes $4$ and $5$, or $5$ and $6$ is returned.
%provided the number of participants is at most $200$.

%We make regular reference to the Handbook of Combinatorial Designs \cite{codihandbook}. To aid the reader we identify the specific chapter, chapter section, or item, using references of the form: \cite[Chapter~III]%{codihandbook}, \cite[II.2.8]{codihandbook} and \cite[Table~VI.16.54]{codihandbook}  respectively. 

%this is taken from the original paper - just edited slightly
%removed this because it was here to satisfy Symmetry Journal, no longer required
\begin{comment} 
Social Golfer type problems are inherently linked to symmetry. Combinatorial search for solutions is hard due to the number of equivalent partial solutions at each level of search. As we discuss in \cref{sect:SGP}, there are four types of symmetries present for any SGP instance. Symmetry breaking techniques \cite{smith01,fomi2001,petrie2004,domical2004,gekelimcmism2005,gent2006handbook} must be used to eliminate these equivalent structures and reduce the search space to a tractable size. In addition, the constructions we refer to in  \cref{sect:maxSGPSolutions} to generate designs for maximal SGP solutions, often involve finite (Galois) groups to create parallel classes of blocks from a set of base blocks. The nature of finite groups prevents the repetition of pairs in different blocks by avoiding undesired periodicity due to symmetry.
\end{comment}

\section{Social Golfer problem and the Social Golfer problem with adjacent block sizes}\label{sect:SGP}
\subsection{The Social Golfer Problem}\label{subsect:SGP} The Social Golfer problem (which we refer to as SGP) \cite{csplib:prob010,mathGames,TriskaThesis2008} originated from the following question posed in 1998 to \texttt{sci.op-research}: 
\begin{quote}
    $32$ golfers play golf once a week, and always in groups of $4$. For how many weeks can they play such that no two players play together more than once in the same group?
\end{quote}

In fact the answer to this problem is ten weeks. We will show in \cref{SGP:k4} that the solution can be obtained using a combinatorial structure known as a {\it Resolvable Group Divisible Design}. The design for this case was originally constructed by Shen in \cite{shen_RGDD4_JShanghai_1996}
and independently solved by Colbourn in 1999 \cite{colbourn_Steiner2Design_1998}. In 2004, Aguado recognised that the design would solve the problem and published a solution \cite{aguadoSocialGolfer2004}.

The problem can be generalised as follows (see \cite{liuLofflerHofstedt_SG_2019}):

\begin{quote}
The SGP consists of scheduling $n=g*s$ players into $g$ groups of $s$ players for $w$ weeks so that any two players are assigned to the same group at most once in $w$ weeks. 
\end{quote}

Applications of the SGP include group teaching scenarios \cite{mibakavapu2021,CaoPngCaiCenXu2021,limpanuparb_datta_2021}, voting strategies \cite{kanav2024}, and scheduling for sporting competitions \cite{lambers_rothuizen_spieksma_2021, lambers_rothuizen_spieksma_2023, lester2021}. 

%in terms of weeks, groups etc.
Any instance of the SGP exhibits four types of symmetry arising from permuting the weeks, groups within weeks, players within groups and the players within the overall set of players. 

The SGP has been adopted by the Combinatorial Search community as a benchmark for testing advanced search and symmetry breaking  techniques
\cite{fomi2001,gekelimcmism2005,gent2006handbook,babr2002,domical2004,petrie2004,smith01}. Techniques for solving the Social Golfer problem include heuristic approaches \cite{dohe2005,schmand2022}, answer-set programming \cite{eiter2022}, Boolean satisfiability (SAT) encoding \cite{gentlynce2005,jaju2023symmetry,trimu2012a,triska2012} and constraint-based techniques \cite{lawlee2004,liuLofflerHofstedt_SG_2019,smith01}. Excellent surveys can be found in \cite{TriskaThesis2008} and more
recently in \cite{liuLofflerHofstedt_SG_revisited_2019,trimu2012a, trimu2012b}.  In this paper, rather than relying on optimized search, we use existing combinatorial designs to construct solutions. 

We will adopt a notation that is slightly different from the one used in the SGP literature. As we will be using combinatorial designs to generate rounds, we refer to {\it blocks} rather than groups (the term group already plays an important role in designs, as we will see in \cref{sect:combinatorialStructures}). We also refer to {\it rounds} rather than weeks and {\it points} instead of players. We will continue to use the terms {\it players} and {\it groups} until we redefine the problem in \cref{subsect:SGA} (\cref{defn:sgp}).

We focus on {\it optimal} solutions: i.e. solutions to problems for which no solution with more rounds is known. We hence distinguish between {\it optimal} and {\it maximal} solutions:
\begin{definition}\label{defn:SGP}
    A solution to SGP$(v,k)$ with $r$ rounds is said to be an optimal solution if no solution with more than $r$ rounds is known.  A solution is maximal if no solution with more than $r$ rounds is possible. 
\end{definition}

Clearly, if a solution is maximal, it must be optimal. Since a player can be in a group with every other player at most once, and every group that contains a given player contains $s-1$ other players, a player can be in at most $r=\lfloor(n-1)/(s-1)\rfloor$ rounds. Hence, if a solution has $r$ rounds, it is clearly maximal. 

Consider the example shown in  \cref{fig:SGPExample1}. There are $28$ players (labelled $0$ to $27$), and each round consists of $7$ groups of size $4$. There are $9$ rounds and the solution is clearly maximal. 

\begin{figure*}
\begin{equation*}
\boxed{
\begin{aligned}
&\bf{Round}\;0:\\
&[0,9,18,27],[1,2,14,16],[3,6,13,17],[4,8,21,24],[5,7,19,20],[10,11,23,25],\\
&[12,15,22,26]\\
% \\
&\bf{Round}\;1:\\
&[0,2,12,17],[1,10,19,27],[3,8,18,20],[4,7,14,15],[5,6,22,25],[9,11,21,26],\\
&[13,16,23,24]\\
% \\
&\bf{Round}\;2:\\
&[0,6,11,16],[1,8,22,23],[2,7,18,24],[3,12,21,27],[4,5,10,17],[9,15,20,25],\\
&[13,14,19,26]\\
% \\
&\bf{Round}\;3:\\
&[0,5,24,26],[1,4,11,12],[2,3,19,22],[6,8,9,14],[7,16,25,27],[10,13,20,21],\\
&[15,17,18,23]\\
% \\
&\bf{Round}\;4:\\
&[0,4,20,23],[1,3,24,25],[2,5,9,13],[6,7,10,12],[8,17,26,27],[11,14,18,22],\\
&[15,16,19,21]\\
% \\
&\bf{Round}\;5:\\
&[0,1,13,15],[2,11,20,27],[3,7,23,26],[4,6,18,19],[5,8,12,16],[9,10,22,24],\\
&[14,17,21,25]\\
% \\
&\bf{Round}\;6:\\
&[0,3,10,14],[1,5,18,21],[2,4,25,26],[6,15,24,27],[7,8,11,13],[9,12,19,23],\\
&[16,17,20,22]\\
% \\
&\bf{Round}\;7:\\
&[0,7,21,22],[1,6,20,26],[2,8,10,15],[3,4,9,16],[5,14,23,27],[11,17,19,24],\\
&[12,13,18,25]\\
% \\
&\bf{Round}\;8:\\
&[0,8,19,25],[1,7,9,17],[2,6,21,23],[3,5,11,15],[4,13,22,27],[10,16,18,26],\\
&[12,14,20,24]
\end{aligned}
}
\end{equation*}

\caption{Maximal Social Golfer solution for v=28 and groups of size 4. There are $9$ rounds, each containing $7$ groups. Every pair of players appears exactly once.
\label{fig:SGPExample1}
}
\end{figure*}

Now consider the example shown in \cref{fig:SGPExample2}. There are $24$ players and each round consists of $6$ groups of size $4$. There are $7$ rounds and the solution is again maximal. 

\begin{figure}
\begin{equation*}
\boxed{
\begin{aligned}
&\bf{Round}\;0:\\
&[0,1,3,14],[2,7,12,21],[4,13,18,23],[5,10,15,20],[6,16,17,19],[8,9,11,22]\\
\\
&\bf{Round}\;1:\\
&[0,17,18,20],[1,2,4,8],[3,7,13,22],[5,14,19,23],[6,11,15,21],[9,10,12,16]\\
\\
&\bf{Round}\;2:\\
&[0,12,15,22],[1,18,19,21],[2,3,5,9],[4,7,14,16],[6,8,20,23],[10,11,13,17]\\
\\
&\bf{Round}\;3:\\
&[0,9,21,23],[1,13,15,16],[2,19,20,22],[3,4,6,10],[5,7,8,17],[11,12,14,18]\\
\\
&\bf{Round}\;4:\\
&[0,4,5,11],[1,10,22,23],[2,14,15,17],[3,16,20,21],[6,7,9,18],[8,12,13,19]\\
\\
&\bf{Round}\;5:\\
&[0,7,10,19],[1,5,6,12],[2,11,16,23],[3,8,15,18],[4,17,21,22],[9,13,14,20]\\
\\
&\bf{Round}\;6:\\
&[0,2,6,13],[1,7,11,20],[3,12,17,23],[4,9,15,19],[5,16,18,22],[8,10,14,21]
\end{aligned}
}
\end{equation*}
\caption{Maximal Social Golfer solution for v=24 and groups of size 4. There are $7$ rounds, each containing $6$ groups. Not all pairs of players appear.
\label{fig:SGPExample2}
}
\end{figure}

Note that, in the example shown in  \cref{fig:SGPExample1}, all pairs of players appear in the groups. However, in \cref{fig:SGPExample2}, this is not the case (for example, player $0$ is not in a group with either player $8$ or player $16$). As we will see in \cref{sect:constructions}, this is because the example in \cref{fig:SGPExample1} is from a \emph{Resolvable Balanced Incomplete Block Design} (RBIBD), and the example in  \cref{fig:SGPExample2} is from a \emph{Resolvable Group Divisible Design} (RGDD).

\subsection{The Social Golfer Problem with adjacent block sizes}\label{subsect:SGA} 
In \cite{mibakavapu2021} we introduced a problem closely related to SGP, which we call the social golfer problem with adjacent block sizes (SGA). This problem arises when allocating players to groups when the group size doesn't exactly divide the number of players, in particular when the number of players is prime. In this case we don't insist that all groups have the same sizes. Obviously, if all group sizes are possible, there would be far too many solutions. A reasonable first restriction is to have group sizes that differ by as little as possible, i.e., they differ by at most $1$. We denote the set of group sizes by $K$.

As we will show in \cref{sect:adjacentBlockSizes}, solutions can be obtained in these cases by examining a solution with more players and equal sized groups and removing certain sets of players (or, in some cases, by adding players to solutions with fewer players).

For example, if we remove players $0$ and $8$ from the example in \cref{fig:SGPExample2} we would obtain the allocation shown in \cref{fig:SGPExample3}. The groups now have sizes $3$ and $4$ (as we already knew that $0$ and $8$ were not in an existing group, it comes as no surprise that no group loses more than one player).
\begin{figure}
\begin{equation*}
\boxed{
\begin{aligned}
&\bf{Round}\;0:\\
&[1,3,14],[2,7,12,21],[4,13,18,23],[5,10,15,20],[6,16,17,19],[9,11,22]\\
\\
&\bf{Round}\;1:\\
&[17,18,20],[1,2,4],[3,7,13,22],[5,14,19,23],[6,11,15,21],[9,10,12,16]\\
\\
&\bf{Round}\;2:\\
&[12,15,22],[1,18,19,21],[2,3,5,9],[4,7,14,16],[6,20,23],[10,11,13,17]\\
\\
&\bf{Round}\;3:\\
&[9,21,23],[1,13,15,16],[2,19,20,22],[3,4,6,10],[5,7,17],[11,12,14,18]\\
\\
&\bf{Round}\;4:\\
&[4,5,11],[1,10,22,23],[2,14,15,17],[3,16,20,21],[6,7,9,18],[12,13,19]\\
\\
&\bf{Round}\;5:\\
&[7,10,19],[1,5,6,12],[2,11,16,23],[3,15,18],[4,17,21,22],[9,13,14,20]\\
\\
&\bf{Round}\;6:\\
&[2,6,13],[1,7,11,20],[3,12,17,23],[4,9,15,19],[5,16,18,22],[10,14,21]
\end{aligned}
}
\end{equation*}
\caption{Maximal Social Golfer solution for v=22 and groups of size 4 and 3, obtained by removing players $0$ and $8$ from \cref{fig:SGPExample3}. Note that players are not sequentially labelled. 
\label{fig:SGPExample3}
}
\end{figure}

Notice that the players are now from the set $\{1,2,\ldots,7,9,10,\ldots 23\}$ (i.e. they are not sequential). To obtain a legitimate allocation (with players in set 
$\{x:0\leq x \leq 21\}$ we would need to renumber the players in the obvious way, and by doing so we obtain the allocation shown in \cref{fig:SGPExample4}.
\begin{figure}
\begin{equation*}
\boxed{
\begin{aligned}
&\bf{Round}\;0:\\
&[0,2,12],[1,6,10,19],[3,11,16,21],[4,8,13,18],[5,14,15,17],[7,9,20]\\
\\
&\bf{Round}\;1:\\
&[15,16,18],[0,1,3],[2,6,11,20],[4,12,17,21],[5,9,13,19],[7,8,10,14]\\
\\
&\bf{Round}\;2:\\
&[10,13,20],[0,16,17,19],[1,2,4,7],[3,6,12,14],[5,18,21],[8,9,11,15]\\
\\
&\bf{Round}\;3:\\
&[7,19,21],[0,11,13,14],[1,17,18,20],[2,3,5,8],[4,6,15],[9,10,12,16]\\
\\
&\bf{Round}\;4:\\
&[3,4,9],[0,8,20,21],[1,12,13,15],[2,14,18,19],[5,6,7,16],[10,11,17]\\
\\
&\bf{Round}\;5:\\
&[6,8,17],[0,4,5,10],[1,9,14,21],[2,13,16],[3,15,19,20],[7,11,12,18]\\
\\
&\bf{Round}\;6:\\
&[1,5,11],[0,6,9,18],[2,10,15,21],[3,7,13,17],[4,14,16,20],[8,12,19]
\end{aligned}
}
\end{equation*}
\caption{Maximal Social Golfer solution for v=22 and groups of size 4 and 3, obtained by removing players $0$ and $8$ from \cref{fig:SGPExample3}, and with appropriate renaming so that players are sequentially labelled. 
\label{fig:SGPExample4}
}
\end{figure}

As the problem initially arose in the context of assigning students to small groups, for which group sizes between $4$ and $6$ are considered optimal \cite{Barr_Nabi_Somerville_2020}, in \cite{mibakavapu2021} we further restricted ourselves to the cases $K=\{4,5\}$ and $K=\{5,6\}$. We will continue to focus on these cases in this paper, although the techniques we introduce can easily be extended to further block sizes.

The link between the Social Golfer problem and combinatorial designs, specifically mutually 
orthogonal Latin squares 
\cite[Table~VI.16.54]{codihandbook}, Kirkman systems \cite[II.2.8]{codihandbook}, \cite{kirkmanKTS-1847,ray_chaudhuriWilson_RBIBD3_1971} and Resolvable Group Divisible Designs \cite[IV.5.3]{codihandbook} is well known 
\cite{babr2002,harveyWinterer_SGP_MOLR_2005,TriskaThesis2008}. 
Some optimal solutions to the Social Golfer problem can be found at \cite{mathGames} and \cite{csplib:prob010_results}.
\begin{comment}both for allocation schedules where there is only one group size, and to construct allocation schedules with two different group sizes by adding or removing points to/from existing Social Golfer solutions. 
\end{comment}

Because we will derive our solutions to SGP and SGA from suitable combinatorial designs, to avoid confusion we will henceforth adopt  \cref{defn:sgp,defn:sga}  for SGP and SGA, in terms of \emph{points} and \emph{blocks} (rather than players and groups):

\begin{definition}\label{defn:sgp}
The SGP$(v,k,r)$ consists of arranging $v=k*n$ points into $r$ rounds of $n$ non-intersecting blocks of size $k$ in such a way that two points are assigned to the same block at most once.
\end{definition}

\begin{definition}\label{defn:sga}
The SGA$(v,n_{1},k_{1},n_{2},k_{2},r)$ consists of arranging $v=n_{1}*k_{1} + n_{2}*k_{2}$ points into $r$ rounds of $n_{1}+n_{2}$ non-intersecting blocks in such as way that:
\begin{enumerate}
    \item $k_2=k_{1}+1$,
    \item each round contains $n_{1}$ blocks of size $k_{1}$ and $n_{2}$ blocks of size $k_{2}$,
    \item two points are assigned to the same block at most once.
\end{enumerate}
\end{definition}

Generally the {\it Social Golfer Problem} refers to fixed values of $v$ and $k$ and finding the largest $r$ for which SGP$(v,k,r)$ has a solution. In our context of allocation schedules, we may not require a maximal solution (as often only a few rounds are required, or the maximal solution would provide too many rounds to feasibly implement). However, finding the maximal solution (or the optimal solution) will allow any number of rounds $r^{\prime}$ up to this value to be obtained by choosing only the first $r^{\prime}$ rounds:

\begin{lemma}
If there is a solution for SG$(v,k,r)$ then there is a solution for SG$(v,k,r^{\prime})$, for any $r^{\prime}$ less than $r$.
\end{lemma}

\section{Combinatorial structures - PBDs, BIBDs, GDDs and URDs}\label{sect:combinatorialStructures}
In this section we introduce the combinatorial structures that will be required in the remainder of the paper. In order to maintain consistency, we adapt the definitions given in \cite{codihandbook} where possible, although they appear in slightly different form in various publications. 
A combinatorial design is an arrangement of a set of elements into defined substructures ({\it blocks}) in such a way that desired combinatorial properties are satisfied. They have many applications in, for example communications, cryptography and networking \cite{colbournDinitzStinson_applications_1999}, optical orthogonal codes \cite{djordjevic2003}, cancer trial design \cite{baileyCameron2019} and genetic screening algorithms \cite{huber_genetic_2013}. Our starting point is pairwise balanced designs. 

\begin{definition}\label{defn:pbd}
Let $K$ be a subset of positive integers and let $\lambda$ be a positive integer. A {\it pairwise balanced design} (PBD) of order $v$ with block sizes from $K$ is a pair $(V,B)$, where $V$ is a finite set (the {\it point set}) of cardinality $v$ and $B$ is a family of subsets ({\it blocks}) of $V$ that satisfy (1) if $b\in B$ then $|b|\in K$ and (2) every pair of distinct elements of $V$ occurs in exactly $\lambda$ blocks of $B$. When $|V|=v$ and $\lambda=1$ we refer to such a PBD as a PBD$(v,K)$. 
\end{definition}

In this paper, we only consider PBDs for which $\lambda=1$, so in the remaining definitions we will only include this case.  The following definition  is for a BIBD, which  is a  special type of PBD with $|K|=1$.

\begin{definition}\label{defn:bibd}
A {\it balanced incomplete block design} (BIBD) is a PBD where all blocks have the same size. If the block size is $k$ and the design is on $v$ points (i.e. has order $v$) we refer to a ($k$-BIBD) and to a specific BIBD with parameters $v$ and $k$ as a BIBD$(v,k)$. If the blocks can be arranged into classes (or rounds)  in which every class contains every point exactly once (i.e. {\it parallel classes}), then the BIBD is said to be resolvable, and is called a {resolvable balanced incomplete block design} (referred to as an RBIBD, a $k$-RBIBD or as an RBIBD$(v,k)$).
\end{definition}

We give some examples of BIBDs below. In all cases we present blocks horizontally, a convention we will adopt throughout this paper. 

\begin{example} \label{example1:introductory}
Suppose that $v=7$ and $k=3$, the following is a BIBD that is not resolvable. In fact, all blocks intersect and the design is an example of a {\it finite projective plane} \cite{battenBeutelspackerLinerSpaces1993}.
$$(0,1,2),\;(0,3,4),\;(0,5,6),\;(1,3,5),\;(1,4,6),\;(2,3,6) $$
\end{example}

\begin{example} The example shown in  \cref{fig:rbibd} is resolvable. To demonstrate this we arrange the blocks into rounds (unlike the SGP and SGA examples seen previously, we no longer label the rounds).
\end{example}
\begin{figure}
\begin{equation*}
\boxed{
\begin{aligned}
&[0,4,5,6],[1,7,11,15],[2,8,12,13],[3,9,10,14]\\
\\
&[0,1,2,3],[4,7,10,13],[5,8,11,14],[6,9,12,15]\\
\\
&[0,7,8,9],[1,4,12,14],[2,5,10,15],[3,6,11,13]\\
\\
&[0,10,11,12],[1,5,9,13],[2,6,7,14],[3,4,8,15]\\
\\
&[0,13,14,15],[1,6,8,10],[2,4,9,11],[3,5,7,12]
\end{aligned}
}
\end{equation*}
\caption{An RBIBD$(16,4)$. 
\label{fig:rbibd}
}
\end{figure}
A $3$-RBIBD is known as a Kirkman triple system \cite{kirkmanKTS-1847,ray_chaudhuriWilson_RBIBD3_1971,stinson-KTS-survey-1991,reesWallis2002}. 

The following lemma combines some trivial conditions for BIBDs, which can be found in 
\cite[II.1.1]{codihandbook}. The parameters $v$ and $k$ are as defined above, and  $b$ and $r$ denote the number of blocks and the number of blocks containing any given point (the {\it regularity}) respectively. 

\begin{lemma}\label{lemmma:bibdBasic}
If $\mathcal{B}$ is a BIBD$(v,k)$ then the following hold:
\begin{itemize}
\item $r=\frac{v-1}{k-1}$
\item $b=\frac{v(v-1)}{k(k-1)}=\frac{vr}{k}$
\item If $\mathcal{B}$ is resolvable, there are $r$ rounds, each of which contain $v/k$ blocks.
\end{itemize}
\end{lemma}

\begin{comment}
In Theorem \ref{thm:abelsResult} we present a construction for a BIBD$(2479,7)$. The construction relies on some definitions that we have yet to introduce. 
\end{comment}

\cref{lem:bibdNec}, stated in \cite[Theorem~1.1]{hananiBIBDRelated_1975}, gives a necessary condition for the existence of a BIBD$(v,k)$. The sufficiency of this condition, for $3\leq k\leq 7$ is considered in \cref{lemma:bibdExistence}.

\begin{lemma}\label{lem:bibdNec}
A BIBD$(v,k)$ exists only if $v\equiv 1 \pmod{k-1}$ and $v\equiv 1 \pmod{k-1}$.
\end{lemma}

\begin{lemma}\label{lemma:bibdExistence}
If $3\leq k \leq 7$ then a BIBD$(v,k)$ exists if:
\begin{itemize}
\item $k\leq 5$ and $v\equiv 1$ or $k \pmod{k(k-1)}$ \cite[Lemmas~5.4,~5.11~and~5.19]{hananiBIBDRelated_1975},
\item $k=6$,  $v\equiv 1$ or $6 \pmod{15}$ and $v\not\in \{16,21,36,46\}$ or one of $29$ other possible exceptions \cite{hananiBIBDRelated_1975, abelmillsBIBD6_1995, abelGreig_RBIBD5_BIBD6_1997, houghtenNoBIBD_46_6_2001, abelBluskovGreig2007},
\item $k=7$, $v\equiv 1$ or $7 \pmod{42}$ and $v$ not equal to one of $21$  exceptions (\cite{abel_BIBDk6Tok10_1996,abelGreig_BIBD7_1998,jankoTonchev_BIBD7_1998,abel_BIBD7_2000}.
%and Theorem \ref{thm:abelsResult}).
\end{itemize}
\end{lemma}

Surveys of existence results for BIBDs with block size at most $9$ are given in \cite{abelGreigBIBD_2007} and \cite[Table~1]{adams_surveyGDesigns_2008}.

We focus here on RBIBDs as they are of most use to us in the context of solutions for the Social Golfer problem. 
We now give existence results for $k$-RBIBDs, for $3\leq k\leq 6$.

\begin{lemma}\label{lemma:rbibdExistence}
If $3\leq k \leq 5$ and an RBIBD$(v,k)$ exists then:
\begin{itemize}
\item $k=3$ and $v\equiv 3 \pmod{6}$ \rm\cite{ray_chaudhuriWilson_RBIBD3_1971},
\item $k=4$ and $v\equiv 4 \pmod{12}$ \cite{hananiRay_ChaudhuriWilson2006-orig1972},
\item $k=5$ and $v\equiv 5 \pmod{20}$, 
$v\not\in\{45,345,465,645\}$ \cite{chenZhu_RBIBD5_1987, zhuChenDu_RBIBD5_1987, zhuDuZhang_RBIBD5_1991, abelGreig_RBIBD5_BIBD6_1997, abelGeGreigZhu_RBIBD5_2001, abel_RBIBD5_2007}.
\end{itemize}
\end{lemma}
Note that for an RBIBD$(v,6)$ to exist, we require $v=30t+6$ for some integer $t\geq 0$. Many of the constructions for RBIBDs for $6\leq k\leq 9$ rely on a theorem using difference families (see \cref{defn:diff_fam}) due to Ray-Chaudhuri and Wilson \cite{ray_chaudhuriWilson_RBIBD3_1971}, which we state as \cref{thm:RayChouWilson}. 

\begin{definition}\label{defn:diff_fam}
A $(v,k,\lambda)$ difference family is a set of subsets $B=\{B_{1}, B_{2}, . . . , B_{t}\}$ of $G$, such that the order of $G$ is $v$, each subset $B_i$ has size $k$, and the differences $(a-b\mid a,b \in B_i; a \neq b, i=1 \dots t)$ contain each non-zero member of $G$ exactly $\lambda$ times.
\end{definition}

\begin{Theorem}\label{thm:RayChouWilson}
If $q$ is an odd prime power, and there is a
$(v, q, 1)$ difference family $\{B_{1}, B_{2}, . . . , B_{t}\}$ in $GF (q)$ such that the base
blocks are mutually disjoint, then an RBIBD$(kq,k)$ exists.
\end{Theorem}

Some results for the existence of RBIBD$(v,6)$ are given in \cref{theorem:k6_existence}
\begin{Theorem}\label{theorem:k6_existence}.

\noindent
\begin{enumerate}
\item An RBIBD$(36,6)$ does not exist 
\cite[Theorem~4.1.4]{furinoMiaoYin_book_1996}.
\item An RBIBD$(156,6)$ is equivalent to a $PG(3,5)$ \cite{lorimer_RBIBDFromPG3_1973}.
\item If $q$ is a prime power less than $2000$, not equal to $61$ or $121$, and $q\equiv 1 \pmod{30}$, then a  block-disjoint $(v, q, 1)$ difference family and thus an RBIBD$(6q,6)$ exists \cite[Theorem~11]{wilson_DifferenceFamilies_1972}, see also \cite[Table~VI.16.54]{codihandbook}.
\item RBIBD$(1716,6)$ exists \cite[Theorem~12.8]{greig_BIBDk=7-9_2001}.
\item If $q$ is a prime power, $q\equiv 1 \pmod{10}$, $q\not\equiv 1 \pmod{50}$, $q\neq 11$ and $q<1500$ then RBIBD$(6q,6)$ exists \cite[Lemma~3.3]{abel_BIBDk6Tok10_1996}. 
\item If $t \neq \{2,6,8,12\}$ and $t\leq 832$, and $6t+1$ is a prime power, and $t$ is even, then RBIBD$(30t+6,6)$ exists \cite[Theorem~5.1]{greig_GDD_constructions_1998}.
\item There exists an RBIBD$(125q+1,6,1)$ for any prime $q\equiv 7 \pmod{12}$ and $q>43$, and for any prime $q\equiv 1 \pmod{12}$ and $q>37$ \cite[Theorem~18]{costaFengWang_FDFs_2018}.
\end{enumerate}
\end{Theorem}
 Example sources for constructions of known RBIBD$(v,6)$ designs for $v=30t+6$, $t\leq 100$ are shown in \cref{table:RBIBD_6}. Parenthesised values of $t$ are those for which a construction is available in both 
\cite{abel_BIBDk6Tok10_1996} and \cite{greig_GDD_constructions_1998}.
Constructions for some of the values of $t$ shown in the table were previously given in \cite{wilson_DifferenceFamilies_1972}. Note that although there are very few RBIBD$(v,6)$ designs for $v=30t+6$ and $t$ odd, for $t\leq 100$, some larger such designs with $t$ odd can be constructed by noting that: if an RBIBD$(w,6)$ and an RTD$(6,w)$
(see \cref{defn:TD}) exist, then so does an RBIBD$(6w,6)$. For example, for $v=936$ ($t=31$) an RBIBD$(v,6)$ can be constructed this way. 

\begin{table}
\centering
%\begin{footnotesize}
% \begin{tabular}{|p{7.5cm}|p{4.5cm}|l|}
\begin{tabular}{|p{7cm}|l|}
\hline
t & Source\\
\hline
$0$ & single block \\ \hline
$1$ & does not exist \\ \hline
$5$ & \cite{lorimer_RBIBDFromPG3_1973} \\ \hline
$57$ & \cite{greig_BIBDk=7-9_2001} \\ \hline
odd and not in $\{1,5,57\}$, $2$, $22$, $34$, $44$, $50$, $64$, $74$, $78$, $80$, $82$
& unknown \\
\hline
$6$, $8$, $12$, $14$, $16$, $24$, $26$, $36$, $38$, $42$, $48$, $54$, $56$, $62$, $66$, $72$, $84$, $86$, $92$, $98$ & \cite{abel_BIBDk6Tok10_1996}\\ \hline
$4$, $10$, ($16$), $18$, $20$, ($26$), $28$, $30$, $32$, ($38$), $40$, $46$,  ($48$), $52$, ($56$), $58$, $60$, ($62$), ($66$), $68$, $70$, ($72$), $76$, $88$, $90$, $96$, $100$& \cite{greig_GDD_constructions_1998}\\
\hline
\end{tabular}
\caption{Sources for RBIBD$(v,6)$ constructions, $v=30t+6$, $0\leq t\leq 100$\label{table:RBIBD_6}
}
\end{table}

Constructions for $v=6q$ where $q$ is a prime power, $q\equiv 1 \pmod{10}$, $q\not\equiv 1 \pmod{50}$, $q\neq 11$, $q<1500$, are given in \cite{abel_BIBDk6Tok10_1996}, but for a small number of values of $q$ the parameter sets listed in appendix D in the paper appear to contain small errors. For $q$ a prime there are four such cases. These are: $q=41$, $q=71$, $q=431$ and $q=1321$. For $q=41$ change $z_{\infty}$ from $1$ to $15$ leaving all other values as presented in the paper. In all other cases leave $x$ and $z_{\infty}$ as presented in the paper and change $z[1][0]$, $z[1][1]$, $z[1][2]$, $z[2][1]$ and $z[2][2]$ as follows. For  $q=71$, $z[1][0]=2$, $z[1][1]=18$, $z[1][2]=16$, $z[2][1]=31$, $z[2][2]=25$. For $q=431$,  $z[1][0]=63$, $z[1][1]=3$, $z[1][2]=69$, $z[2][1]=65$, $z[2][2]=186$. For $q=1321$, $z[1][0]=404$, $z[1][1]=248$, $z[1][2]=260$, $z[2][1]=5$, $z[2][2]=84$.

For the constructions in \cite{greig_GDD_constructions_1998} (Theorem 5.1), we had to generate the values of some of the required parameters
%$\gamma_{ij}$, $0\leq i,j \leq 1$ 
ourselves, as they were not provided in the table of parameters (\cref{ownSG2}).

A summary of methods for RBIBD$(v,k)$, for $k\in\{7,8,9\}$ is given in \cite{greig_BIBDk=7-9_2001}. We summarise known results for these cases in \cref{theorem:rbibd7-9}.

\begin{Theorem}\label{theorem:rbibd7-9}
\noindent
\begin{enumerate}
\item Lists of $(q,k,1)$ block disjoint difference families for odd prime powers $q$, and $7\leq k \leq 9$ are given in \cite[Tables~C.I-C.III]{abel_BIBDk6Tok10_1996} and \cite[Tables~A.1-A.5]{greig_BIBDConstructions_ge_7_1990}.
\item A construction for a $(577,9,1)$ block disjoint difference family is given in \cite{buratti_powerfulMethodForDiffFamilies_1995}.
\item If $t$ is one of the following $25$ integers less than $100$, $RB(42t+7,7)$ exists (\cite[Table~C.1]{greig_BIBDk=7-9_2001})
$0,1,8,9,17,28,33,41,49,56,57,63,64,65,70,72,73,77,80,81,88,89,91,96,97$.
\item Constructions for $RB(56t+8,8)$ for all but $66$ values of $t$ are given in \cite{greigAbel_RBIBD8_1997} and constructions for $8$ further values of $t$ are given in \cite{costaFengWang_FDFs_2018}. 
\item If $t$ is {\bf not} one of the following values, and $t\leq 100$, then $RB(56t+8,8)$ exists:
$3,13,20,22,23,25,26,27,31,38,43,47,52,58,59,61,67,69,76,79,93$.
\item If $t$ is one of the following $16$ integers less than $100$, $RB(72t+9,9)$ exists (\cite[Table~C.1]{greig_BIBDk=7-9_2001})
$0,1,7,8,9,10,37,54,64,71,72,73,81,82,90,91$.
\end{enumerate}

\end{Theorem}

\begin{comment}
Additional methods for $k=8$ can be found in  
\cite{abel_BIBDk6Tok10_1996}, \cite{greigAbel_RBIBD8_1997}, \cite{greig_GDD_constructions_1998}, and \cite{costaFengWang_FDFs_2018}.
\end{comment}

A PBD that is particularly useful to us is one in which a distinguished set of blocks, called {\it groups} form a parallel class. This is known as a group divisible design (GDD). We only consider the case where all blocks have the same size (a $k$-GDD) and in most cases where the groups have the same size (uniform $k$-GDDs). More general group divisible designs are used in the constructions for GDDs whose blocks are resolvable  - namely resolvable group divisible designs (RGDDs). We provide the necessary definition below.

\begin{definition}\label{defn:gdd}
 A {\it group divisible design} is a triple $(V,B,G)$, where $V$ is a set of points, and $B$ and $G$ sets of subsets called {\it blocks} and {\it groups} respectively and where: (1)  every pair of distinct elements of $V$ occurs in exactly one block or group, but not both, and (2) the groups partition the points. A $k$-GDD of type $g_{0}^{u_{0}},g_{1}^{u_{1}},\ldots g_{s-1}^{u_{s-1}}$ is a GDD with blocks of size $k$ and $u_{i}$ groups of size $g_{i}$, for $0\leq i \leq s-1$. A {\it uniform} $k$-GDD of type $g^{u}$ is a $k$-GDD with $u$ groups of size $g$. If the blocks can be arranged into parallel classes, then the GDD is called a {\it resolvable group divisible design} (a $k\!-\!\!$-RGDD).  If the blocks can be 
 arranged into partial parallel classes, each of which misses just the points in a single group, then the GDD is called a frame.

\end{definition}
Since, in an RGDD, every point appears in a block with every other point  that does not lie in the same group as it, 
a $k\!-\!\!$-RGDD with groups of size $g$ has $(v-g)/(k-1)$ rounds. 
A $4\!-\!\!$-RGDD of type $7^{4}$ is shown in \cref{fig:RGDDExample}. Note that we present the blocks horizontally, and the groups vertically. 
\begin{figure}
\begin{equation*}
\boxed{
\begin{aligned}
&Blocks:\\
&[0,7,14,21],[1,8,16,24],[2,9,18,27],[3,10,20,23],[4,11,15,26],[5,12,17,22],%\\
[6,13,19,25]\\
\\
&[0,8,15,22],[1,9,17,25],[2,10,19,21],[3,11,14,24],[4,12,16,27],[5,13,18,23],%\\
[6,7,20,26]\\
\\
&[0,9,16,23],[1,10,18,26],[2,11,20,22],[3,12,15,25],[4,13,17,21],[5,7,19,24],%\\
[6,8,14,27]\\
\\
&[0,10,17,24],[1,11,19,27],[2,12,14,23],[3,13,16,26],[4,7,18,22],[5,8,20,25],%\\
[6,9,15,21]\\
\\
&[0,11,18,25],[1,12,20,21],[2,13,15,24],[3,7,17,27],[4,8,19,23],[5,9,14,26],%\\
[6,10,16,22]\\
\\
&[0,12,19,26],[1,13,14,22],[2,7,16,25],[3,8,18,21],[4,9,20,24],[5,10,15,27],%\\
[6,11,17,23]\\
\\
&[0,13,20,27],[1,7,15,23],[2,8,17,26],[3,9,19,22],[4,10,14,25],[5,11,16,21],%\\
[6,12,18,24]\\
\\
&Groups:\\
% &[0,1,2,3,4,5,6],[7,8,9,10,11,12,13],[14,15,16,17,18,19,20],[21,22,23,24,25,26,27]\\
&\begin{bNiceMatrix}
0 \\ 1 \\ 2 \\ 3 \\ 4 \\ 5 \\ 6 \\
\end{bNiceMatrix}
\begin{bNiceMatrix}
7 \\ 8 \\ 9 \\ 10 \\ 11 \\ 12 \\ 13 \\
\end{bNiceMatrix}
\begin{bNiceMatrix}
14 \\ 15 \\ 16 \\ 17 \\ 18 \\ 19 \\ 20 \\
\end{bNiceMatrix}
\begin{bNiceMatrix}
21 \\ 22 \\ 23 \\ 24 \\ 25 \\ 26 \\ 27 \\
\end{bNiceMatrix}
\end{aligned}
}
\end{equation*}
\caption{$4\!-\!$-RGDD of type $7^{4}$.
\label{fig:RGDDExample}
}
\end{figure}

Notice that a $k$-RBIBD is a $k\!-\!\!$-RGDD of type $k^{v/k}$ (the last round of blocks provide the groups).

\cref{lemma:RGDDUniform}, from \cite[Lemma~1.1]{vanstone_doublyResolvable_1980}, implies that if we are trying to use the blocks of an RGDD to generate rounds of blocks of equal size, we must use a uniform RGDD.

\begin{lemma}\label{lemma:RGDDUniform}
If an RGDD has blocks of the same size, then its groups are of the same size.   
\end{lemma}

A $3\!-\!\!$-RGDD of type $2^{v/2}$ is known as a nearly Kirkman triple system (NKTS$(v)$) \cite{kotzigRosa_NKTS_1974,bakerWilson77,brouwer_NewNKTS_1978,reesStinson_RGDD3_1987,assha1989,shen_NKS_1990,reesWallis2002,abelChanColbournLamken_DRNKTS_2013}. The existence result in \cref{lemma:nktsExistence} is from \cite[Theorem~27]{reesWallis2002}.

\begin{lemma}\label{lemma:nktsExistence}
There exists an NKTS$(v)$ if and only if $v\equiv 0 \pmod{6}$ and $v \notin \{6,12\}$.
\end{lemma}

When $k=4$ and $v\equiv 0$ or $8 \pmod{12}$, we will use $4\!-\!\!$-RGDDs with group size $3$ or $2$ respectively to produce maximal allocations. 
%Moving comment to the todos

The existence result in \cref{lem:rggd4_1} for $4\!-\!\!$-RGDDs with groups of size $3$ follows from constructions in 
\cite{shen_RGDD4_JCMCC_1987,shen_NKS_1990,reesStinson_Frames4_1992,shen_RGDD4_JShanghai_1996,shenShen_RGDD4_2002,ge_RGDD4_2002,geLingRGDD4_Survey_2004,geLamLingShen_ResolvableMaximalPackingsQuads_2005}.
\begin{lemma}\label{lem:rggd4_1}
If $v\equiv 0 \pmod{12}$ and $u=v/3$ then a $4\!-\!\!$-RGDD of type $3^{u}$ exists, except when $u=4$.
\end{lemma}

The existence result in \cref{lem:rggd4_2} for $4\!-\!\!$-RGDDs with groups of size $2$ follows from constructions in 
\cite{shen_RGDD4_JShanghai_1996,kreherLingReesLam_Note_4RGDDs_2003,
geLingRGDD4_Survey_2004,
geLamLingShen_ResolvableMaximalPackingsQuads_2005,
schusterGe_URD_3_4_2010,
weiGe_5GDD_4Frame_4RGDD_2014,duAbelWang_RGGD4_2015,weiGe_URD_2_4_2017}.
\begin{lemma}\label{lem:rggd4_2}
If $v\equiv 8 \pmod{12}$ and $u=v/2$ then a $4\!-\!\!$-RGDD of type $2^{u}$ exists, except when $u=4$, $u=10$ and possibly when $u\in\{46,70,82,94,100,118,130,202,214\}$.
\end{lemma}

\begin{figure}
    \centering
    \includegraphics[width=\textwidth]{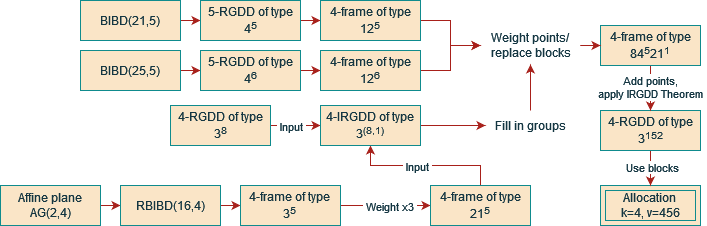}
    \caption{Steps involved in the construction of a $4-$RGDD of type $3^{152}$ and its use in the construction of an allocation with $v=456$, $k=4$}
    \label{fig:allocation-diagram}
\end{figure}

We use $4\!-\!\!$-RGDDs for some of our allocations in \cref{sect:maxSGPSolutions}. The constructions for them are often recursive and can involve the construction of many smaller, intermediate designs.  \cref{fig:allocation-diagram} illustrates this for an example, with $v=456$ and $g=3$. Even here we have omitted some of the necessary stages (for example constructing suitable sets of mutually orthogonal Latin squares (see \cref{defn:orthogonal}) and constructing the small input BIBDs).

The existence of many of the combinatorial structures we use to provide SGP solutions rely on the existence of sufficient numbers of {\it mutually orthogonal Latin squares} (MOLS). The following definitions are from 
\cite[Chapter~III]{codihandbook}:

\begin{definition}\label{defn:latinSquare}
    A Latin square $L$ of side $n$ (or order $n$) is an $n\times n$ array in which each cell $L(a,b)$ contains a single symbol from an $n$-set $S$, such that each symbol occurs exactly once in each row and each column.
\end{definition}

\begin{definition}\label{defn:orthogonal}
Two Latin squares $L$ and $L^{\prime}$ of the same order are {\it orthogonal} if $L(a,b)=L(c,d)$ and $L^{\prime}(a,b)=L^{\prime}(c,d)$ implies $a=c$ and $b=d$. A set of Latin squares $L_{1},L_{2},\ldots, L_{m}$ is mutually orthogonal (or is a set of $m$ MOLS), if whenever $1\leq i < j \leq m$, $L_{i}$ and $L_{j}$ are orthogonal. 
\end{definition}

For a given $n$, it is common to use the notation $N(n)$ to denote the maximum  number of MOLS of order $n$, and we say that a set of MOLS of order $n$ is optimal if no larger set of MOLS of order $n$ is known. Some useful results (originally from \cite{macNeishEulerSquares_1922}, and also stated in 
\cite[III.3.1]{codihandbook}) are given in  \cref{thm:macNeish}.

\begin{Theorem}\label{thm:macNeish} 

\noindent
\begin{enumerate}
    \item $N(n\times m)\geq min \{N(n),N(m)\}$.
    \item If $n = p^{e}$, for prime $p$ and some $e>0$, then $N(n)=n-1$.
    \item If $n=p_{1}^{e_{1}}p_{2}^{e_{2}}\ldots p_{k}^{e_{k}}$ where each $p_{i}$ is a prime, then $N(n)\geq min\{p_{i}^{e_{i}}-1|i \in \{1,2,\ldots, k\}\}$.
\end{enumerate}
\end{Theorem}

An original comprehensive list of the then greatest known lower bounds on the number of MOLS of side $n$, where $n\leq 10,000$ was given in 1979 in \cite{brouwer_MOLSTable_1979}. The table was updated in 1996 in \cite[Table~II.2.72]{codihandbook_old}, and again in 2007 in 
\cite[Table~III.3.87]{codihandbook}. 
%(specifically in the chapter on MOLS %\cite{abelColbournDinitz_MOLSTable_2007}). 

Recently the list was updated for $n\leq 500$ in \cite{millerAbelValkovFraser_MOLS_2024}. There a table listing the sources of the constructions for optimal sets of MOLS was provided for each $n$ and a few corrections to constructions were given. 
In \cref{table:molsTable} we provide similar information, but for $n$ up to $100$ only. By convention, N$(n)$ for $n \in \{0,1\}$ is said to be $\infty$, so these cases are omitted.
In addition, in the light of \cref{thm:macNeish}, we do not consider cases where $n$ is a prime power. 

\begin{table*}
\centering
%\begin{footnotesize}
% \begin{tabular}{|p{7.5cm}|p{4.5cm}|l|}
\begin{tabular}{|l|l|l|}
\hline
n & $N(n)$ (lower bound) & Source\\
\hline
$6$, $10$, $15$, $20$, $21$, $22$, $24$, $26$, $28$, $30$&
$1$, $2$, $4$, $4$, $5$, $3$, $7$, $4$, $5$, $4$&
\cite{abelColbournDinitz_MOLSTable_2007}\\
$33$, $34$, $36$, $38$, $39$, $40$, $42$, $44$, $46$&
$5$, $4$, $8$, $4$, $5$, $7$, $5$, $5$, $4$&\\
$50$, $51$, $52$, $55$, $56$, $62$, $75$, $80$&
$6$, $5$, $5$, $6$, $7$, $5$, $7$, $9$&\\
\hline
$12$&$5$&\cite{robertsMOLS_1992}\\
\hline
$14$&$4$&\cite{Todorov2012}\\
\hline
$18$, $60$&$5$, $5$&\cite{abelMOLS_18_60_2015}
\\
\hline
$35$, $48$, $63$&$6$, $10$, $8$&
\cite{janiszczakStaszewskiCodesLS_2019}\\
\hline
$45$, $54$, $96$& $6$, $8$, $10$& \cite{abelJaniszczakStaszewski_improvedMOLS_2024}\\
\hline
$57$, $69$, $70$, $74$, $78$, $84$, $90$&$7$, $6$, $6$, $5$, $6$, $6$, $6$&\cite{brouwerRecursiveMols_1991}\\
\hline
$58$, $66$, $68$& all $5$&\cite{hananiBIBDRelated_1975}\\
\hline
$65$ & 7 & \cite{boseShrikhandeparkerMOLS_Further_1960}\\
\hline
 $72$, $77$, $88$, $99$& $7$, $6$, $7$, $8$
& \cref{thm:macNeish}\\
\hline
$76$&$6$&\cite{colbournYinZhu_MOLS_76_1995}\\
\hline
$82$, $100$&
both $8$&\cite{colbourn_ITDConstructions_1996}\\
\hline
$85$, $86$, $87$, $92$, $93$, $94$, $95$, $98$& all $6$ &\cite{hananiBIBDRelated_1975}\\
\hline
$91$& $7$&\cite{brouwer_MOLSTable_1979}\\
\hline
\end{tabular}
%\end{footnotesize}
\caption{Sources for MOLS constructions for $2\leq n\leq 100$ where $n$ is not a prime power.} \label{table:molsTable}
\end{table*}

Further constructions (for $n>100$) can be found in \cite{boseShrikhandeFalsity_1960,brouwer_MOLSTable_1979, brouwer_separable_1980, 
colbournDinitzWojtas_thwarts_1995,abel_thesis_1995, 
colbournDinitzStinson_moreThwarts_1996,greigDesignsFromProjectivePlanes_1999,abelColbournWojtas_7and8MOLS_2004,millerAbelValkovFraser_MOLS_2024}.

\begin{comment}
Where to find constructions. OAs and their relation to TDs.
Define N$(n)$. Also remind them that we are only looking at $\lambda=1$. 
\end{comment}

A type of RGDD that is particularly useful for our Social Golfer solutions (as discussed in \cref{sect:SGP}) is a {\it resolvable transversal design}. These are straightforward to construct from a suitable set of MOLS and will provide a ready source for solutions. 

\begin{definition}\label{defn:TD}
    A transversal design of order $n$ and block size $k$, denoted TD$(k,n)$, is a triple $(V,G,B)$, where 
    $V$ is a set of $kn$ elements (or points); $G$ is a partition of $V$ into $k$ classes, called groups, each of size $n$; $B$ is a collection of $k$-subsets of $V$, called blocks; and every unordered pair of elements from $V$ is contained in one group,  or one block, but not both. If the blocks can be resolved into $n$ rounds of $n$ blocks, the design is said to be a resolvable transversal design, RTD$(k,n)$. An RTD$(k,n)$ can be obtained by removing one of the groups from a TD$(k+1,n)$.
\end{definition}

\begin{lemma}
    The existence of a TD$(k,n)$  is equivalent to the existence of $k-2$ MOLS of order $n$. Similarly, the existence of an RTD$(k,n)$
    is equivalent to the existence of $k-1$ MOLS of order $n$.
\end{lemma} 

\vspace{.2cm}
The RGDD shown in \cref{fig:RGDDExample} is an RTD$(4,7)$. Some useful RGDDs are given in \cref{lem:moreRGDDs}. The first result is from \cite{rees_GDD_k_with_k+1_groups_2000} and the second follows from a generalisation of a theorem of a result due to Shen Hao \cite[Theorem 2]{shen_RGDD4_JCMCC_1987}. The $9$-RGDD is from \cite{mathon_DiviibleSemiplanes_2007} and is constructed by removing $3$ parallel blocks from a $12$-GDD of type $3^{45}$.

\begin{lemma}\label{lem:moreRGDDs}
There exist: 
\begin{enumerate}
\item a $5\!-\!\!$-RGDD of type $g^{6}$ if and only if $g\equiv 0 \pmod{20}$,
\item a $5\!-\!\!$-RGDD of type $4^{30}$,
\item a $9\!-\!\!$-RGDD of type $3^{33}$.
\end{enumerate}
\end{lemma}

When there are insufficient MOLS to construct an RTD$(k,n)$, we can still create some parallel rounds by removing a group from an incomplete transversal design, ITD. An ITD is a particular type of incomplete GDD (IGDD), which is a group divisible design for which there is a subset of points, no two of which appear together in any block  (although subsets of them may appear together in the groups). 

\begin{definition}
An {\it incomplete transversal design} ITD$(n_{1},n_{2};k)$ is a 
is a tuple $(V,G,B,H)$, where 
    $V$ is a set of $kn_{1}$ points; $G$ is a partition of $V$ into $k$ groups, each of size $n_{1}$; $B$ is a collection of blocks of size $k$; and $H$ is a set of $kn_{2}$ points, called a hole. Every unordered pair of elements from $V$ is contained in a block or a group (but not both), or in the hole.  Pairs of points in $H$ may appear in a group, but not in a block. 
\end{definition}

Incomplete transversal designs can be constructed via a variety of means both direct and recursive. An ITD$(n_{1},n_{2};k)$ can be thought of as a (hypothetical) TD$(k,n_{1})$ with a (hypothetical) TD$(k,n_{1})$ removed. However, often an ITD$(n_{1},n_{2};k)$ exists when either or both of the transversal designs TD$(k,n_{1})$ and TD$(k,n_{2})$ do not exist.  

The example in \cref{fig:ITD} is an ITD$(10,2;6)$ originally from \cite{brouwer_4MOLSOrder10_1984}, neither a TD$(6,10)$ nor a TD$(6,2)$ exist. Note that the blocks are not arranged into rounds in this case (as the design is not resolvable). Two points of the hole appear in each group (they have been highlighted in the groups for reference).

\begin{figure}
\begin{equation*}
\boxed{
\begin{aligned}
&Blocks:\\
&[0,12,24,35,46,58],[0,14,26,38,41,52],[0,16,22,34,48,57],[0,11,28,32,42,55],\\
&[0,18,23,31,44,51],[0,10,25,30,43,59],[0,17,27,39,47,56],[0,15,20,33,49,50],\\
&[0,13,29,37,45,53],[0,19,21,36,40,54],[1,13,25,36,47,58],[1,15,27,38,42,53],\\
&[1,17,23,35,48,54],[1,12,28,33,43,56],[1,18,20,32,45,52],[1,11,26,31,40,59],\\
&[1,14,24,39,44,57],[1,16,21,30,49,51],[1,10,29,34,46,50],[1,19,22,37,41,55],\\
&[2,10,26,37,44,58],[2,16,24,38,43,50],[2,14,20,36,48,55],[2,13,28,30,40,57],\\
&[2,18,21,33,46,53],[2,12,27,32,41,59],[2,15,25,39,45,54],[2,17,22,31,49,52],\\
&[2,11,29,35,47,51],[2,19,23,34,42,56],[3,11,27,34,45,58],[3,17,25,38,40,51],\\
&[3,15,21,37,48,56],[3,10,28,31,41,54],[3,18,22,30,47,50],[3,13,24,33,42,59],\\
&[3,16,26,39,46,55],[3,14,23,32,49,53],[3,12,29,36,44,52],[3,19,20,35,43,57],\\
&[4,16,20,31,42,58],[4,10,22,38,45,56],[4,12,26,30,48,53],[4,15,28,36,46,51],\\
&[4,18,27,35,40,55],[4,14,21,34,47,59],[4,13,23,39,43,52],[4,11,24,37,49,54],\\
&[4,17,29,33,41,57],[4,19,25,32,44,50],[5,17,21,32,43,58],[5,11,23,38,46,57],\\
&[5,13,27,31,48,50],[5,16,28,37,47,52],[5,18,24,36,41,56],[5,15,22,35,44,59],\\
&[5,10,20,39,40,53],[5,12,25,34,49,55],[5,14,29,30,42,54],[5,19,26,33,45,51],\\
&[6,14,22,33,40,58],[6,12,20,38,47,54],[6,10,24,32,48,51],[6,17,28,34,44,53],\\
&[6,18,25,37,42,57],[6,16,23,36,45,59],[6,11,21,39,41,50],[6,13,26,35,49,56],\\
&[6,15,29,31,43,55],[6,19,27,30,46,52],[7,15,23,30,41,58],[7,13,21,38,44,55],\\
&[7,11,25,33,48,52],[7,14,28,35,45,50],[7,18,26,34,43,54],[7,17,20,37,46,59],\\
&[7,12,22,39,42,51],[7,10,27,36,49,57],[7,16,29,32,40,56],[7,19,24,31,47,53],\\
&[8,10,21,35,42,52],[8,11,22,36,43,53],[8,12,23,37,40,50],[8,13,20,34,41,51],\\
&[8,14,25,31,46,56],[8,15,26,32,47,57],[8,16,27,33,44,54],[8,17,24,30,45,55],\\
&[9,10,23,33,47,55],[9,11,20,30,44,56],[9,12,21,31,45,57],[9,13,22,32,46,54],\\
&[9,14,27,37,43,51],[9,15,24,34,40,52],[9,16,25,35,41,53],[9,17,26,36,42,50]
\\
&Groups:\\
&\begin{bNiceMatrix}
0\\ 1\\ 2\\ 3\\ 4\\ 5\\ 6\\ 7\\ {\bf 8}\\ {\bf 9}\\
\end{bNiceMatrix} 
\begin{bNiceMatrix}
10\\ 11\\ 12\\ 13\\ 14\\ 15\\ 16\\ 17\\ {\bf 18}\\ {\bf 19}\\
\end{bNiceMatrix} 
\begin{bNiceMatrix}
20\\ 21\\ 22\\ 23\\ 24\\ 25\\ 26\\ 27\\ {\bf 28}\\ {\bf 29}\\
\end{bNiceMatrix} 
\begin{bNiceMatrix}
30\\ 31\\ 32\\ 33\\ 34\\ 35\\ 36\\ 37\\ {\bf 38}\\ {\bf 39}\\
\end{bNiceMatrix} 
\begin{bNiceMatrix}
40\\ 41\\ 42\\ 43\\ 44\\ 45\\ 46\\ 47\\ {\bf 48}\\ {\bf 49}\\
\end{bNiceMatrix} 
\begin{bNiceMatrix}
50\\ 51\\ 52\\ 53\\ 54\\ 55\\ 56\\ 57\\ {\bf 58}\\ {\bf 59}\\
\end{bNiceMatrix} 
\end{aligned}
}
\end{equation*}
\caption{ITD$(10,2;6)$. Neither a TD$(6,10)$ nor a TD$(6,2)$ exist. Intersection of each group with the hole is emphasised. 
\label{fig:ITD}
}
\end{figure}

Examples for $k=5$, $k=6$ and $k=7$ can be found in   \cite{ abel_Vmt_2008,abelDuIdempotentMOLS_2003,abelcolbournYinZhangITD_k5_1997,
brouwerRees_MoreMols_1982, colbourn_ITDConstructions_1996, horton_IncompleteLatinSquares_1974, stinson_EquivalenceFramesITDs_1986}. A number of general recursive constructions from \cite{brouwerRees_MoreMols_1982,wilson_LatinSquares_1974} are summarised in  \cite{furinoMiaoYin_book_1996}.

\begin{lemma}\label{lem:ITDExistence}
For an ITD$(n_{1},n_{2};k)$ to exist we must have that $n_{1}\geq (k-1)n_{2}$.
\end{lemma}

\begin{comment}
Theorem \ref{thm:directProductForBIBDs} (see for example 
 \cite[Theorem ~9.1]{greig_BIBDk=7-9_2001})  is a generalisation of the Singular Indirect Product Theorem \cite{mullinSchellenbergVanstoneWallisFrames_1981} that is useful for the construction of BIBDs. Recall from \cref{defn:pbd} that a PBD with $\lambda=1$, $v$ points and block sizes in $K$ is referred to as a PBD$(v,K)$. In particular, when $K=\{k\cup f\}$ and (when $f\neq k$ there is a single distinguished block of size $f$), we refer to this PBD as a PBD$(v, \{k,f^{\ast}\})$.

\begin{Theorem}\label{thm:directProductForBIBDs}
Suppose that an ITD$(v-f+a, a;k)$ and a PBD$(v, \{k,f^{\ast}\})$ exist. Suppose also, that $f\geq a$ and  a BIBD$(f+(k-1)a,k)$ exists.
Then there exists a BIBD$(v^{\prime},k)$ where $v^{\prime}=kv-(k-1)(f-a)$ \,$=$ \, $k(v-f) + (f+ (k-1)a)$.
\end{Theorem}
\end{comment}

\begin{comment}
The following construction for a BIBD$(2479,7)$ appears to be unpublished:
\begin{Theorem}\label{thm:abelsResult}
A  BIBD$(2479,7)$ exists. 
\end{Theorem}

\vspace{.2cm}
\noindent
{\bf Proof} Use \cref{thm:directProductForBIBDs} with $v=367$, $f=61$, $a=46$. The required PBD$(367,\{7,61^{\ast}\})$ is constructed by adding an infinite point to each of the $61$ rounds of an RBIBD$(306,6)$, and a block consisting of the infinite points. The required ITD$(352,46;7)$ exists by Table 4.14, \cite{codihandbook}. \qed{}

\todo{Explain how the ITD is constructed?}
\end{comment}

By removing the points from a group of an ITD$(n_{1},n_{2};k+1)$, we achieve an incomplete group divisible design on $n_{1}k$ points with blocks of size $k$. There are $(n_{1}-n_{2})$ complete parallel rounds (one for each of the removed points that are not in the hole) and $n_{2}$ partial parallel rounds (one for each of the removed points that are in the hole).

\cref{example:IGDD_50_110} follows from the fact that there is an ITD$(6,10;2)$ \cite{brouwer_4MOLSOrder10_1984} and an ITD$(6,22;3)$ \cite{colbourn_ITDConstructions_1996}.

\begin{example}
\label{example:IGDD_50_110}
There exists an incomplete group divisible design on $50$ points with blocks of size $5$, and eight disjoint complete parallel classes. Similarly, there is an incomplete group divisible design on $110$ points with blocks of size $5$ and $19$ disjoint complete parallel classes.    
\end{example}

We refer to the design which consists of the blocks and groups formed by removing a group from an ITD$(n_{1},n_{2};k+1)$ as an RITD$(n_{1},n_{2};k)$.

If we view the groups of a group divisible design as blocks, we can think of an RGDD as a resolvable design for which all of the blocks (apart from those in one parallel class) have block size $k_{1}$, and the blocks in one parallel class have size $k_{2}$ (i.e. the group size). A similar structure is a {\it uniformly resolvable pairwise balanced design} ($URD$) \cite{rees_URD_2_3_1996}. In this case there may be multiple block sizes, but all blocks in each parallel class have the same size.

\begin{definition}
A uniformly resolvable  design (URD) is a PBD whose blocks can be resolved into parallel classes in such a way that all blocks in a given parallel class have the same size.   
\end{definition}

 We refer to a URD with $v$ points and two block sizes, $k_{1}$ and $k_{2}$, where $k_{1}<k_{2}$ as a URD$(\{k_{1},k_{2}\}; v)$. The number of parallel classes of blocks of size $k_{1}$ and $k_{2}$ are denoted by $r_{k_{1}}$ and $r_{k_{2}}$ respectively. 
 Some URDs with two small block sizes have been studied by various authors, for example, for $(k_{1},k_{2})=(2,3)$
\cite{rees_URD_2_3_1996}, $(2,4)$ \cite{dinitzLingDanziger_URD_2_4_2009, weiGe_URD_2_4_2017},
$(3,4)$ \cite{schusterGe_URD_3_4_2010,schuster_URD_3_4_2013, weiGe_URD_3_4_2016} and $(3,5)$ \cite{schuster_URD_3_5_2009}.

Some URDs that will be useful for our allocations are those for which $v\equiv 8 \pmod{12}$ and no $4\!-\!\!$-RGDD of type $2^{v/2}$ exists (see \cref{lem:rggd4_2}). The results in \cref{lem:urd} are from \cite{weiGe_URD_2_4_2017} and \cite{schusterGe_URD_3_4_2010} (for $n=100$).
%$u\in\{46,70,82,94,100,118,130,202,214\}$.
% v in 92 140 164 188 200 236 260 404 428

\begin{lemma}\label{lem:urd}
There exists a URD$(\{2,4\};2n)$ with $r_{2}=4$ and $r_{4}= (2n-5)/3$ for $n\in\{46,70,82,94,100,118,130,202,214\}$. 
\end{lemma}

\cref{example:urd5} is from the published best known Social Golfer solution for $v=30$, $k=5$ \cite{mathGames,csplib:prob010_results} \,\,(we added the rounds of blocks of size $2$ ourselves). 

\begin{example}\label{example:urd5}
There exists a $URD(\{2,5\};30)$ with $r_{2}=5$ and $r_{5}= 6$ (shown in \cref{fig:urdExample}). 
\end{example}

\begin{figure}
\begin{equation*}
\boxed{
\begin{aligned}
&[0,1,2,3,4],[5,6,7,8,9],[10,11,12,13,14],[15,16,17,18,19],[20,21,22,23,24],\\
&[25,26,27,28,29]\\
\\
&[0,5,10,15,20],[1,6,12,18,29],[2,7,13,23,28],[3,8,19,24,27],[4,14,17,22,26],\\
&[9,11,16,21,25]\\
\\
&[0,7,12,19,21],[1,5,11,22,27],[2,8,10,17,29],[3,14,15,23,25],[4,6,16,24,28],\\
&[9,13,18,20,26]\\
\\
&[0,6,17,23,27],[1,9,10,19,28],[2,5,14,18,21],[3,13,16,22,29],[4,8,12,20,25],\\
&[7,11,15,24,26]\\
\\
&[0,9,14,24,29],[1,8,16,23,26],[2,6,11,19,20],[3,5,12,17,28],[4,13,15,21,27],\\
&[7,10,18,22,25]\\
\\
&[0,8,11,18,28],[4,5,19,23,29],[3,6,10,21,26],[7,14,16,20,27],[2,9,12,15,22],\\
&[1,13,17,24,25]\\
\\
&[0,13],[1,7],[2,16],[3,9],[4,10],[5,24],[6,14],[8,15],[11,23],[12,26],[17,20],[18,27],[19,25],\\
&[21,29],[22,28]\\
\\
&[0,16],[1,14],[2,25],[3,7],[4,9],[5,13],[6,15],[8,22],[10,23],[11,17],[12,27],[18,24],[19,26],\\
&[20,29],[21,28]\\
\\
&[0,22],[1,15],[2,24],[3,11],[4,18],[5,26],[6,25],[7,29],[8,13],[9,23],[10,27],[12,16],[14,19],\\
&[17,21],[20,28]\\
\\
&[0,25],[1,20],[2,26],[3,18],[4,11],[5,16],[6,13],[7,17],[8,21],[9,27],[10,24],[12,23],[14,28],\\
&[15,29],[19,22]\\
\\
&[0,26],[1,21],[2,27],[3,20],[4,7],[5,25],[6,22],[8,14],[9,17],[10,16],[11,29],[12,24],[13,19],\\
&[15,28],[18,23]
\end{aligned}
}
\end{equation*}
\caption{A $URD(\{2,5\};30)$ with $r_{2}=5$ and $r_{5}=6$.
\label{fig:urdExample}
}    
\end{figure}

\section{(Incomplete) Transversal design constructions: (Incomplete) orthogonal arrays, difference matrices and quasi-difference arrays}

Most of the definitions in this section are (or are similar to those) from \cite[III.3.2,~III.4.1,~VI.17.1~and~VI.17.44]{codihandbook} unless stated otherwise. 
A structure that is closely related to a transversal design is an orthogonal array (OA):

\begin{definition}
    An orthogonal array OA$(k,s)$ is a $k\times s^{2}$ array with entries from an $s$-set having the property that in any two rows, each (ordered) pair of symbols from $S$ occurs exactly once. 
\end{definition}

Existence of an orthogonal array OA$(k,n)$ is equivalent to existence of a set of $k-2$ MOLS$(n)$ or a TD$(k,n)$. See 
\cite[Remark~III.3.10]{codihandbook}.

\begin{lemma}
    Let $A$ be an OA$(k,n)$ on the $n$ symbols in $X$. On $V=X\times \{0,\ldots, k-1\}$ (of size $kn$), form a set $B$ of $k$-sets as follows. For $0\leq j < n^{2}$, include $(a_{i,j},i):0\leq i < k\}$ in $b$. Then let $G$ be the partition of $V$ whose classes are $\{X\times \{i\}:0\leq i < k\}$. Then $(V,G,B)$ is a TD$(k,n)$. This process can be reversed to recover an OA$(k,n)$ from a TD$(k,n)$.
\end{lemma}

Orthogonal arrays (and thus, MOLS) are often constructed from difference matrices. First we define a difference matrix:

\begin{definition}\label{defn:diffMatrix}
Let $T$ be an abelian group of order $t$. A $(t,k;\lambda)$ difference matrix (or DM)  over $T$ is a $k\times t\lambda$ matrix $D=(d_{x,y})$ with entries from $T$ such that for each $i,j$ satisfying $0\leq i<j< k$, the multiset $\{d_{i,l}-d_{j,l}:0\leq l < t\lambda\}$ (the {\it difference list}) contains every element of $T$ $\lambda$ times.
\end{definition}

\begin{example}
Six MOLS of order $35$ were given in \cite{janiszczakStaszewskiCodesLS_2019}. These were obtained by using permutation codes. Later, Ingo Janiszczak noted that this construction gives a set of 6 MOLS$(35)$ which can  be obtained 
more simply from the following $(35,7;1)$ difference matrix.  Let $D_1$ and $D_2$ be respectively,  the following $7 \times 1$ and 
$7 \times 17$ arrays, and then let $D_3$ be the array obtained by interchanging the following pairs of rows 
in $D_2$:  $2$ and $3$,  $4$ and $5$,  $6$ and $7$. Then  $D= [\, D_1 \,  | \, D_2 \, | \, D_3]$ is the required $(35,7;1)$ difference matrix.

\begin{center} 
 $D_1 = \left (
\begin{array}{cc}
   0 \\
   22 \\
   22 \\
   11 \\
   11 \\
    4 \\
    4 \\
\end{array} \right ), \ \ \ D_2 = \left ( \begin{array}{cccccccccccccccccc}
    0 & 0  &  0  &  0 &  0 &  0 &  0 &  0 &  0 &  0  &  0  &  0 &  0 &  0 &  0 &  0 &  0 \\ 
    4 & 21 &  3  & 20 &  2 & 19 &  1 & 18 &  0 & 17  & 34  & 16 & 33 & 15 & 32 & 14 & 31  \\
    5 & 23 &  6  & 24 &  7 & 25 &  8 & 26 &  9 & 27  & 10  & 28 & 11 & 29 & 12 & 30 & 13  \\
   19 & 34 & 15  &  1 & 13 & 33 & 18 & 12 & 27 &  8  & 29  &  6 & 26 & 25 &  5 & 20 &  0  \\
    3 & 23 &  7  & 21 &  9 & 24 &  4 & 10 & 30 & 14  & 28  & 16 & 31 & 32 & 17 &  2 & 22  \\
   23 &  7 & 33  & 32 & 13 & 21 &  9 & 28 & 29 &  2  &  0  & 16 &  5 & 24 & 25 & 17 & 12  \\
   20 &  1 & 10  & 11 & 30 & 22 & 34 & 15 & 14 &  6  &  8  & 27 &  3 & 19 & 18 & 26 & 31  \\
\end{array} \right ) . $ \end{center}
\end{example}

If $S$ is  a symmetric $v \times v$ matrix, then five MOLS$(v)$  of the form  $A_1$, $A_1^{T}$, $A_2$, $A_2^{T}$ and $S$ are called a 2-SOLSSOM$(v)$. In 
\cite[Theorem~4.1]{abgz2002}, it is noted that if $v$ is odd and there exists
a $(v,7;1)$ DM with the property that it has an order $2$ automorphism
that leaves one row unaltered and permutes the other six rows in pairs, 
then a $2$-SOLSSOM$(v)$ exists  (this method was used in \cite{abgz2002} to obtain 
a $2$-SOLSSOM$(55)$).  Hence the above $(35,7,1)$ DM gives a new
$2$-SOLSSOM$(35)$. In \cref{2solssomodd} we update 
\cite[Theorem~4.1]{abgz2002} for $2$-SOLSSOMs of odd order:

\begin{Theorem} \label{2solssomodd} 
 If $v$ is an odd positive integer, then there exists a 2-SOLSSOM$(v)$ except for 
 $v \in \{3,5\}$ and possibly for  $v \in \{15,21,33,39,45,51,65,87,123,135\}$.
    
\end{Theorem}

An OA$(k,n)$ is obtained from a $(v,k;1)$ difference matrix $D$ by developing the columns of $D$ over $T$. Difference matrices with $\lambda>1$ can 
sometimes be used  to obtain OAs with $\lambda = 1$.
If $N$ is an order $n$ subgroup of an abelian group $T$ where $v=\lambda n$,  then a $(v,k;\lambda)$ difference matrix  $D = (d_{x,y})$ over $T$
is also said to be a difference matrix over $(T,N)$ if the equality $d_{i,l} -d_{j,l}   =d_{i,s} - d_{j,s}$ implies that $d_{i,l}$ and $d_{i,s}$ are from different cosets of $N$ in $T$. In \cite{Todorov2012}, an OA$(6,14)$ is obtained from a $(14,6;2)$ difference matrix. A modified form of this difference matrix together with the corresponding six MOLS of order $14$ are  displayed in \cite{AbelLi2015}.
%JA: Note that if $\lambda =1$, the OA$(k,n)$ can be extended to an OA$(k+1,n)$ (by adding %an extra row - see \cite{abelMOLS_18_60_2015}).
% \lambda = 1 not > 1. If we want to state this, we shouldn't do it after talking
% about \lambda > 1.

\begin{figure}
    \centering
    \includegraphics[scale=0.55]{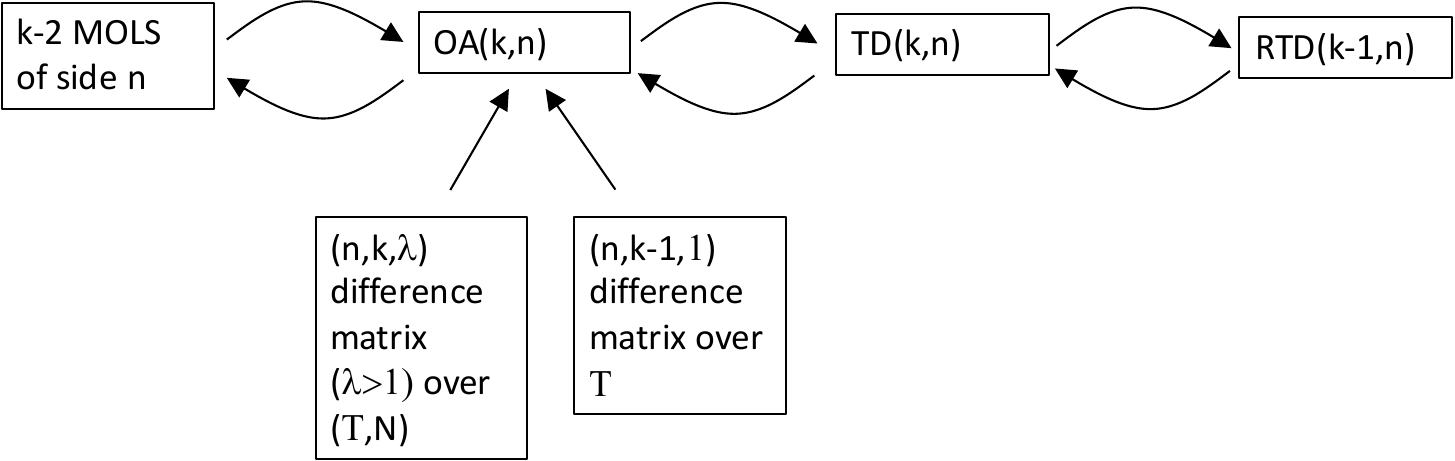}
    \caption{MOLS, OAs, TDs and Difference Matrices}
    \label{fig:complete}
\end{figure}

The relationship between MOLS, OAs, DMs and TDs is illustrated in \cref{fig:complete}.

As we have previously seen, we can obtain an RTD$(k,n)$ by removing the points in a single group of a TD$(k+1,n)$. 
\cref{const:DM14_rounds} will be useful when attempting to add points to rounds of blocks formed by an RTD$(5,14)$ for an SGA solution (see \cref{sect:adjacentBlockSizes}).

\begin{construction}\label{const:DM14_rounds}
The blocks of an RTD$(5,14)$ obtained from the $(14,6,2)$ difference matrix 
in {\rm\cite{Todorov2012}}  can be arranged in a $14\times 14$ grid (labelled $A$) where the blocks in each row of $A$ form a parallel class and the blocks in each of the first $10$ columns of $A$  form a parallel class. 
\end{construction}

{\bf Proof}
Let $D$ denote the  $(14,6;2)$ difference matrix  given in \cite{Todorov2012}. $D$ has $28$ columns with entries from $Z_{14}$, labelled
as $0,1, \ldots, 27$. Fourteen of these columns have first element $0$ and fourteen have first element $1$.

The blocks formed by expanding ten specific pairs of columns in $D$ form a parallel class. 
The appropriate ten pairs of columns are given in \cite{Todorov2012};  no column lies in more than one of these pairs.
 
Order the first $10$ columns of $D$ so that each of the columns labelled as $0,1,  \ldots, 9$ is the member of one of those $10$ pairs of columns with first entry $0$,
and for $x=0,1, \ldots  9$, the column labelled as $14+x$ is the column paired with the one labelled as $x$.  
Also choose the remaining $8$ columns so that (for each $y=0,1,2,3)$ the one labelled as $10+y$ has first element $0$,
and the one labelled as $24+y$ has first element $1$.

For the RTD$(5,14)$ obtained from $D$,  the groups are given as 
$H_{s} = \{ (1,s), (2,s),   \ldots , (13,s)\}$ for $1\leq s \leq 5$.

Now, for each column $c_{j} = (c_{0,j}, c_{1,j}, \ldots, c_{5,j})^{T}$ $(0 \leq j \leq 27)$ in $D$, 
we construct a block  $B_{j}  = \{(c_{1,j} ,1), (c_{2,j} ,2), 
c_{3,j} ,3), (c_{4,j} ,4), (c_{5,j},5)\}$.
Note that the values $c_{0,j}$  which all equal $0$ or $1$ do not appear in  the entries of $B_j$.

Let $A$ be the required $14 \times 14$ grid. Label the rows and columns of $A$ as $0,1, \ldots 13$.  % that is by the elements of $Z_{14}$.
For $i=0,1, \ldots, 13$, place $B_i$ and $B_{i+14}$ 
respectively in the $(i,0)$ and $(i,1)$ cells of $A$. 

Finally, for  $j=0,1, \ldots, 27$ and $z=2,4, \ldots, 12$, let $B_{j} +z  =$ 
$\{ (c_{1,j}+z,1),$  $(c_{2,j}+z,2),$  $(c_{3,j}+z,3),$ $(c_{4,j}+z,4),$  $(c_{5,j}+ z,5) \}$. Here, the additions $c_{s,j}+z$ are done within $Z_{14}$.
For $i=0,1, \ldots, 13$, place  $B_{i} + z$ and $B_{14+i} + z$ respectively,
in the $(i,z)$ and $(i,z+1)$  cells of $A$.
This gives the required grid.   \qedsymbol{} 

\vspace{0.2cm}

We can extend the concept of orthogonal arrays and difference matrices to incomplete orthogonal arrays and incomplete difference matrices (known as quasi-difference matrices). The following definition is from \cite{colbourn_ITDConstructions_1996}:

\begin{definition}
An IOA$(k,n;h)$ (which we call $A$)   is a $k\times (n^{2}-h^{2})$ array over a set $X$ of symbols of size $n$. There is an $h$-set $Y\subset X$ such that for every pair of distinct elements $i,j\in X$ and every pair of distinct rows $s,t$ in $A$, there is a unique column $u$ for which $A_{s,u}=i$ and $A_{t,u} = j$ unless both $i$ and $j$ lie in  $Y$, in which case no such column exists.
\end{definition}

Existence of an  IOA$(k,n;h)$ is equivalent to the existence of an ITD$(n,h;k)$. It follows from \cref{lem:ITDExistence}, that the existence of an IOA$(k,n;h)$ implies that $(k-1)h\leq n$.

Incomplete orthogonal arrays can be constructed from quasi-difference matrices  \cite{wilsonFewMoreSquares_1974,abelcolbournYinZhangITD_k5_1997}. 

\begin{definition}\label{defn:QDA}
    Let $G$ be an abelian group of order $n$. An $(n,k;\lambda,\mu; u)$-quasi-difference matrix (QDM) is a matrix $Q=(q_{ij})$ with $k$ rows and $\lambda(n-1+2u)+\mu$ columns, where each entry is either empty (usually denoted by $-$) or contains a single element of $G$. Each row contains exactly $\lambda u$ empty entries, and each column contains at most one empty entry. Furthermore, for each $0\leq i < j< k$, the multiset 
    $$\{q_{il} - q_{jl}: 0\leq l < \lambda(n-1+2u)+\mu,\;{\rm with}\;q_{il}\;{\rm and}\;q_{jl}\;{\rm not}\;{\rm empty}\}$$
    contains every nonzero element of $G$ $\lambda$ times and contains $0$ $\mu$ times.
\end{definition}

%Construction VI.17.47 in \cite{codihandbook} 
The following lemma states the link between a QDM and an incomplete orthogonal array, and is proved by following the construction in 
\cite[Construction~VI.17.47]{codihandbook}.

\begin{lemma}\label{lem:ITDFromQDM}
    If an $(n,k;\lambda, \mu;u)$-QDM exists and $\mu\leq \lambda$, then an ITD$_{\lambda}(k,n+u;u)$ exists. 
\end{lemma}

%{\bf Proof} Start with an $(n,k;\lambda, \mu;u)$-QDM $A$ over the group $G$. Append %$\lambda-\mu$ columns of zeros. Then select $u$ elements $\infty_{1},\ldots,\infty_{u}$ 
%    not in $G$ and replace the empty entries ($-$), each with one of these infinite %symbols, so that $\infty_{i}$ appears exactly once in each row, for $1\leq i\leq u$. %Develop the resulting matrix over the group $G$ (leaving infinite symbols fixed), to %obtain a $k\times (n^{2}+2nu)$ matrix $T$. Then $T$ is an incomplete orthogonal array %with $k$ rows and index $\lambda$, having $n+u$ symbols and one hole of size $u$. %\qedsymbol{}

% JA: Should we give this proof or cite a reference?

\vspace{0.2cm}
\cref{lem:QDM_rounds} is also useful when attempting to add points to rounds of blocks formed by a QDM for a SGA solution (see \cref{sect:adjacentBlockSizes}).

\begin{construction}\label{lem:QDM_rounds}
If an $(n,k+1;1,\mu;u)$-QDM exists and $\mu\leq 1$, then we can construct an $n\times (n+u)$ array of blocks of size $k$, where the blocks in each row of blocks form a parallel class and 
$n-(k-1)u$ of the columns consist of parallel blocks. 
\end{construction}

%JA: I need to check this proof
%AM: I have just changed this to hold for mu in {0,1}, not just 1 as we had it before
\noindent
{\bf Proof}  The proof is similar to that of \cref{const:DM14_rounds} and follows the steps of the ITD construction described in \cite[Construction~VI.17.47]{codihandbook}.
After adding the infinite points and $1-\mu$ columns of zeros, there are $n+2u$ columns of which $u$ contain an infinite point in the first position, and $k\times u$ contain an infinite point in another position. Reorder the columns so that those containing an infinite point in position $i$, for $0\leq i\leq k$ are 
$\{c_{j}:iu\leq j<(i+1)u\}$,  and those that contain no infinite points are $\{c_{j}:(k+1)u\leq j<n+2u\}$. Call the resulting matrix QDM$^{\prime}$ and construct the ITD$(k+1,n+u,u)$. The groups are $H_{j} = ((0,j), (1,j), \ldots , (n+u-1,j))$ for $0\leq j\leq k$. Each of the columns $c_{j}$ of QDM$^{\prime}$ leads to $n$ blocks (by adding elements of $G$) thus:
$c_{j} = (c_{0,j}, c_{1,j}, \ldots, c_{k,j})^{T}$ leads to the set of blocks $$B_{j} = \{b_{j}^{g}:b_{j}^{g} = ((c_{0,j}+g,0), (c_{1,j}+g,1), \ldots, (c_{k,j}+g,k)), \;g\in G\}$$

Place the blocks from the ITD derived from the last $n+u$ columns of QDM$^{\prime}$ in an $n\times (n+u)$ array so that the block arising from $c_{j}$ with first element $(i,0)$ is in position $(i,(j-u))$ in the array. 

Removing elements of group $H_{0}$ from each block, the blocks in each row form a parallel class. In addition, all of the columns in the array that are not formed from columns QDM$^{\prime}$ that contain infinite points, consist of parallel blocks (i.e. the last $n-((k-1)u)$ columns). \qedsymbol{}

\begin{construction}\label{cor:extraRoundQDA}
If an $(n,k+1;1,\mu;u)$-QDM exists and $\mu\leq 1$, then we can construct an $(n+1)\times (n+u)$ array of blocks of size $k$ where: each row of blocks is a parallel class and $u$ of the columns consist of parallel blocks. 
\end{construction}

\noindent
{\bf Proof} Construct the first $n$ rows as for the proof of \cref{lem:QDM_rounds}. A further round, from the blocks of the TD containing $(\infty_{1},0)$ should be added, with the $u$ such blocks which are from filling the hole of the IOA placed in columns $k*u,(k*u)+1,\ldots, ((k+1)*u)-1$. Those $u$ columns (of size $n+1$) will consist of parallel blocks. \qedsymbol{}
 
%A V$(m,t)$ vector, which is defined as follows:
\begin{definition}
    Let $q=mt+1$ (where $m, t$ are integers) be a prime power, and let $\omega$ be a primitive element in the finite field $F_{q}$. A V$(m,t)$ vector is a vector $(a_{1},\ldots, a_{m+1})$ for which, for each $k \in \{1,2, \ldots, m+1\}$ the
    differences $\{a_{i+k}-a_{i}: 1\leq i \leq m+1, i+k\neq m+2\}$ represent the $m$ cyclotomic     classes of $F_{mt+1}$.   Here the subscripts $i+k$ are calculated $\pmod{m+2}$.    In other words, for fixed $k$, 
    if $a_{i+k}-a_{i} = \omega^{mx+\alpha}$ and  $a_{j+k} - a_{j}=\omega^{my+\beta}$ where     $0 \leq \alpha, \beta  \leq m-1$
    and $ i \neq j$, then $\alpha \neq \beta$.
\end{definition}

\begin{Remark} A QDM can be obtained from a $V(m,t)$ vector (\cite[Construction~VI.17.49]{codihandbook}, \cite{abel_Vmt_2008}).
\end{Remark}

\begin{figure}
    \centering
    \includegraphics[scale=0.55]{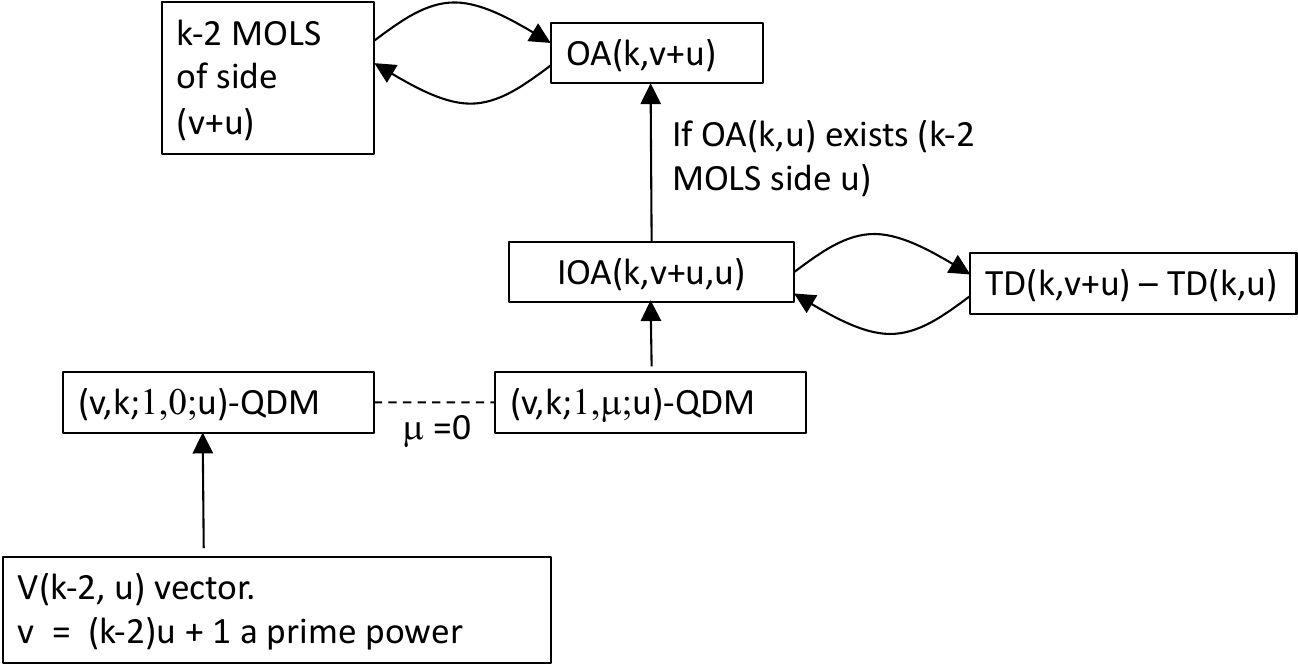}
    \caption{IOAs, ITDs, QDMs and $V(m,t)$ vectors.} 
    \label{fig:incomplete}
\end{figure}

The relationship between IOAs, ITDs and V$(m,t)$s is illustrated in \cref{fig:incomplete}.

\section{Allocations: Constructions and Examples}\label{sect:constructions}
We now return to optimal solutions to the SGP as stated in \cref{defn:SGP}. The problem is to create as many parallel rounds of $v/k$ blocks of size $k$. We refer to a set of rounds of $v/k$ blocks of size $k$ as a $(v,k)$ allocation. If there are $r$ rounds we refer to an $r$-round $(v,k)$ allocation, and if this is optimal, we refer to an optimal $r$-round $(v,k)$ allocation (or, more generally, as an optimal allocation). An optimal SGP solution is given by an optimal allocation, so from henceforth we will refer to allocations rather than SGP solutions. 
We will refer to the number of rounds in an optimal $(v,k)$ allocation as optR$(v,k)$.

The following construction, originally due to Sharma and Das \cite{sharmaDas_RBIBD_1985} 
is presented in \cite{harveyWinterer_SGP_MOLR_2005}, and provides a 
lower bound for $r$. We refer to this construction as MOLRs$(k,n)$.

\begin{construction}\label{const:sharmaAndDas}
%{\bf Sharma and Das Construction}
If there exist $g$ MOLS$(n)$ where $n=v/k$, then there is a $g+1$ round $(v,k)$ allocation. If $\Gamma$ is the complete graph on $v$ points, there are $k$ cliques of size $n$ in $\Gamma$ whose edges do not appear in the blocks. These cliques are $$C_{i}=\{(n*i) + j: 0\leq j < n\},\;0\leq i <k$$
If $k$ divides $n$ then we can add further rounds from the cliques, and there is an optR$(n,k)+g+1$ round $(v,k)$ allocation. 
%it works when there aren't enough MOLS for an RTD
%but when there are, we get the same number of rounds
%better to use an RTD because then we have the groups (to use as cliqueS) %for reduction purposes
%It gives a lower bound in all cases

Consider the $g$ rectangles $L_{\alpha}$, $\alpha \in \{0,1,\ldots, g-1\}$ formed by taking the first $k$ rows of the $g$ MOLS$(n)$. These are mutually orthogonal Latin rectangles (MOLRs) as defined in \cite{harveyWinterer_SGP_MOLR_2005}. Let $R_{t}=\{b_{t,j}, j=0,1,\ldots, n-1\}$ 
denote the $t$'th round of blocks, for $t=0,1,\ldots, g$.

Construct a $k\times n$ array $G$ from the elements $0,1,\ldots, v-1$ as shown in \cref{fig:MOLR_array}. 

\begin{figure}
$\left[\begin{array}{ccccc}
0&1&2&\ldots&n-1\\
n&n+1&n+2&\ldots&2n-1\\
\vdots&\vdots&\vdots&\ldots&\vdots\\
(k-1)n&(k-1)n+1&(k-1)n+2&\ldots&kn-1
\end{array}
\right]
$
\caption{Array $G$ for Sharma and Das Construction \label{fig:MOLR_array}}
\end{figure}

The first round of blocks, $R_{0}$ consists of the columns of $G$.
Subsequent rounds, $R_{t}$, $1\leq t\leq g$ are formed by superimposing a Latin rectangle, $L_{t-1}$ on $G$ to form an array of pairs $G^{t}$, where, for $0\leq i <k$, $0\leq j<n$, 
$G^{t}[i][j]=(G[i][j], L_{t-1}[i][j])$. Then, for $0\leq j<n$, $b_{t,j}$ is the list of all values $val$ for which $(val,j)$ is an element of $G^{t}$.

The rows of $G$ form unused cliques of size $n$. If $k$ divides $n$, they can be broken up to form optR$(n,k)$ extra rounds of blocks. 
\end{construction}

\begin{example}\label{example:sharmaDas}
For $v=36$ and $k=6$, $n=6$, since there are $1$ MOLS$(6)$ and $k$ divides $n$, we can construct $3$ rounds using MOLRs$(6,6)$.

Array $G$, the single Latin rectangle formed from the first $6$ rows of the Latin square of order $6$ (i.e. the Latin square itself), $L_{0}$, and $G^{1}$ formed by superimposing $G$ and $L_{0}$, are shown in \cref{fig:MOLR_array36}.

\begin{figure}[t!]
\centering
\begin{subfigure}[t]{0.26\textwidth}
\centering
$\left[\arraycolsep=2.5pt\begin{array}{cccccc}
0&1&2&3&4&5\\
6&7&8&9&10&11\\
12&13&14&15&16&17\\
18&19&20&21&22&23\\
24&25&26&27&28&29\\
30&31&32&33&34&35\\
\end{array}
\right]
$
\caption{$G$}
\end{subfigure}%
~
\begin{subfigure}[t]{0.18\textwidth}
\centering
$\left[\arraycolsep=2.5pt\begin{array}{cccccc}
0&1&2&3&4&5\\
1&0&3&2&5&4\\
2&3&4&5&0&1\\
3&2&5&4&1&0\\
4&5&0&1&2&3\\
5&4&1&0&3&2\\
\end{array}
\right]
$
\caption{$L_{0}$}
\end{subfigure}%
~ 
\begin{subfigure}[t]{0.3\textwidth}
\centering
$\left[\arraycolsep=1.4pt\begin{array}{cccccc}
(0,0)&(1,1)&(2,2)&(3,3)&(4,4)&(5,5)\\
(6,1)&(7,0)&(8,3)&(9,2)&(10,5)&(11,4)\\
(12,2)&(13,3)&(14,4)&(15,5)&(16,0)&(17,1)\\
(18,3)&(19,2)&(20,5)&(21,4)&(22,1)&(23,0)\\
(24,4)&(25,5)&(26,0)&(27,1)&(28,2)&(29,3)\\
(30,5)&(31,4)&(32,1)&(33,0)&(34,3)&(35,2)\\
\end{array}
\right]
$
\caption{$G^{1}$}
\end{subfigure}
\caption{$G$, $L_{0}$ and $G^{1}$ for MOLRs$(6,6)$\label{fig:MOLR_array36}}
\end{figure}

The $4$ rounds obtained via the Sharma and Das construction are shown in 
\cref{fig:MOLRs}.

\begin{figure}
\begin{equation*}
\boxed{
\begin{aligned}
&\bf{Round}\;R_{0}\; ({\rm columns\;of\;}G):\\
&[0,6,12,18,24,30],[1,7,13,19,25,31],[2,8,14,20,26,32],[3,9,15,21,27,33],\\&[4,10,16,22,28,34],[5,11,17,23,29,35]\\
\\
&\bf{Round}\;R_{1}\; ({\rm from\;}G^{1}):\\
&[0,7,16,23,26,33],[1,6,17,22,27,32],[2,9,12,19,28,35],[3,8,13,18,29,34],\\&[4,11,14,21,24,31],[5,10,15,20,25,30]\\
\\
&\bf{Round}\;R_{2}: ({\rm extra\;round\;from\;rows\;of\;}G)\\
&[0,1,2,3,4,5],[6,7,8,9,10,11],[12,13,14,15,16,17],[18,19,20,21,22,23],\\&[24,25,26,27,28,29],[30,31,32,33,34,35]\\
\\
\end{aligned}
}\end{equation*}
\caption{$3$ round $(36,6)$ allocation using Sharma and Das construction (MOLRs$(6,6)$).
\label{fig:MOLRs}
}
\end{figure}

\end{example}

%ADD LEMMAs starterLem and starterLemExamples when finished
\begin{construction}\label{const:starterblocks}
{\bf Rounds from starter blocks} 
For some $(v,k)$ the optimal allocation is obtained by expanding  a set of $r$ starter blocks. 
    Let $g=v/k$. If we have a set of $r$ starter blocks which have the property that:
 \begin{enumerate}
     \item each starter block contains $0$ and none of the differences within any block is a multiple of $k$,
     \item no two starter blocks intersect at any point other than $0$ and
     \item for any two starter blocks $a=\{a_{i}:0\leq i <k \}$ and $b=\{b_{j}:0\leq j < k\}$, if there are two pairs $(a_{i_{1}}, b_{j_{1}})$ and $(a_{i_{2}}, b_{j_{2}})$ where $a_{i_{1}}, a_{i_{2}}\in a$ and $b_{j_{1}}, b_{j_{2}}\in b$ and the differences $(a_{i_{1}} - b_{j_{1}})$ and $ (a_{i_{2}}- b_{j_{2}})$ are equal modulo $v$, then the difference is not a multiple of $k$, 
 \end{enumerate}
 then we can construct a set of $r$ parallel rounds. In addition, the graph on the unused pairs contains $k$ cliques of size $(v/k)$. If $n=v/k$ and $k$ divides $n$, then we can add a further optR$(n,k)$ rounds from the cliques.

\vspace{.2cm}
\noindent
The rounds are generated by adding $k*m$ to each starter block, for $0\leq m\leq g-1$.
Each round contains distinct points (this follows from condition (1)). No pair is contained in more than one block. 
 Suppose for a contradiction that pair $(p_{0},p_{1})$ were to appear in two blocks, $a^{\prime}$ and $b^{\prime}$ say, derived from starter blocks $a$ and $b$.  Then $a^{\prime}=a+km_{1}$ and $b^{\prime}=b+km_{2}$ for some values $0\leq m_{1},m_{2}<g$ and $(p_{1},p_{2})=(a_{i_{1}}+km_{1},a_{i_{2}} + km_{1})$ and $(b_{j_{1}}+km_{2},b_{j_{2}} + km_{2})$. But this contradicts condition (3). The cliques in each case are on the sets of points $S_{i}=\{i+jk:0<=j<(v/k)\}$ for $0\leq i < k$.
 \end{construction}
 
\vspace{.2cm}
\noindent
We refer to \cref{const:starterblocks} as ownSG$(v,k)$, and to the cliques as {\it groups}.

\begin{example}
A $7$-round $(60,6)$ allocation can be obtained from a set of $7$ starter blocks. 
\end{example}

\cref{table:ownSGResults} presents  values $(v.k)$ for cases where an $r$-round ownSG$(v,k)$ allocation is an optimal allocation. Starter blocks for each of these cases are given in \cref{appendix:ownSG}. Note the only case for which $k$ divides $v/k=n$ is $(v,k)=(98,7)$. There is one extra round for this example. 
\begin{table}
\centering
%\begin{footnotesize}
% \begin{tabular}{|p{7.5cm}|p{4.5cm}|l|}
\begin{tabular}{|l|l||l|l||l|l||l|l||l|l|}
\hline
$v,k$ & $r$ & $v,k$ & $r$& $v,k$ & $r$&$v,k$ & $r$&$v,k$ & $r$\\
\hline
$60,6$ & $7$ & $70,7$ & $7$ & $80,8$ & $5$  &$84,6$ &$9$& $84,7$ &$7$  \\
$90,6$  &$10$ & $90,9$ &$5$   &$96,8$ &$6$ &$98,7$ &$8$  &$105,7$ &$9$  \\
$112,8$ &$8$   &$120,6$ &$13$ & $126,7$&$10$ &$126,9$ &$7$  &$132,6$ &$14$  \\
 $135,9$ &$7$&$156,6$ &$15$&&&&&&\\
\hline
\end{tabular}
\caption{Pairs $(v,k)$ for which  ownSG$(v,k)$ is an optimal $r$-round allocation. \label{table:ownSGResults}}
\end{table}

\section{Optimal (v,k) allocations}\label{sect:maxSGPSolutions}
The maximum number of rounds possible in a $(v,k)$ allocation is $r^{\prime}(v,k)=\lfloor\frac{(v-1)}{(k-1)}\rfloor$. We say that a $(v,k)$ allocation is maximal if it has $r^{\prime}(v,k)$ rounds (it is clearly optimal in this case). A maximal allocation is a resolvable maximum packing as defined in \cite{geLamLingShen_ResolvableMaximalPackingsQuads_2005}.

We will show that when $v\equiv 0 \pmod{3}$ and $k=3$ a maximal allocation is always possible. For $v\equiv 0 \pmod{4}$ and $k=4$ a maximal allocation is almost always possible, and in the few cases where it is not, an optimal allocation that is close to maximal exists. We also discuss cases where a maximal allocation is possible for $v\equiv 0 \pmod{5}$ and $k=5$. 

\subsection{When all blocks have size three}
For an allocation to exist, we must have $v\equiv 0 \pmod{3}$. Hence, either $v\equiv 3 \pmod{6}$ or $v\equiv 0 \pmod{6}$. 

If $v\equiv 3 \pmod{6}$ then we can use the blocks of a RBIBD$(v,3)$ (see \cref{defn:bibd}), i.e. a Kirkman triple system on $v$ points (KTS$(v)$), as our allocation. Since a KTS$(v)$ has $(v-1)/2$ classes, the allocation is maximal. For every $v\equiv 3 \pmod{6}$ there is a KTS$(v)$ \cite{ray_chaudhuriWilson_RBIBD3_1971}, so a maximal allocation exists in all cases.

If $v\equiv 0 \pmod{6}$ then we can use the blocks of a $3\!-\!\!$-RGDD of type $2^{\frac{v}{2}}$ (see \cref{defn:gdd}), i.e. a nearly Kirkman triple system on $v$ points (NKTS$(v)$), as our allocation. Since NKTS$(v)$ has $(v-1)/2$ classes, the allocation is maximal. By \cref{lemma:nktsExistence}, an NKTS$(v)$ exists for every $v\equiv 3 \pmod{6}$ apart from $v=6$ or $v=12$, a maximal solution exists for all $v\equiv 0 \pmod{6}$, $v\geq18$.  

\begin{comment}
The fact that no NKTS$(v)$ exists for $v=6$ or $v=12$ was proved by Kotzog and Rosa in \cite{kotzigRosa_NKTS_1974}, where constructions for infinitely many examples were given. Constructions for all cases apart from $v=84$, $v=102$ and $174$ were found by Baker and Wilson \cite{bakerWilson77} (a slight correction to starter blocks for case $v=36$ was later given in \cite{reesWallis2002}).  Cases $v=102$ and $v=174$ were solved by Brouwer \cite{brouwer_NewNKTS_1978} and case $v=84$ by Rees and Stinson 
\cite{reesStinson_RGDD3_1987}. A good summary of constructions of Nearly Kirkman Triple Systems is given in \cite{reesWallis2002}. 
\end{comment}

%
\subsection{When all blocks have size four}\label{SGP:k4}
Since $v\equiv 0 \pmod{4}$ either $v\equiv 4 \pmod{12}$, $v\equiv 0 \pmod{12}$ or $v\equiv 8 \pmod{12}$. It has been previously observed in \cite{stinsonWeiYin_Packings_2007} that a resolvable maximum packing is equivalent to an RBIBD$(v,4)$, a $4\!-\!\!$-RGDD with groups of size $3$, and a $4\!-\!\!$-RGDD with groups of size $2$ respectively. We consider these cases below. 

If $v\equiv 4 \pmod{12}$ we can use the blocks of an RBIBD$(v,4)$ for a maximal allocation. For every $v\equiv 4 \pmod{12}$ there is an RBIBD$(v,4)$ \cite{hananiRay_ChaudhuriWilson2006-orig1972}, so a maximal allocation exists for all cases. 

If $v\equiv 0 \pmod{12}$ then the number of rounds in a maximal allocation is $(v-3)/3$ which can be achieved by using the rounds of blocks of a $4\!-\!\!$-RGDD of type $3^{\frac{v}{3}}$. By \cref{lem:rggd4_1} such an RGDD exists for all $v\geq 24$, so a maximal allocation exists for all cases $v\geq 24$.

If $v\equiv 8 \pmod{12}$ then a maximal allocation can be achieved by using the rounds of blocks of a $4\!-\!\!$-RGDD of type $2^{\frac{v}{2}}$. By \cref{lem:rggd4_2} such an RGDD exists for all $v=2u$, except for $u=4$, $u=10$ and possibly when $u\in\{46,70,82,94,100,118,130,202,214\}$.

For the exceptions listed above, for $u>10$, by \cref{lem:urd} there exists a $URD(\{2,4\};v)$ with $r_{2}=4$ and $r_{4}= (v-5)/3$.
Hence in these cases we use the $(v-5)/3$ rounds of blocks of size $4$ for our allocations. 

%\cite{harachwi2006}
\subsection{When all blocks have size five}
Since $v\equiv 0 \pmod{5}$, $v\equiv 5$, $0$, $10$ or $15 \pmod{20}$.

If $v\equiv 5 \pmod{20}$ we can use the blocks of an RBIBD$(v,5)$ when it exists for a maximal allocation. An RBIBD$(v,5)$ exists in this case unless $v\in\{45,345,465,645\}$ (see \cref{lemma:rbibdExistence}), so a maximal solution exists for all other cases. 

When $v\in\{45,345,465,645\}$ we can use an RTD$(5,v/5)$ to obtain an allocation. The allocation is not maximal in any case, indeed for the larger values since there are only $v/5$ rounds this is far from the case. However, currently no better solution is known. 

For the remaining cases, ideally we would use a  $5\!-\!\!$-RGDD with small group size, however there are not many known examples of these. We can always use an RTD$(5,v/5)$ to obtain an allocation with $r=v/5$ rounds, except when $N(v/5)<4$ (i.e. when $v\in \{30,50,110\}$). 

When $v=30$ there is a solution with $6$ rounds, namely the rounds of blocks of size $5$ of a URD$(\{2,5);30)$ (see \cref{example:urd5}). 

When $v=50$ or $v=110$ we use the $8$ ($19$) parallel rounds of blocks of size $5$ of the incomplete group divisible design on $50$ ($110$) points given in \cref{example:IGDD_50_110}.   

If $v \equiv 0 \pmod{20}$ we can use the $v/5$ rounds of an RTD$(5,v/5)$ if there are a sufficient number of MOLS. 

When $v\equiv 0 \pmod{120}$, a better solution is to 
use the  $5\!-\!\!$-RGDD of type $g^{6}$ from \cref{lem:moreRGDDs}
which has $5v/24$ rounds. Since the group size is divisible by $5$ in this case, we can use \cref{lem:increasingRounds} to increase the number of rounds further. 

\begin{lemma}\label{lem:increasingRounds}
    If $g$ is divisible by $k$, then a uniform $k$-RGDD with groups of size $g$ gives an allocation with optR$(g,k)+((v-g)/(k-1))$ rounds. 
\end{lemma}

\begin{comment}
\noindent
{\bf Proof} In addition to the $((v-g)/(k-1))$ rounds of blocks from the RGDD, divide the groups into optR$(g,k)$ rounds of blocks of size $k$ and combine them to form  optR$(g,k)$ additional rounds of the allocation.   \qedsymbol{}
\end{comment}

When $v=120$, a better solution still is to use the blocks of the $5\!-\!\!$-RGDD of type $4^{30}$ in  \cref{lem:moreRGDDs} (which has $29$ rounds).

\subsection{When $v\leq 150$ and all blocks have size k, where k is at least 3}
The choice of method to generate optimal $(v,k)$ allocations is governed by \cref{alg1}. Our results for $12\leq v\leq 150$, where there are at least three rounds, are given in \cref{appendix:tables}. Note that the method in each case (as denoted in the tables) is provided as a comment (i.e. enclosed within \texttt{/*...*/}) in \cref{alg1}. Note that we use the shorthand RGDD$(v,k,g)$ to denote a $k$-RGDD with groups of size $g$. If the blocks of an RTD, or an RITD, or rounds from the MOLRs$(k,n)$ or ownSG$(v,k)$ constructions provide a solution, there are $k$ disjoint unused cliques of size $n$ (the groups). If $k$ divides $n$, then the groups themselves can be divided up to form extra rounds (in a similar way to the proof of \cref{lem:increasingRounds}). 

\begin{algorithm}
\caption{\enskip For selecting construction for optimal allocation, $v$ points and single block size $k$}\label{alg1}
\begin{verbatim}
optimal(v,k){
//assume that k divides v
   n=v/k;
   EXCEPTIONSA = {8, 20, 92, 140, 164, 188, 200, 236, 260, 404, 428}
   if(k==3){
      if(v(mod 6)==3) return blocks of KTS(v); /*KTS(v)*/
      if(v(mod 6)==0){
         if(v>=18) return blocks of NKTS(v) /*NTKS(v)*/
         if(v==6) return single round of two blocks
         if(v==12) return blocks of RTD(3,4) /*RTD(3,4)*/
      }
   }
   if(k==4){
      if(v(mod 12)==4) return blocks of RBIBD(v,4); /*RBIBD(v,4)*/
      if(v(mod 12) == 0 and v>12)
         return blocks of 4-RGDD(3^{v/3}) /*RGDD(v,4,3)*/
      
      if(v(mod 12) == 8 and v>8){
         if(v is not in EXCEPTIONSA)
             return blocks of 4-RGDD(2^{v/2}) /*RGDD(v,4,2)*/
         if(v in EXCEPTIONSA)
             return blocks of size 4 of a URD(4,2) with r_{2} = 4
             /*URD(v,4,2)*/
      }
   
      if(v==8 or v==12) return single round of v/4 blocks
   }

   else{
      if(v(mod k(k-1))==k and RBIBD(v,k) exists) 
         return blocks of RBIBD(v,k); /*RBIBD(v,k)*/
      else{
          if(k-RGDD of type g^{v/g} exists for g<=v/k){
              choose smallest such g and 
              if(k does not divide g) return blocks of RGDD;
                 /*RGDD(v,k,g)*/
              else{ 
                    return blocks of RGDD plus opt(g,k) extra rounds 
                    from the groups /*RGDD(v,k,g) + G(opt(g,k))*/
              }
          }
          else{
             return special(v,k);
          }
    }

  }      
   
}
\end{verbatim}
\end{algorithm}

\begin{algorithm}
\caption{\enskip  Supplement to \cref{alg1} for special cases}\label{alg2}
\begin{verbatim}
special(v,k){
//assume that k>4 divides v and v<=200
   n=v/k;
   let r be maximum of:
      choice 1. number of rounds from MOLRs(n,k) construction including 
                extra rounds from groups if k divides n 
                /*MOLRs(k,n) or MOLRs(k,n)+G(opt(n,k))*/
      choice 2. number of rounds from ownSG(v,k) construction including 
                extra rounds from groups if k divides n 
                /*ownSG(v,k) or ownSG(v,k)+G(opt(n,k))*/
      choice 3. number of rounds of blocks of size k from a
                URD({k1,k};v) /*URD(v,k,k1)*/
      choice 4a.(n-n2) where ITD(n, n2; k+1) exists and k does not    
                divide n /*RITD(n,n2;k)*/
      choice 4b.(n-n2+optR(n,k)) where ITD(n, n2; k+1) exists and k 
                divides n 
                /*RITD(n,n2;k) + G(opt(n,k))*/
      
    return r rounds of blocks using min choice that gives r rounds
}
\end{verbatim}
\end{algorithm}

%removed this as we don't seem to need it any more?
%choice 5. number of rounds from an existing, published SGP  
%solution  /*SG(v,k)*/

%
\section{When blocks have adjacent sizes - SGA}\label{sect:adjacentBlockSizes}
We now consider the problem where not all blocks have the same size. This problem arises when an allocation is required for $v$ participants, but $v$ has no divisors greater than $2$. If we permitted all possible block sizes, the number of solutions would be very large. We therefore restrict the number of block sizes to $2$ and insist that, for block size sets $K=\{k_{1},k_{2}\}$, $k_{2}=k_{1}+1$ (so the block sizes are {\it adjacent}). 
We refer to the problem of maximising the number of rounds in this case as the Social Golfer Problem with adjacent block sizes (SGA) and, as explained in \cref{subsect:SGA}, in this paper we only consider $K=\{4,5\}$ and $\{5,6\}$. Results for $v$ up to $150$ are contained in our tables in \cref{appendix:tables}.  For completeness, we find solutions for all $v$ (regardless of whether $v$ has other suitable divisors) and for all values $m_{1}$, $m_{2}$ for which $k_{1}m_{1}+k_{2}m_{2}=v$. For brevity, we refer to a $(v,k_{1},k_{2},m_{1},m_{2})$-allocation. 

\begin{definition}
For values $v$, $k_{1}$, $k_{2}$, $m_{1}$ and $m_{2}$ where $k_{2}=k_{1}+1$ and $k_{1}m_{1}+k_{2}m_{2}=v$, we refer to the optimal $(v+m_{1}, k_{2})$ allocation and the optimal $(v-m_{2},k_{1})$ allocation as the superior $(v+m_{1}, k_{2})$-allocation and the inferior $(v-m_{2},k_{1})$-allocation respectively. 
\end{definition}

In most cases, an optimal $(v,k_{1},k_{2},m_{1},m_{2})$-allocation is obtained by removing $m_{1}$ points from a single unused clique in the superior $(v+m_{1}, k_{2})$-allocation. 

\subsection{Removing points from an unused clique in the superior allocation}
If the superior allocation is from the blocks of an RTD$(k,n)$,  then any group provides a suitable clique. If $k$ does not divide $n$, the subsequent solution is denoted RTD$(k,n)-m_{1}$ and there are $n$ rounds. If $k$ divides $n$, optR$(n,k)$ additional rounds from the groups are included in the superior allocation, RTD$(k,n)+G($optR$(n,k))$. If $m_{1}=1$, we simply remove a point from the final block in the final round. The solution is denoted RTD$(k,n)+G($optR$(n,k))-1$ and there are $n+$optR$(n,k)$ rounds. If $1<m_{1}\leq k$ we remove the final round, and remove $m_{1}$ points of the final block of that round from the rest of the design. The solution is RTD$(k,n)+G($optR$(n,k)-1)$ and there are $n+$optR$(n,k)-1$ rounds. If $m_{1}>k$, remove the additional rounds, reinstate the cliques, and remove $m_{1}$ points from one of them. The solution is RTD$(k,n)-m_{1}$ and there are $n$ rounds. A similar argument applies when the superior allocation is constructed using MOLRs$(k,n)$, ownSG$(v,k)$ or RITD$(n,n_{2};k)$. 

\subsection{Adding points to an inferior allocation}
In some cases (when the superior design does not provide many rounds, e.g. if the MOLRs$(k,n)$ construction was used to construct it) it is preferable to add points to the inferior design instead, if possible. This involves identifying $m_{2}$ \emph{columns} of parallel blocks in the inferior allocation and adding a new point to every block in each of these columns. For some simple designs is is fairly easy to do this for one additional point (e.g. to create a $(41,5,6,9,1)$-allocation from an RBIBD$(40,4)$). However, a more general approach that can be effective for larger values of $m_{2}$, can be used when the inferior allocation 
 is an RTD that is constructed from a difference matrix (DM) or a quasi-difference matrix (QDM).

Consider the RTD$(5,14)$ constructed from a $(14,6;2)$ DM
as described in \cref{const:DM14_rounds}. The blocks can be arranged in a $14\times 14$ grid where the blocks of each row form a parallel class and the blocks in each of the first $10$ columns form a parallel class.

Since the first $10$ columns of blocks each form a parallel class, we can add up to $10$ infinite points, one to each block in one of these columns, and still maintain the fact that each row of blocks is a parallel class. Hence we can construct a $(70+m_{2},5,6,14-m_{2},m_{2})$-allocation for $1\leq m_{2}\leq 10$, with $14$ rounds in each case.

Now consider RTD$(5,15)$ which is constructed using a $(14,6;1,0;1)$-QDM  from \cite[Theorem~III.3.47]{codihandbook}.
By \cref{lem:QDM_rounds}, we can construct a $14 \times 15$ array of blocks of size $5$ where each row of blocks is a parallel class and $10$ of the columns consist of parallel blocks. By adding infinite points to the columns of parallel blocks we can construct a $(75+m_{2},5,6,15-m_{2},m_{2})$-allocation for $1\leq m_{2}\leq 10$, with $14$ rounds in each case. In fact, for $m_{2}=1$ we can obtain an extra round, by \cref{cor:extraRoundQDA}.

Since the  RTD$(5,20)$ in \cite[Lemma III.3.49]{codihandbook} is constructed using a $(19,6;1,1;1)$-QDM, by a similar argument we can construct 
$(100+m_{2},5,6,20-m_{2},m_{2})$-allocations for $1\leq m_{2}\leq 15$, with $20$ rounds when $m_{2}=1$ and $19$ rounds otherwise.

Similarly, since the RTD$(5,26)$  in \cite[Lemma III.3.53]{codihandbook}
is constructed from a $(21,6;1,0;5)$-QDM, we can construct a 
$(101,5,6,25,1)$-allocation with $22$ rounds. 

Note that in all cases, when we run out of parallel columns of blocks, we must resort to removing points from the superior allocation. 

\section{Conclusion}
We have surveyed the combinatorial structures that can be used to find optimal solutions to the Social Golfer problem (SGP) and the Social Golfer problem with adjacent group sizes (SGA). We have provided a number of new constructions and an algorithm to find optimal $(v,k)$ allocations (i.e. solutions to SGP) for any pair $(v,k)$, when $k>2$ divides $v$. We have also shown how optimal $(v,k_{1},k_{2},m_{1},m_{2})$-allocations, (i.e. optimal solutions for SGA) can be found by either adding points to an inferior allocation, or removing points from a clique associated with a superior allocation. Results for all suitable values of $v\leq 150$ and $k$ (or $k_{1}$ and $k_{2}$) are given in \cref{appendix:tables}, and the allocations can be downloaded from our BoRAT website, the link for which can be found at \cite{miller_research}.
%
%

%\backmatter
\bmsection*{Author contributions}

Original draft preparation, A.M., I.V. and R.J.R.A.; Review and editing, A.M., I.V. and R.J.R.A.; Software, A.M.,  I.V.; Data curation, I.V., Investigation, A.M., I.V. and R.J.R.A..
All authors have read and agreed to the published version of the manuscript.

%\bmsection*{Acknowledgments}

\bmsection*{Financial disclosure}
Ivaylo Valkov was supported by the EPSRC Doctoral Training Partnership award EP/N007565/1 and the UKRI funded grant EP/V026607/1. Alice Miller was funded by a Research Fellowship awarded by the Leverhume Trust. 

\bmsection*{Conflict of interest}
The authors declare that they have no potential conflict of interest.
The authors declare that they have no potential conflict of interest.
\bibliography{combinedForArxiv}

\appendix
%\bmsection{Appendix}
% A $URD(\{2,5\};30)$ with $r_{2}=5$ and $r_{5}=6$ is shown in \cref{fig:urdExample}.
% \begin{figure}
% \begin{equation*}
% \boxed{
% \begin{aligned}
% &[0,1,2,3,4],[5,6,7,8,9],[10,11,12,13,14],[15,16,17,18,19],[20,21,22,23,24],\\
% &[25,26,27,28,29]\\
% \\
% &[0,5,10,15,20],[1,6,12,18,29],[2,7,13,23,28],[3,8,19,24,27],[4,14,17,22,26],\\
% &[9,11,16,21,25]\\
% \\
% &[0,7,12,19,21],[1,5,11,22,27],[2,8,10,17,29],[3,14,15,23,25],[4,6,16,24,28],\\
% &[9,13,18,20,26]\\
% \\
% &[0,6,17,23,27],[1,9,10,19,28],[2,5,14,18,21],[3,13,16,22,29],[4,8,12,20,25],\\
% &[7,11,15,24,26]\\
% \\
% &[0,9,14,24,29],[1,8,16,23,26],[2,6,11,19,20],[3,5,12,17,28],[4,13,15,21,27],\\
% &[7,10,18,22,25]\\
% \\
% &[0,8,11,18,28],[4,5,19,23,29],[3,6,10,21,26],[7,14,16,20,27],[2,9,12,15,22],\\
% &[1,13,17,24,25]\\
% \\
% &[0,13],[1,7],[2,16],[3,9],[4,10],[5,24],[6,14],[8,15],[11,23],[12,26],[17,20],[18,27],[19,25],\\
% &[21,29],[22,28]\\
% \\
% &[0,16],[1,14],[2,25],[3,7],[4,9],[5,13],[6,15],[8,22],[10,23],[11,17],[12,27],[18,24],[19,26],\\
% &[20,29],[21,28]\\
% \\
% &[0,22],[1,15],[2,24],[3,11],[4,18],[5,26],[6,25],[7,29],[8,13],[9,23],[10,27],[12,16],[14,19],\\
% &[17,21],[20,28]\\
% \\
% &[0,25],[1,20],[2,26],[3,18],[4,11],[5,16],[6,13],[7,17],[8,21],[9,27],[10,24],[12,23],[14,28],\\
% &[15,29],[19,22]\\
% \\
% &[0,26],[1,21],[2,27],[3,20],[4,7],[5,25],[6,22],[8,14],[9,17],[10,16],[11,29],[12,24],[13,19],\\
% &[15,28],[18,23]
% \end{aligned}
% }
% \end{equation*}
% \caption{A $URD(\{2,5\};30)$ with $r_{2}=5$ and $r_{5}=6$.
% \label{fig:urdExample}
% }    
% \end{figure}

\bmsection{Starter Blocks for optimal ownSG allocations}\label{appendix:ownSG} The starter blocks for optimal ownSG$(v,k)$ allocations (see \cref{const:starterblocks}) for a range of values of $v$ and $k$, together with the number of rounds, $r$, obtained from the construction, are given in \cref{ownSG2}. %\cref{table:starterBlocks}.

\begin{table}
\centering
%\begin{footnotesize}
% \begin{tabular}{|p{10cm}|p{4.5cm}|l|}
% \begin{tabular}{|l|p{12cm}|}
\begin{tabular}{|l|l|}
\hline
$(v,k)$, $r$ & starter blocks \\ \hline
$(60,6)$, $7$ & 
$(0,1,2,3,4,5)$,
$(0,7,14,21,28,35)$,
$(0,8,13,23,33,46)$,
$(0,9,17,19,34,50)$,\\
&$(0,11,15,32,49,58)$,
$(0,16,25,41,44,57)$,
$(0,20,29,40,43,51)$\\
\hline
$(70,7)$, $7$ & 
$(0,1,2,3,4,5,6)$, 
$(0,8,16,24,32,40,48)$,
$(0,9,15,26,38,60,69)$,
$(0,10,19,25,30,43,55)$,\\
&$(0,12,27,44,53,59,64)$,
$(0,13,18,22,31,61,65)$,
$(0,17,23,34,54,57,67)$\\ \hline
$(80,8)$, $5$&
$(0,1,2,3,4,5,6,7)$, 
$(0,9,18,27,36,45,54,63)$,
$(0,10,17,29,39,43,68,78)$,\\
&$(0,11,22,25,44,55,69,74)$,
$(0,19,31,38,42,53,65,76)$\\ \hline
$(84,6)$, $9$ &
$(0,1,2,3,4,5)$, 
$(0,7,14,21,28,35)$,
$(0,8,13,23,33,46)$,\\
&$(0,9,17,19,32,58)$,
$(0,10,15,20,31,47)$,
$(0,11,16,25,44,57)$,\\
&$(0,22,26,45,71,79)$,
$(0,27,55,64,77,80)$,
$(0,37,52,63,74,83)$\\
\hline
$(84,7)$, $7$ &
$(0,1,2,3,4,5,6)$,
$(0,8,16,24,32,40,48)$,
$(0,9,15,26,38,53,76)$,
$(0,10,19,22,37,55,67)$,\\
&$(0,11,17,30,43,69,82)$,
$(0,12,41,45,51,74,78)$,
$(0,18,23,61,66,71,83)$\\
\hline
$(90,6)$, $10$ &
$(0,1,2,3,4,5)$,
$(0,7,14,21,28,35)$,
$(0,8,13,23,33,46$,
$(0,9,17,19,32,58)$,
$(0,10,15,20,31,47)$,\\
&$(0,11,16,25,44,63)$,
$(0,22,26,55,77,81)$,
$(0,27,49,64,80,89)$,
$(0,29,38,61,69,88)$,
$(0,40,59,62,75,79)$\\
\hline
$(90,9)$, $5$ &
$(0,1,2,3,4,5,6,7,8)$,
$(0,10,20,30,40,50,60,70,80)$,
$(0,11,19,33,44,52,57,68,76)$,\\
&$(0,13,26,32,37,48,61,69,74)$,
$(0,17,22,34,38,51,59,66,73)$\\
\hline
$(96,8)$, $6$ &
$(0,1,2,3,4,5,6,7)$,
$(0,9,18,27,36,45,54,63)$,
$(0,10,17,29,39,43,60,86)$,\\
&$(0,11,21,25,42,55,62,84)$,
$(0,12,31,46,50,61,67,89)$,
$(0,19,38,44,49,69,74,95)$\\
\hline
$(98,7)$, $8$ &
$(0,1,2,3,4,5,6)$,
$(0,8,16,24,32,40,48)$,
$(0,9,15,26,38,53,76)$,
$(0,10,19,22,37,55,67)$,\\
&$(0,11,17,23,36,54,83)$,
$(0,12,18,27,29,51,66)$,
$(0,13,47,58,81,85,94)$,
$(0,25,30,69,82,87,92)$
\\
\hline
$(105,7)$, $9$ &
$(0,1,2,3,4,5,6)$,
$(0,8,16,24,32,40,48)$,
$(0,9,15,26,38,53,76)$,\\
&$(0,10,19,22,37,55,67)$,
$(0,11,17,23,36,54,83)$,
$(0,12,18,27,29,51,73)$,\\
&$(0,20,30,45,50,81,103)$,
$(0,25,44,71,90,96,101)$,
$(0,43,52,69,72,82,95)$
\\
\hline
$(112,8)$, $8$ &
$(0,1,2,3,4,5,6,7)$,
$(0,9,18,27,36,45,54,63)$,
$(0,10,17,29,39,43,60,86)$,
$(0,11,21,25,42,55,62,84)$,\\
&$(0,12,19,26,41,69,94,111)$,
$(0,14,28,34,59,71,93,105)$,
$(0,15,20,30,33,53,58,107)$,
$(0,31,50,61,76,81,91,110)$
\\
\hline
$(120,6)$, $13$ &
$(0,1,2,3,4,5)$,
$(0,7,14,21,28,35)$,
$(0,8,13,23,33,46)$,
$(0,9,17,19,32,58)$,
$(0,10,15,20,31,47)$,\\
&$(0,11,16,25,44,57)$,
$(0,22,26,45,49,71)$,
$(0,27,38,55,70,89)$,
$(0,29,37,63,80,88)$,\\
&$(0,34,73,81,110,119)$,
$(0,39,59,76,91,116)$,
$(0,40,53,56,93,115)$,
$(0,41,50,79,106,117)$
\\
\hline
$(126,7)$, $10$ &
$(0,1,2,3,4,5,6)$,
$(0,8,16,24,32,40,48)$,
$(0,9,15,26,38,53,76)$,
$(0,10,19,22,37,55,67)$,\\
&$(0,11,17,23,36,54,83)$,
$(0,12,18,27,29,51,66)$,
$(0,13,30,39,57,96,122)$,\\
&$(0,20,25,45,65,75,120)$,
$(0,31,41,44,64,89,123)$,
$(0,34,47,80,85,102,121)$\\
\hline
$(126,9)$, $7$ &
$(0,1,2,3,4,5,6,7,8)$,
$(0,10,20,30,40,50,60,70,80)$,
$(0,11,19,32,43,48,67,96,125)$,
$(0,12,23,28,44,47,61,69,94)$,\\
&$(0,13,24,29,41,88,107,111,118)$,
$(0,17,52,55,78,86,101,112,120)$,
$(0,34,38,51,62,77,84,109,121)$\\
\hline
$(132,6)$, $14$ &
$(0,1,2,3,4,5)$,
$(0,7,14,21,28,35)$,
$(0,8,13,23,33,46)$,
$(0,9,17,19,32,58)$,
$(0,10,15,20,31,47)$,\\
&$(0,11,16,25,44,57)$,
$(0,22,26,45,49,71)$,
$(0,27,38,55,70,89)$,
$(0,29,37,63,80,88)$,
$(0,34,50,59,79,87)$,\\
&$(0,39,77,86,106,121)$,
$(0,40,53,69,103,128)$,
$(0,41,62,93,112,115)$,
$(0,43,52,74,107,129)$\\
% \hline
% \end{tabular}
% \caption{Starter blocks and for optimal ownSG allocations 1\label{ownSG1}}
% \end{table}

% \begin{table}
% \centering
% %\begin{footnotesize}
% % \begin{tabular}{|p{10cm}|p{4.5cm}|l|}
% \begin{tabular}{|l|p{12cm}|}
% \hline
% $(v,k)$, $r$ & starter blocks \\ \hline
\hline
$(135,9)$, $7$ &
$(0,1,2,3,4,5,6,7,8)$,
$(0,10,20,30,40,50,60,70,80)$,
$(0,11,19,32,43,48,67,96,125)$,
$(0,12,23,28,44,47,61,69,94)$,\\
&$(0,13,21,29,52,73,104,116,132)$,
$(0,14,22,37,79,92,111,123,134)$,
$(0,35,39,76,82,97,114,119,131)$\\
\hline
$(156,6)$, $15$ &
$(0,1,2,3,4,5)$,
$(0,7,14,21,28,35)$,
$(0,8,13,23,33,46)$,
$(0,9,17,19,32,58)$,
$(0,10,15,20,31,47)$,\\
&$(0,11,16,25,44,57)$,
$(0,22,26,45,49,71)$,
$(0,27,38,55,70,89)$,
$(0,29,37,63,80,88)$,
$(0,34,50,59,79,87)$,\\
&$(0,39,61,86,106,119)$,
$(0,40,43,74,105,143)$,
$(0,41,52,67,104,147)$,
$(0,53,56,97,124,135)$,
$(0,62,109,118,149,153)$\\
\hline
\end{tabular}
\caption{Starter blocks and for optimal ownSG allocations\label{ownSG2}}
\end{table}

% \nocite{*}% Show all bib entries - both cited and uncited; comment this line to view only cited bib entries;

\begin{comment}
\bmsection*{Author Biography}

\begin{biography}{\includegraphics[width=76pt,height=76pt,draft]{Template/empty}}{
{\textbf{Author Name.} Please check with the journal's author guidelines whether
author biographies are required. They are usually only included for
review-type articles, and typically require photos and brief
biographies for each author.}}
\end{biography}
\end{comment}

\bmsection{Tables of results}\label{appendix:tables}
% \section{Introduction}\label{sect:intro}
This appendix contains optimal solutions to the Social Golfer problem (SGP) and the Social Golfer problem with adjacent group sizes (SGA).
%a draft of which is available at ADD %ARXIV REFERENCE.
An optimal solution is one for which no better solution (i.e. no solution with more rounds) is known. 

Each of \crefrange{table:parallel1}{table:parallel38} contains solutions for an interval of values $v$ (number of players/points), for $v\leq 150$. For each value $v$ the possible sets of block sizes $K$ are considered. If $|K|=1$, $K$ consists of a single divisor of $v$, otherwise (except for \cref{table:parallel1}), $K=\{4,5\}$ or $K=\{5,6\}$. When $|K|=1$, $K=\{k_1\}$, $m_{1}=v/k_{1}$ denotes the number of blocks of size $k_{1}$ in each round, and $m_{2}$ is $0$. When $K=\{k_{1}, k_{2}\}$, $m_{i}$ is the number of blocks of size $k_{i}$ in each round, for $i\in \{1,2\}$. $MAX$ denotes the maximum possible number of rounds and $r$ denotes the number of rounds achieved using the solution provided. \cref{table:additional} contains some additional solutions for $v>150$. Each table caption indicates the sets $K$ that are relevant for that range of $v$.

\crefrange{table:summary1}{table:summary2} include descriptions of the solutions provided. 
Note that in \cref{table:summary1}, $N(v)$ denotes the (current lower bound on the) number of MOLS of order $v$ and $n=v/k$. 

In the \emph{Description} columns of \crefrange{table:summary1}{table:summary2}, KTS, NKTS, RBIBD, RTD and RGDD are abbreviations for: Kirkman triple System, Nearly Kirkman Triple System, resolvable balanced incomplete block design, resolvable transversal design and resolvable group divisible design respectively. Similarly, URD, RITD, DM and QDM stand for: uniformly resolvable pairwise balanced design, resolvable incomplete transversal design, difference matrix, and quasi-difference matrix. 

In the \emph{Solution} columns, designs D$(v,k)$, D$^{\prime}(v,k)$ and E$(v,k)$ represent any element of sets $\mathcal{D}$,  $\mathcal{D}^{\prime}$ and $\mathcal{E}$ where: 
$\mathcal{D} = \{$RTD$(k,n)$, RGDD$(v,k,g)$, 
URD$(v,k,k_{1})$, RITD$(n_{1},n_{2};k)$, MOLRs$(k,n)$, ownSG$(v,k)\}$; 
$\mathcal{D}^{\prime}=\{$D$(v,k)+G(t) :\;$D$(v,k)\in \mathcal{D}\}$, where $G(t)$ denotes that $t$ rounds are added from the groups of the design associated with D$(v,k)$; and 
$\mathcal{E} = \mathcal{D} \cup \mathcal{D}^{\prime}$.

\begin{table}
\centering
%\resizebox{\textwidth}{!}{%
%
\caption{\label{table:summary1} Key to solutions in \crefrange{table:parallel1}{table:parallel38} part 1.}
\end{table}

\begin{table}
\centering
%\resizebox{\textwidth}{!}{%
% %
\caption{\label{table:summary2} Key to solutions in \crefrange{table:parallel1}{table:parallel38} part 2.}
\end{table}

\begin{comment}
Other terms for captialization (to be commented out):
\begin{enumerate}
    \item balanced incomplete block design. 
    \item transversal design
    \item group divisible design
    \item incomplete transversal design
\end{enumerate}
\end{comment}

%removed as don't seem to need any more
%SG$(v,k)$ & $r$ rounds of blocks of size $k$ from an %existing, published $r$-round SGP solution.\\

\begin{table}
\centering
\begin{tabular}{|c|c|c|c|c|c|}
\toprule
    $v$ & $K$ & $m_{1},m_{2}$ & $MAX$ & solution & $r$  \\
    \midrule 
    $12$ & $3$ &$4,0$&$4$&
    RTD$(3,4)$&$4$\\
    \hline
    $13$ & $3,4$ &$3,1$ & $5$ & RBIBD$(16,4) - 3,B$ & $4$ \\
    \hline
    $14$ & $3,4$ & $2,2$ & $5$ & RBIBD$(16,4) - 2,B$  & $4$\\
    \hline
    $15$ & $3$ & $5,0$ & $7$ & KTS$(15)$&$7$\\
         & $3,4$ & $1,3$ & $5$& RBIBD$(16,4) - 1$ & $5$\\
    \hline
    $16$ &$4$ & $4,0$ & $5$ &RBIBD$(16,4)$&$5$\\
         &$3,4$ & $4,1$ &$6$&RTD$(4,5) - 4 $&$5$\\
    \hline
    $17$ & $3,4$ & $3,2$ &$6$&RTD$(4,5) - 3 $&$5$\\
    \hline
    $18$ & $3$ & $6,0$ & $8$ &NKTS$(18,3)$ &$8$\\
         & $3,4$ & $2,3$ & $6$ & RTD$(4,5) - 2 $&$5$\\
    \hline
    $19$& $3,4$ & $1,4$& $6$ & RTD$(4,5) - 1$ &$5$\\
        &       & $5,1$& $8$&NKTS$(18)+1$&$6$\\
    \bottomrule
\end{tabular}
\caption{\label{table:parallel1} SGP and SGA solutions for  $12\leq v\leq 19$. 
\\Block sizes $K=\{3\}, \{4\}$ or $\{3,4\}$.}
\end{table}

% \newpage

\begin{table}[p]
\centering
\begin{tabular}{|c|c|c|c|c|c|c|}
\toprule
    $v$ & $K$ & $m_{1},m_{2}$ & $MAX$ & solution & $r$  \\
    \midrule
    $20$ & $4$ &$5,0$ & $5$ & RTD$(4,5)$& $5$\\
    \hline
    $21$ & $3$ & $7,0$ & $10$&KTS$(21)$&$10$\\
         &$4,5$&$4,1$&$6$ & RBIBD$(25,5)-4,B$& $5$\\
   \hline
    $22$  & $4,5$ & $3,2$ & $6$ & RBIBD$(25,5)-3,B$& $5$\\
    \hline
    $23$ & $4,5$ & $2,3$ & $6$ &RBIBD$(25,5)-2,B$& $5$\\
    \hline
    $24$ & $3$ & $8,0$ & $11$ &NKTS$(24)$ &$11$\\
         & $4$ & $6,0$ & $7$ &RGDD$(24,4,3)$ & $7$\\
\bottomrule
\end{tabular}
\caption{\label{table:parallel2}  SGP and SGA solutions for  $20\leq v\leq 34$. \\Block sizes  $K=\{3\}, \{4\},\{5\},\{4,5\}$, or $\{5,6\}$.}

\end{table}

% \newpage

\begin{table}[p]
\centering
\begin{tabular}{|c|c|c|c|c|c|c|}
\toprule
    $v$ & $K$ & $m_{1},m_{2}$ & $MAX$ & solution & $r$  \\
    \midrule
    $35$ & $5$ & $7,0$ & $8$ &RTD$(5,7)$& $7$\\
         & $4,5$ & $5,3$ & $9$ &RTD$(5,8)-5$& $8$\\
         & $5,6$ & $1,5$ & $7$ &MOLRs$(6,6) +G(1) - 1$ & $3$\\
    \hline
    $36$ & $3$ & $12,0$ & $17$ &NKTS$(36)$ &$17$\\
         & $4$ & $9,0$ & $11$ &RGDD$(36,4,3)$ & $11$\\
         & $6$ & $6,0$ & $3$ &MOLRs$(6,6)+G(1)$&$3$\\
         &$4,5$ & $4,4$ & $9$ &RTD$(5,8)-4$& $8$\\
\bottomrule
\end{tabular}
\caption{\label{table:parallel3} SGP and SGA solutions for  $35\leq v\leq 43$. \\Block sizes $K=\{3\}, \{4\}, \{5\}, \{6\},\{4,5\}$ or $\{5,6\}$}
\end{table}

\begin{table}[p]
\centering
\begin{tabular}{|c|c|c|c|c|c|c|}
\toprule
    $v$ & $K$ & $m_{1},m_{2}$ & $MAX$ & solution& $r$  \\
    \midrule
    $44$ &$4$ & $11,0$ & $14$ &RGDD$(44,4,2)$& $14$\\
         &$4,5$ & $1,8$ & $11$ &RTD$(5,9)-1$& $9$\\
         &      & $6,4$ & $12$ &RITD$(10,2;5)-6$ & $8$\\
         & $5,6$ & $4,4$ & $9$ &RTD$(6,8)-4$& $8$\\
    \hline
    $45$& $3$ & $15,0$ & $22$ &KTS$(45)$&$22$\\
        & $5$ & $9,0$ & $11$ &RTD$(5,9)$& $9$\\
        & $4,5$ & $5,5$ & $12$ &RITD$(10,2;5)-5$ & $8$\\
        &       & $10,1$ & $14$ &RTD$(5,11)-10$& $11$\\
        &$5,6$ & $3,5$ & $9$ &RTD$(6,8)-3$& $8$\\
   \hline
    $46$ & $4,5$&$4,6$&$12$&RITD$(10,2;5)-4$ & $8$\\
         &      &$9,2$&$13$&RTD$(5,11)-9$& $11$\\
         &$5,6$ &$2,6$ & $9$ &RTD$(6,8)-2$& $8$\\
         &      &$8,1$&$10$ &RTD$(6,9)-8$& $9$\\
    \hline
    $47$ &$4,5$ & $3,7$ & $12$ &RITD$(10,2;5)-2$ & $8$\\
         &      & $8,3$ & $13$ &RTD$(5,11)-8$& $11$\\
         &$5,6$ & $1,7$ & $9$ &RTD$(6,8)-1$& $8$\\
         &      & $7,2$ & $10$ &RTD$(6,9)-7$& $9$\\
    \hline
    $48$& $3$ & $16,0$ & $23$ & NKTS$(48)$ &$23$ \\
        & $4$ & $12,0$ & $15$ &RGDD$(44,4,3)$& $15$\\
        & $6$ & $8,0$ & $9$ &RTD$(6,8)$& $8$\\
        & $4,5$&$2,8$&$12$ &RITD$(10,2;5)-2$ & $8$\\
        &      &$7,4$&$13$ &RTD$(5,11)-7$& $11$\\
        &$5,6$ &$6,3$ &$10$ &RTD$(6,9)-6$& $9$\\
    \hline
    $49$ & $7$ & $7,0$ & $8$ &RBIBD$(49,7)$&$8$\\
         &$4,5$&$1,9$ & $12$ &RITD$(10,2;5)+G(1)-1$ & $9$\\
         &     &$6,5$ & $13$ &RTD$(5,11)-6$& $11$\\
         &     &$11,1$ & $15$ &RTD$(5,12)-11$& $12$\\
         &$5,6$ &$5,4$&$10$&RTD$(6,9)-5$& $9$\\
    \hline
    $50$&$5$&$10,0$&$12$&RITD$(10,2;5)+G(1)$ & $9$\\
        &$4,5$&$5,6$&$13$&RTD$(5,11)-5$& $11$\\
        &     &$10,2$&$15$&RTD$(5,12)-10$& $12$\\
        &$5,6$&$4,5$&$10$&RTD$(6,9)-4$& $9$\\
\bottomrule
\end{tabular}
\caption{\label{table:parallel4} SGP and SGA solutions for  $44\leq v\leq 50$. \\Block sizes $K=\{3\}, \{4\}, \{5\}, \{6\}, \{7\},\{4,5\}$ or $\{5,6\}$}
\end{table}

\begin{table}[p]
\centering
\begin{tabular}{|c|c|c|c|c|c|}
\hline
    $v$ & $K$ & $m_{1},m_{2}$ & $MAX$ & solution& $NR$  \\
    \hline
    $51$ & $3$&$17,0$&$25$&KTS$(51)$&$25$\\
         & $4,5$ & $4,7$ & $13$ &RTD$(5,11)-4$&$11$\\ 
         &       & $9,3$ & $15$ &RTD$(5,12)-9$&$12$\\
         &$5,6$ & $3,6$ & $10$ &RTD$(6,9)-3$&$9$\\
         &      & $9,1$ & $12$ &ownSG$(60,6)-9$&$7$\\
    \hline
    $52$ & $4$ & $13$ &$17$ &RBIBD$(52,4)$&$17$ \\
         & $4,5$ & $3,8$ &$13$ &RTD$(5,11)-3$&$11$\\ 
         &       & $8,4$ &$15$ &RTD$(5,12)-8$&$12$\\
         &$5,6$  & $2,7$ &$10$ &RTD$(6,9)-2$&$9$\\
         &       & $8,2$ & $12$ &ownSG$(60,6)-8$&$7$\\
    \hline
    $53$ & $4,5$ & $2,9$ &$13$ &RTD$(5,11)-2$&$11$\\  
         &       & $7,5$ &$14$ &RTD$(5,12)-7$&$12$\\
         &       &$12,1$ & $16$ &RTD$(5,13)-12$&$13$\\
         &$5,6$ & $1,8$ & $10$ &RTD$(6,9)-1$&$9$\\
         &      & $7,3$ & $11$ &ownSG$(60,6)-7$&$7$\\
    \hline
    $54$ & $3$ & $18,0$ & $26$ &NKTS$(54)$&$26$\\
         & $6$ & $9,0$  & $10$ &RTD$(6,9)$&$9$\\ 
         & $4,5$ & $1,10$ &$13$ &RTD$(5,11)-1$&$11$\\ 
         &       & $6,6$  &$14$ &RTD$(5,12)-6$&$12$\\
         &       &$11,2$  &$16$ &RTD$(5,13)-11$&$13$\\
         &$5,6$ & $6,4$ & $11$ &ownSG$(60,6)-6$&$7$\\
    \hline
    $55$ & $5$ &$11,0$ & $13$ &RTD$(5,11)$&$11$\\
         & $4,5$ & $5,7$ &$14$ &RTD$(5,12)-5$&$12$\\
         &       & $10,3$  &$16$ &RTD$(5,13)-10$&$13$\\
         &$5,6$ & $5,5$ & $11$ &ownSG$(60,6)-5$&$7$\\
    \hline
\end{tabular}
\caption{\label{table:parallel5}  SGP and SGA solutions for  $51\leq v\leq 55$. \\Block sizes  $K=\{3\}, \{4\},\{5\},\{6\},\{4,5\}$, or $\{5,6\}$.}
\end{table}

\begin{table}[p]
\centering
\begin{tabular}{|c|c|c|c|c|c|}
\hline
    $v$ & $K$ & $m_{1},m_{2}$ & $MAX$ & solution& $NR$  \\
    \hline
     $56$ & $4$ &$14,0$ & $18$ &RGDD$(56,4,2)$&$18$\\
          & $7$ &$8,0$ & $9$ &RTD$(7,8)$&$8$\\
          & $4,5$ & $4,8$ &$14$ &RTD$(5,12)-4$&$12$\\
          &       & $9,4$ &$16$ &RTD$(5,13)-9$&$13$\\
          &$5,6$ & $4,6$ & $11$ &ownSG$(60,6)-4$&$7$\\
          &      & $10,1$ & $13$ &RTD$(6,11)-10$&$11$\\
    \hline
    $57$ &  $3$ &$19,0$ & $28$ &KTS$(57)$ &$28$\\ 
           & $4,5$ & $3,9$ &$14$ &RTD$(5,12)-3$&$12$\\
           &       & $8,5$ &$16$ &RTD$(5,13)-8$&$13$\\
           &       & $13,1$ &$18$ &RTD$(5,14)-13$&$14$\\
           &$5,6$ & $3,7$ & $11$ &ownSG$(60,6)-3$&$7$\\
           &      & $9,2$ & $13$ &RTD$(6,11)-9$&$11$\\
    \hline
    $58$ & $4,5$ & $2,10$ &$14$ &RTD$(5,12)-2$&$12$\\
           &       & $7,6$ &$16$ &RTD$(5,13)-7$&$13$\\
           &       & $12,2$ &$17$ &RTD$(5,14)-12$&$14$\\
           &$5,6$ & $2,8$ & $11$ &ownSG$(60,6)-2$&$7$\\
           &      & $8,3$ & $13$ &RTD$(6,11)-8$&$11$\\
    \hline
    $59$ & $4,5$ & $1,11$ &$14$ &RTD$(5,12)-1$&$12$\\
           &       & $6,7$ &$16$ &RTD$(5,13)-6$&$13$\\
           &       & $11,3$ &$17$ &RTD$(5,14)-11$&$14$\\
           &$5,6$ & $1,9$ & $11$ &ownSG$(60,6)-1$&$7$\\
           &      & $7,4$ & $13$ &RTD$(6,11)-7$&$11$\\
    \hline
    $60$ & $3$ & $20,0$ & $29$ &NKTS$(60)$&$29$\\
         & $4$ & $15,0$ & $19$ &RGDD$(60,4,3)$&$19$\\
         & $5$ & $12,0$ & $14$ &RTD$(5,12)$&$12$\\
         & $6$ & $10,0$ & $11$ &ownSG$(60,6)$&$7$\\
         & $4,5$ & $5,8$ &$16$ &RBIBD$(65,5)-5,B$&$15$\\
         &       &$10,4$ &$17$ &RTD$(5,14)-10$&$14$\\
         &$5,6$ & $6,5$ & $13$ &RTD$(6,11)-6$&$11$\\
    \hline
\end{tabular}
\caption{\label{table:parallel6} SGP and SGA solutions for  $56\leq v\leq 60$. \\Block sizes  $K=\{3\}, \{4\},\{5\},\{6\},\{7\},\{4,5\}$, or $\{5,6\}$.}
\end{table}

\begin{table}[p]
\centering
\begin{tabular}{|c|c|c|c|c|c|}
\hline
    $v$ & $K$ & $m_{1},m_{2}$ & $MAX$ & solution& $NR$  \\
    \hline
     $61$ & $4,5$ & $4,9$ &$16$ &RBIBD$(65,5)-4,B$&$15$\\
          &       & $9,5$ &$17$ &RTD$(5,14)-9$&$14$\\
          &       & $14,1$ &$19$ &RTD$(5,15)-14$&$15$\\
          &$5,6$ & $5,6$ & $13$ &RTD$(6,11)-5$&$11$\\
          &      & $11,1$ & $14$ &RTD$(6,12)-11$&$12$\\
    \hline
    $62$ & $4,5$ & $3,10$ &$16$ &RBIBD$(65,5)-3,B$&$15$\\
          &       & $8,6$ &$17$ &RTD$(5,14)-8$&$14$\\
          &       & $13,2$ &$19$ &RTD$(5,15)-13$&$15$\\
          &$5,6$ & $4,7$ & $13$ &RTD$(6,11)-4$&$11$\\
          &      & $10,2$ & $14$ &RTD$(6,12)-10$&$12$\\
    \hline
    $63$  & $3$ & $21,0$ & $31$ &KTS$(63)$&$31$\\
          & $7$ & $9,0$ &$10$ &RTD$(7,9)$&$9$\\
          & $4,5$ & $2,11$ &$16$ &RBIBD$(65,5)-2,B$&$15$\\
          &       & $7,7$ &$17$ &RTD$(5,14)-7$&$14$\\
          &       & $12,3$ &$19$ &RTD$(5,15)-12$&$15$\\
          &$5,6$ & $3,8$ & $13$ &RTD$(6,11)-3$&$11$\\
          &      & $9,3$ & $14$ &RTD$(6,12)-9$&$12$\\
    \hline
    $64$ & $4$ & $16,0$ &$21$ &RBIBD$(64,4)$&$21$\\
          & $8$ & $8,0$ &$9$ &RBIBD$(64,8)$&$9$\\
          & $4,5$ & $1,12$ &$16$ &RBIBD$(65,5)-1$&$16$\\
          &       & $6,8$ &$17$ &RTD$(5,14)-6$&$14$\\
          &       & $11,4$ &$19$ &RTD$(5,15)-11$&$15$\\
          &$5,6$ & $2,9$ & $13$ &RTD$(6,11)-2$&$11$\\
          &      & $8,4$ & $14$ &RTD$(6,12)-8$&$12$\\
    \hline
    $65$ & $5$ & $13,0$ &$16$ &RBIBD$(65,5)$&$16$\\
          & $4,5$ & $5,9$ &$17$ &RTD$(5,14)-5$&$14$\\
          &       & $10,5$ &$18$ &RTD$(5,15)-10$&$15$\\
          &       & $15,1$ &$20$ &RTD$(5,16)-15$&$16$\\
          &$5,6$ & $1,10$ & $13$ &RTD$(6,11)-1$&$11$\\
          &      & $7,5$ & $14$ &RTD$(6,12)-7$&$12$\\
    \hline
\end{tabular}
\caption{\label{table:parallel7} SGP and SGA solutions for  $61\leq v\leq 65$. \\Block sizes  $K=\{3\}, \{4\},\{5\},\{7\},\{8\},\{4,5\}$, or $\{5,6\}$.}
\end{table}

\begin{table}[p]
\centering
\begin{tabular}{|c|c|c|c|c|c|}
\hline
    $v$ & $K$ & $m_{1},m_{2}$ & $MAX$ & solution& $NR$  \\
    \hline
     $66$ & $3$ & $22,0$ &$32$ &NKTS$(66)$&$32$\\
          & $6$ & $11,0$ &$13$ &RTD$(6,11)$&$11$\\
          & $4,5$ & $4,10$ &$17$ &RTD$(5,14)-4$&$14$\\
          &       & $9,6$  &$18$ &RTD$(5,15)-9$&$15$\\
          &       & $14,2$ &$20$ &RTD$(5,16)-14$&$16$\\
          & $5,6$ & $6,6$ & $14$ &RTD$(6,12)-6$&$12$\\
          &       & $12,1$ & $15$ &RTD$(6,13)-12$&$13$\\
    \hline
    $67$ & $4,5$ & $3,11$ &$17$ &RTD$(5,14)-3$&$14$\\
          &       & $8,7$  &$18$ &RTD$(5,15)-8$&$15$\\
          &       & $13,3$ &$20$ &RTD$(5,16)-13$&$16$\\
          & $5,6$ & $5,7$ & $14$ &RTD$(6,12)-5$&$12$\\
          &       & $11,2$ & $15$ &RTD$(6,13)-11$&$13$\\
    \hline
     $68$ & $4$ & $17,0$ &$22$ & RGDD$(68,4,2)$&$22$\\
          & $4,5$ & $2,12$ &$17$ &RTD$(5,14)-2$&$14$\\
          &       & $7,8$  &$18$ &RTD$(5,15)-7$&$15$\\
          &       & $12,4$ &$20$ &RTD$(5,16)-12$&$16$\\
          & $5,6$ & $4,8$ & $14$ &RTD$(6,12)-4$&$12$\\
          &       & $10,3$ & $15$ &RTD$(6,13)-10$&$13$\\
    \hline    
    $69$ & $3$ & $23,0$ &$34$ &KTS$(69)$&$34$\\
          & $4,5$ & $1,13$ &$17$ &RTD$(5,14)-1$&$14$\\
          &       & $6,9$  &$18$ &RTD$(5,15)-6$&$15$\\
          &       & $11,5$ &$20$ &RTD$(5,16)-11$&$16$\\
          &       & $16,1$ &$22$ &RTD$(5,17)-16$&$17$\\
          & $5,6$ & $3,9$ & $14$ &RTD$(6,12)-3$&$12$\\
          &       & $9,4$ & $15$ &RTD$(6,13)-9$&$13$\\
    \hline
     $70$ & $5$ & $14,0$ &$17$ &RTD$(5,14)$&$14$\\
          & $7$ & $10,0$ &$11$ &ownSG$(70,7)$&$7$\\
          & $4,5$ & $5,10$ &$18$ &RTD$(5,15)-5$&$15$\\
          &       & $10,6$  &$20$ &RTD$(5,16)-10$&$16$\\
          &       & $15,2$ &$21$ &RTD$(5,17)-15$&$17$\\
          & $5,6$ & $2,10$ & $14$ &RTD$(6,12)-2$&$12$\\
          &       & $8,5$ & $15$ &RTD$(6,13)-8$&$13$\\
\hline
\end{tabular}
\caption{\label{table:parallel8} SGP and SGA solutions for  $66\leq v\leq 70$. \\Block sizes  $K=\{3\}, \{4\},\{5\},\{6\},\{7\},\{4,5\}$, or $\{5,6\}$.}
\end{table}

\begin{table}[p]
\centering
\begin{tabular}{|c|c|c|c|c|c|}
\hline
    $v$ & $K$ & $m_{1},m_{2}$ & $MAX$ & solution& $NR$  \\
    \hline
     $71$ & $4,5$ & $4,11$ &$18$ &RTD$(5,15)-4$&$15$\\
          &       & $9,7$  &$20$ &RTD$(5,16)-9$&$16$\\
          &       & $14,3$ &$21$ &RTD$(5,17)-14$&$17$\\
          & $5,6$ & $1,11$ & $14$ &RTD$(6,12)+G(1) -1$&$13$\\
          &       & $7,6$ & $15$ &RTD$(6,13)-7$&$13$\\
          &       & $13,1$ & $17$ &RTD$(5,14)\_D+1$&$14$\\
    \hline
    $72$ & $3$ & $24,0$ &$35$ &NKTS$(72)$& $35$\\
\hline
\end{tabular}
\caption{\label{table:parallel9} SGP and SGA solutions for  $71\leq v\leq 75$. \\Block sizes  $K=\{3\}, \{4\},\{5\},\{6\},\{8\},\{4,5\}$, or $\{5,6\}$.}
\end{table}

\begin{table}[p]
\centering
\begin{tabular}{|c|c|c|c|c|c|}
\hline
    $v$ & $K$ & $m_{1},m_{2}$ & $MAX$ & solution& $NR$  \\
    \hline
     $76$ & $4$ & $19,0$ &$25$ &RBIBD$(76,4)$&$25$\\
          & $4,5$ & $4,12$ &$19$ &RTD$(5,16)-4$&$16$\\
          &       & $9,8$  &$21$ &RTD$(5,17)-9$&$17$\\
          &       & $14,4$ &$22$ &RTD$(5,18)-14$&$18$\\
          & $5,6$ & $2,11$ & $15$ &RTD$(6,13)-2$&$13$\\
          &       & $8,6$ & $16$ &RTD$(5,14)\_D+6$&$14$\\
          &       & $14,1$ & $18$ &RTD$(5,15)\_Q(14,1)+1$&$15$\\
    \hline
    $77$ & $7$ & $11,0$ &$12$ &RTD$(7,11)$&$11$\\
         & $4,5$ & $3,13$ &$19$ &RTD$(5,16)-3$&$16$\\
         &       & $8,9$  &$21$ &RTD$(5,17)-8$&$17$\\
         &       & $13,5$ &$22$ &RTD$(5,18)-13$&$18$\\
         &       & $18,1$ &$24$ &RTD$(5,19)-18$&$19$\\
         & $5,6$ & $1,12$ & $15$ &RTD$(6,13)-1$&$13$\\
         &       & $7,7$ & $16$ &RTD$(5,14)\_D+7$&$14$\\
         &       & $13,2$ & $18$ &RTD$(5,15)\_Q(14,1)+2$&$14$\\
    \hline
     $78$ & $3$ & $26,0$ &$38$ &NKTS$(78)$&$38$\\
          & $6$ & $13,0$ &$15$ &RTD$(6,13)$&$13$\\
          & $4,5$ & $2,14$ &$19$ &RTD$(5,16)-2$&$16$\\
          &       & $7,10$  &$21$ &RTD$(5,17)-7$&$17$\\
          &       & $12,6$ &$22$ &RTD$(5,18)-12$&$18$\\
          &       & $17,2$ &$24$ &RTD$(5,19)-17$&$19$\\
          & $5,6$ & $6,8$ & $16$ &RTD$(5,14)\_D+8$&$14$\\
          &       & $12,3$ & $18$ &RTD$(5,15)\_Q(14,1)+3$&$14$\\
    \hline    
    $79$ & $4,5$ & $1,15$ &$19$ &RTD$(5,16)-1$&$16$\\
          &       & $6,11$  &$21$ &RTD$(5,17)-6$&$17$\\
          &       & $11,7$ &$22$ &RTD$(5,18)-11$&$18$\\
          &       & $16,3$ &$24$ &RTD$(5,19)-16$&$19$\\
          & $5,6$ & $5,9$ & $16$ &RTD$(5,14)\_D+9$&$14$\\
          &       & $11,4$ & $18$ &RTD$(5,15)\_Q(14,1)+4$&$14$\\
\hline
\end{tabular}
\caption{\label{table:parallel10} SGP and SGA solutions for  $76\leq v\leq 79$. \\Block sizes  $K=\{3\}, \{4\},\{6\},\{7\},\{4,5\}$, or $\{5,6\}$.}
\end{table}

\begin{table}[p]
\centering
\begin{tabular}{|c|c|c|c|c|c|}
\hline
    $v$ & $K$ & $m_{1},m_{2}$ & $MAX$ & solution& $NR$  \\
    \hline
    $80$ & $4$ & $20,0$ &$26$ &RGDD$(20,4,2)$&$26$\\
          & $5$ & $16,0$ &$19$ &RTD$(5,16)$&$16$\\
          & $8$ & $10,0$ &$11$ &ownSG$(80,8)$&$5$\\
          & $4,5$ & $5,12$ &$21$ &RBIBD$(85,5)-5,B$&$20$\\
          &       & $10,8$  &$22$ &RTD$(5,18)-10$&$18$\\
          &       & $15,4$ &$24$ &RTD$(5,19)-15$&$19$\\
          & $5,6$ & $4,10$ & $16$ &RTD$(5,14)\_D+10$&$14$\\
          &       & $10,5$ & $18$ &RTD$(5,15)\_Q(14,1)+5$&$14$\\
    \hline
     $81$ & $3$ & $27,0$ &$40$ &KTS$(81)$&$40$\\
          & $9$ & $9,0$ &$10$ & RBIBD$(81,9)$&$10$\\
          & $4,5$ & $4,13$ &$21$ &RBIBD$(85,5)-4,B$&$20$\\
          &       & $9,9$  &$22$ &RTD$(5,18)-9$&$18$\\
          &       & $14,5$ &$24$ &RTD$(5,19)-14$&$19$\\
          &       & $19,1$ &$26$ &RTD$(5,20)-19$&$20$\\
          & $5,6$ & $3,11$ & $16$ &ownSG$(84,6)-3$&$9$\\
          &       & $9,6$ & $18$ &RTD$(5,15)\_Q(14,1)+6$&$14$\\
          &       & $15,1$ & $19$ &RTD$(6,16)-15$&$16$\\
    \hline
    $82$ & $4,5$ & $3,14$ &$21$ &RBIBD$(85,5)-3,B$&$20$\\
         &       & $8,10$  &$22$ &RTD$(5,18)-8$&$18$\\
         &       & $13,6$ &$24$ &RTD$(5,19)-13$&$19$\\
         &       & $18,2$ &$25$ &RTD$(5,20)-18$&$20$\\
         & $5,6$ & $2,12$ & $16$ &ownSG$(84,6)-2$&$9$\\
         &       & $8,7$ & $17$ &RTD$(5,15)\_Q(14,1)+7$&$14$\\
         &       & $14,2$ & $19$ &RTD$(6,16)-14$&$16$\\
    \hline
     $83$ & $4,5$ & $2,15$ &$21$ &RBIBD$(85,5)-2,B$&$20$\\
          &       & $7,11$  &$22$ &RTD$(5,18)-7$&$18$\\
          &       & $12,7$ &$23$ &RTD$(5,19)-12$&$19$\\
          &       & $17,3$ &$25$ &RTD$(5,20)-17$&$20$\\
          & $5,6$ & $1,13$ & $16$ &ownSG$(84,6)-1$&$9$\\
          &       & $7,8$ & $17$ &RTD$(5,15)\_Q(14,1)+8$&$14$\\
          &       & $13,3$ & $19$ &RTD$(6,16)-13$&$16$\\
    \hline    
\end{tabular}
\caption{\label{table:parallel11} SGP and SGA solutions for  $80\leq v\leq 83$. \\Block sizes  $K=\{3\}, \{4\},\{5\},\{8\},\{9\},\{4,5\}$, or $\{5,6\}$.}
\end{table}

\begin{table}[p]
\centering
\begin{tabular}{|c|c|c|c|c|c|}
\hline
    $v$ & $K$ & $m_{1},m_{2}$ & $MAX$ & solution& $NR$  \\
    \hline
    $84$  & $3$ & $28,0$ &$41$ &NKTS$(84)$&$41$\\
          & $4$ & $21,0$ &$27$ &RGDD$(84,4,3)$& $27$\\
          & $6$ & $14,0$ &$16$ &ownSG$(84,6)$&$9$\\
          & $7$ & $12,0$ &$13$ &ownSG$(84,7)$&$7$\\
          & $4,5$ & $1,16$ &$21$ &RBIBD$(85,5)-1$&$21$\\
          &       & $6,12$  &$22$ &RTD$(5,18)-6$&$18$\\
          &       & $11,8$ &$23$ &RTD$(5,19)-11$&$19$\\
          &       & $16,4$ &$25$ &RTD$(5,20)-16$&$20$\\
          & $5,6$ & $6,9$ & $17$ &RTD$(5,15)\_Q(14,1)+9$&$14$\\
          &       & $12,4$ & $19$ &RTD$(6,16)-12$&$16$\\
    \hline
    $85$ & $5$ & $17,0$ &$21$ &RBIBD$(85,5)$&$21$\\
          & $4,5$ & $5,13$ &$22$ &RTD$(5,18)-5$&$18$\\
          &       & $10,9$  &$23$ &RTD$(5,19)-10$&$19$\\
          &       & $15,5$ &$25$ &RTD$(5,20)-15$&$20$\\
          &       & $20,1$ &$27$ &RTD$(5,21)-20$&$21$\\
          & $5,6$ & $5,10$ & $17$ &RTD$(5,15)\_Q(14,1)+10$&$14$\\
          &       & $11,5$ & $19$ &RTD$(6,16)-11$&$16$\\
    \hline
    $86$ & $4,5$ & $4,14$ &$22$ &RTD$(5,18)-4$&$18$\\
          &       & $9,10$  &$23$ &RTD$(5,19)-9$&$19$\\
          &       & $14,6$ &$25$ &RTD$(5,20)-14$&$20$\\
          &       & $19,2$ &$27$ &RTD$(5,21)-19$&$21$\\
          & $5,6$ & $4,11$ & $17$ &ownSG$(90,6)-4$&$10$\\
          &       & $10,6$ & $19$ &RTD$(6,16)-10$&$16$\\
          &       & $16,1$ & $20$ &RTD$(6,17)-16$&$17$\\
    \hline
    $87$ &  $3$  & $29,0$ & $43$ &KTS$(87)$&$43$\\
         & $4,5$ & $3,15$ &$22$ &RTD$(5,18)-3$&$18$\\
         &       & $8,11$  &$23$ &RTD$(5,19)-8$&$19$\\
         &       & $13,7$ &$25$ &RTD$(5,20)-13$&$20$\\
         &       & $18,3$ &$27$ &RTD$(5,21)-18$&$21$\\
         & $5,6$ & $3,12$ & $17$ &ownSG$(90,6)-3$&$10$\\
         &       & $9,7$ & $19$ &RTD$(6,16)-9$&$16$\\
         &       & $15,2$ & $20$ &RTD$(6,17)-15$&$17$\\
    
\hline
\end{tabular}
\caption{\label{table:parallel12} SGP and SGA solutions for  $84\leq v\leq 87$. \\Block sizes  $K=\{3\}, \{4\},\{5\},\{6\},\{7\},\{4,5\}$, or $\{5,6\}$.}
\end{table}

\begin{table}[p]
\centering
\begin{tabular}{|c|c|c|c|c|c|}
\hline
    $v$ & $K$ & $m_{1},m_{2}$ & $MAX$ & solution& $NR$  \\
    \hline
    $88$ &  $4$  & $22,0$ & $29$ &RBIBD$(88,4)$&$29$\\
          &  $8$  & $11,0$ & $12$ &RTD$(8,11)$&$11$\\
          & $4,5$ & $2,16$ &$22$ &RTD$(5,18)-2$&$18$\\
          &       & $7,12$  &$23$ &RTD$(5,19)-7$&$19$\\
          &       & $12,8$ &$25$ &RTD$(5,20)-12$&$20$\\
          &       & $17,4$ &$26$ &RTD$(5,21)-17$&$21$\\
          & $5,6$ & $2,13$ & $17$ &ownSG$(90,6)-2$&$10$\\
          &       & $8,8$ & $19$ &RTD$(6,16)-8$&$16$\\
          &       & $14,3$ & $20$ &RTD$(6,17)-14$&$17$\\
    \hline    
$89$  & $4,5$ & $1,17$ &$22$ &RTD$(5,18)-1$&$18$\\
          &       & $6,13$  &$23$ &RTD$(5,19)-6$&$19$\\
          &       & $11,9$ &$25$ &RTD$(5,20)-11$&$20$\\
          &       & $16,5$ &$26$ &RTD$(5,21)-16$&$21$\\
          &       & $21,1$ &$28$ &RITD$(22,3;5)-21$&$19$\\
          & $5,6$ & $1,14$ & $17$ &ownSG$(90,6)-1$&$10$\\
          &       & $7,9$ & $19$ &RTD$(6,16)-7$&$16$\\
          &       & $13,4$ & $20$ &RTD$(6,17)-13$&$17$\\
    \hline
     $90$ & $3$ & $30,0$ &$44$ &NKTS$(90)$&$44$\\
          & $5$ & $18,0$ &$22$ &RTD$(5,18)$&$18$\\
          & $6$ & $15,0$ &$17$ &ownSG$(90,6)$&$10$\\
          & $9$ & $10,0$ &$11$ &ownSG$(90,9)$&$5$\\
          & $4,5$ & $5,14$ &$23$ &RTD$(5,19)-5$&$19$\\
          &       & $10,10$  &$25$ &RTD$(5,20)-10$&$20$\\
          &       & $15,6$ &$26$ &RTD$(5,21)-15$&$21$\\
          &       & $20,2$ &$28$ &RITD$(22,3;5)-20$&$19$\\
          & $5,6$ & $6,10$ & $19$ &RTD$(6,16)-6$&$16$\\
          &       & $12,5$ & $20$ &RTD$(6,17)-12$&$17$\\
         \hline
\end{tabular}
\caption{\label{table:parallel13} SGP and SGA solutions for  $88\leq v\leq 90$. \\Block sizes  $K=\{3\}, \{4\},\{5\},\{6\},\{8\},\{9\},\{4,5\}$, or $\{5,6\}$.}
\end{table}

\begin{table}[p]
\centering
\begin{tabular}{|c|c|c|c|c|c|}
\hline
    $v$ & $K$ & $m_{1},m_{2}$ & $MAX$ & solution& $NR$  \\
     \hline
     $91$ & $7$ & $13,0$ &$15$ &RTD$(7,13)$&$13$\\
          & $4,5$ & $4,15$ &$23$ &RTD$(5,19)-4$&$19$\\
          &       & $9,11$  &$24$ &RTD$(5,20)-9$&$20$\\
          &       & $14,7$ &$26$ &RTD$(5,21)-14$&$21$\\
          &       & $19,3$ &$28$ &RITD$(22,3;5)-19$&$19$\\
          & $5,6$ & $5,11$ & $19$ &RTD$(6,16)-5$&$16$\\
          &       & $11,6$ & $20$ &RTD$(6,17)-11$&$17$\\
\hline     
\end{tabular}
\caption{\label{table:parallel14} SGP and SGA solutions for  $91\leq v\leq 94$. \\Block sizes  $K=\{3\}, \{4\},\{7\},\{4,5\}$, or $\{5,6\}$.}
\end{table}

\begin{table}[p]
\centering
\begin{tabular}{|c|c|c|c|c|c|}
\hline
    $v$ & $K$ & $m_{1},m_{2}$ & $MAX$ & solution& $NR$  \\
    \hline
     $95$ & $5$ & $19,0$ &$23$ &RTD$(5,19)$&$19$\\
          & $4,5$ & $5,15$ &$24$ &RTD$(5,20)-5$&$20$\\
          &       & $10,11$  &$26$ &RTD$(5,21)-10$&$21$\\
          &       & $15,7$ &$27$ &RITD$(22,3;5)-15$&$19$\\
          &       & $20,3$ &$29$ &RTD$(5,23)-20$&$23$\\
          & $5,6$ & $1,15$ & $19$ &RTD$(6,16)-1$&$16$\\
          &       & $7,10$ & $20$ &RTD$(6,17)-7$&$17$\\
          &       & $13,5$ & $21$ &RTD$(6,18)-13$&$18$\\
    \hline
    $96$ & $3$ & $32,0$ &$47$ &NKTS$(96)$&$47$\\
          & $4$ & $24,0$ &$31$ &RGDD$(96,4,3)$&$31$\\
          & $6$ & $16,0$ &$19$ &RTD$(6,16)$&$16$\\
          & $8$ & $12,0$ &$13$ &MOLRs$(8,12)$&$6$\\
          & $4,5$ & $4,16$ &$24$ &RTD$(5,20)-4$&$20$\\
          &       & $9,12$  &$26$ &RTD$(5,21)-9$&$21$\\
          &       & $14,8$ &$27$ &RITD$(22,3;5)-14$&$19$\\
          &       & $19,4$ &$29$ &RTD$(5,23)-19$&$23$\\
          & $5,6$ & $6,11$ & $20$ &RTD$(6,17)-6$&$17$\\
          &       & $12,6$ & $21$ &RTD$(6,18)-12$&$18$\\
          &       & $18,1$ & $23$ &RTD$(6,19)-18$&$19$\\
     \hline
    $97$  & $4,5$ & $3,17$ &$24$ &RTD$(5,20)-3$&$20$\\
          &       & $8,13$  &$26$ &RTD$(5,21)-8$&$21$\\
          &       & $13,9$ &$27$ &RITD$(22,3;5)-13$&$19$\\
          &       & $18,5$ &$29$ &RTD$(5,23)-18$&$23$\\
          &       & $23,1$ &$31$ &RTD$(5,24)-23$&$24$\\
          & $5,6$ & $5,12$ & $20$ &RTD$(6,17)-5$&$17$\\
          &       & $11,7$ & $21$ &RTD$(6,18)-11$&$18$\\
          &       & $17,2$ & $23$ &RTD$(6,19)-17$&$19$\\
    \hline
     
\end{tabular}
\caption{\label{table:parallel15} SGP and SGA solutions for  $95\leq v\leq 97$. \\Block sizes  $K=\{3\}, \{4\},\{5\},\{6\},\{8\},\{4,5\}$, or $\{5,6\}$.}
\end{table}

\begin{table}[p]
\centering
\begin{tabular}{|c|c|c|c|c|c|}
\hline
    $v$ & $K$ & $m_{1},m_{2}$ & $MAX$ & solution& $NR$  \\
    \hline
    $98$ & $7$ & $14,0$ &$16$ &ownSG$(98,7)+G(1)$&$9$\\
          & $4,5$ & $2,18$ &$24$ &RTD$(5,20)-2$&$20$\\
          &       & $7,14$  &$26$ &RTD$(5,21)-7$&$21$\\
          &       & $12,10$ &$27$ &RITD$(22,3;5)-12$&$19$\\
          &       & $17,6$ &$29$ &RTD$(5,23)-17$&$23$\\
          &       & $22,2$ &$31$ &RTD$(5,24)-22$&$24$\\
          & $5,6$ & $4,13$ & $20$ &RTD$(6,17)-4$&$17$\\
          &       & $10,8$ & $21$ &RTD$(6,18)-10$&$18$\\
          &       & $16,3$ & $23$ &RTD$(6,19)-16$&$19$\\
     \hline
    $99$ & $3$ & $33,0$ &$49$ &KTS$(99)$&$49$\\
          & $9$ & $11,0$ &$12$ &RGDD$(99,9,3)$&$12$\\
          & $4,5$ & $1,19$ &$24$ &RTD$(5,20)+G(1)-1$&$21$\\
          &       & $6,15$  &$26$ &RTD$(5,21)-6$&$21$\\
          &       & $11,11$ &$27$ &RITD$(22,3;5)-11$&$19$\\
          &       & $16,7$ &$29$ &RTD$(5,23)-16$&$23$\\
          &       & $21,3$ &$31$ &RTD$(5,24)-21$&$24$\\
          & $5,6$ & $3,14$ & $20$ &RTD$(6,17)-3$&$17$\\
          &       & $9,9$ & $21$ &RTD$(6,18)-9$&$18$\\
          &       & $15,4$ & $23$ &RTD$(6,19)-15$&$19$\\
     \hline
    $100$ & $4$ & $25,0$ &$33$ &RBIBD$(100,4)$&$33$\\
          & $5$ & $20,0$ &$24$ &RTD$(5,20)+G(1)$&$21$\\
          & $10$ & $10,0$ &$11$ &MOLRs$(10,10)+G(1)$&$4$\\
          & $4,5$ & $5,16$ &$26$ &RBIBD$(105,5)-5,B$&$25$\\
          &       & $10,12$  &$27$ &RITD$(22,3;5)-10$&$19$\\
          &       & $15,8$ &$29$ &RTD$(5,23)-15$&$23$\\
          &       & $20,4$ &$30$ &RTD$(5,24)-20$&$24$\\
          & $5,6$ & $2,15$ & $20$ &RTD$(6,17)-2$&$17$\\
          &       & $8,10$ & $21$ &RTD$(6,18)-8$&$18$\\
          &       & $14,5$ & $23$ &RTD$(6,19)-14$&$19$\\
     \hline
\end{tabular}
\caption{\label{table:parallel16} SGP and SGA solutions for  $98\leq v\leq 100$. \\Block sizes  $K=\{3\}, \{4\},\{5\},\{7\},\{9\},\{10\},\{4,5\}$, or $\{5,6\}$.}
\end{table}

\begin{table}[p]
\centering
\begin{tabular}{|c|c|c|c|c|c|}
\hline
    $v$ & $K$ & $m_{1},m_{2}$ & $MAX$ & solution& $NR$  \\
    \hline
    $101$  & $4,5$ & $4,17$ &$26$ &RBIBD$(105,5)-4,B$&$25$\\
          &       & $9,13$  &$27$ &RITD$(22,3;5)-6$&$19$\\
          &       & $14,9$ &$29$ &RTD$(5,23)-14$&$23$\\
          &       & $19,5$ &$30$ &RTD$(5,24)-19$&$24$\\
          &       & $24,1$ &$32$ &RTD$(5,25)-24$&$25$\\
          & $5,6$ & $1,16$ & $20$ &RTD$(6,17)-1$&$17$\\
          &       & $7,11$ & $21$ &RTD$(6,18)-7$&$18$\\
          &       & $13,6$ & $22$ &RTD$(6,19)-13$&$19$\\
          &       & $19,1$ & $24$ &RTD$(5,20)\_Q(19,1)+1$&$20$\\
    \hline
     $102$ & $3$ & $34,0$ &$50$ &NKTS$(102)$&$50$\\
           & $6$ & $17,0$ &$20$ &RTD$(6,17)$&$17$\\
          & $4,5$ & $3,18$ &$26$ &RBIBD$(105,5)-3,B$&$25$\\
          &       & $8,14$  &$27$ &RITD$(22,3;5)-8$&$19$\\
          &       & $13,10$ &$28$ &RTD$(5,23)-13$&$23$\\
          &       & $18,6$ &$30$ &RTD$(5,24)-18$&$24$\\
          &       & $23,2$ &$32$ &RTD$(5,25)-23$&$25$\\
          & $5,6$ & $6,12$ & $21$ &RTD$(6,18)-6$&$18$\\
          &       & $12,7$ & $22$ &RTD$(6,19)-12$&$19$\\
          &       & $18,2$ & $24$ &RTD$(5,20)\_Q(19,1)+2$&$19$\\
     \hline
    $103$ & $4,5$ & $2,19$ &$26$ &RBIBD$(105,5)-2,B$&$25$\\
          &       & $7,15$  &$27$ &RITD$(22,3;5)-7$&$19$\\
          &       & $12,11$ &$28$ &RTD$(5,23)-12$&$23$\\
          &       & $17,7$ &$30$ &RTD$(5,24)-17$&$24$\\
          &       & $22,3$ &$32$ &RTD$(5,25)-22$&$25$\\
          & $5,6$ & $5,13$ & $21$ &RTD$(6,18)-5$&$18$\\
          &       & $11,8$ & $22$ &RTD$(6,19)-11$&$19$\\
          &       & $17,3$ & $24$ &RTD$(5,20)\_Q(19,1)+3$&$19$\\
     \hline
     
\end{tabular}
\caption{\label{table:parallel17} SGP and SGA solutions for  $101\leq v\leq 103$. \\Block sizes  $K=\{3\}, \{6\},\{4,5\}$, or $\{5,6\}$.}
\end{table}

\begin{table}[p]
\centering
\begin{tabular}{|c|c|c|c|c|c|}
\hline
    $v$ & $K$ & $m_{1},m_{2}$ & $MAX$ & solution& $NR$  \\
    \hline
    $104$ & $4$ & $26,0$ &$34$ &RGDD$(104,4,2)$&$34$\\
          & $8$ & $13,0$ &$14$ &RTD$(8,13)$&$13$\\
          & $4,5$ & $1,20$ &$26$ &RBIBD$(21,5)-1$&$26$\\
          &       & $6,16$  &$27$ &RITD$(22,3;5)-6$&$19$\\
          &       & $11,12$ &$28$ &RTD$(5,23)-11$&$23$\\
          &       & $16,8$ &$30$ &RTD$(5,24)-16$&$24$\\
          &       & $21,4$ &$32$ &RTD$(5,25)-21$&$25$\\
          & $5,6$ & $4,14$ & $21$ &RTD$(6,18)-4$&$18$\\
          &       & $10,9$ & $22$ &RTD$(6,19)-10$&$19$\\
          &       & $16,4$ & $24$ &RTD$(5,20)\_Q(19,1)+4$&$19$\\
     \hline
    $105$  
          & $3$ & $35,0$ &$52$ &KTS$(105)$&$52$\\
          & $5$ & $21,0$ &$26$ &RBIBD$(105,5)$&$26$\\
          & $7$ & $15,0$ &$17$ &ownSG$(105,7)$&$9$\\
          & $4,5$ & $5,17$ &$27$ &RITD$(22,3;5)-5$&$19$\\
          &       & $10,13$  &$28$ &RTD$(5,23)-10$&$23$\\
          &       & $15,9$ &$30$ &RTD$(5,24)-15$&$24$\\
          &       & $20,5$ &$32$ &RTD$(5,25)-20$&$25$\\
          &       & $25,1$ &$34$ &RTD$(5,26)-25$&$26$\\
          & $5,6$ & $3,15$ & $21$ &RTD$(6,18)-3$&$18$\\
          &       & $9,10$ & $22$ &RTD$(6,19)-9$&$19$\\
          &       & $15,5$ & $24$ &RTD$(5,20)\_Q(19,1)+5$&$19$\\
    \hline
    $106$ 
          & $4,5$ & $4,18$ &$27$ &RITD$(22,3;5)-4$&$19$\\
          &       & $9,14$  &$28$ &RTD$(5,23)-9$&$23$\\
          &       & $14,10$ &$30$ &RTD$(5,24)-14$&$24$\\
          &       & $19,6$ &$31$ &RTD$(5,25)-19$&$25$\\
          &       & $24,2$ &$33$ &RTD$(5,26)-24$&$26$\\
          & $5,6$ & $2,16$ & $21$ &RTD$(6,18)-2$&$18$\\
          &       & $8,11$ & $22$ &RTD$(6,19)-8$&$19$\\
          &       & $14,6$ & $24$ &RTD$(5,20)\_Q(19,1)+6$&$19$\\
          &       & $20,1$ & $25$ &RTD$(6,21)-20$&$21$\\
     \hline
\end{tabular}
\caption{\label{table:parallel18} SGP and SGA solutions for  $104\leq v\leq 106$. \\Block sizes  $K=\{3\}, \{4\},\{5\},\{7\},\{8\},\{4,5\}$, or $\{5,6\}$.}
\end{table}

\begin{table}[p]
\centering
\begin{tabular}{|c|c|c|c|c|c|}
\hline
    $v$ & $K$ & $m_{1},m_{2}$ & $MAX$ & solution& $NR$  \\
    \hline
    $107$ 
          & $4,5$ & $3,19$ &$27$ &RITD$(22,3;5)-3$&$19$\\
          
          &       & $8,15$  &$28$ &RTD$(5,23)-8$&$23$\\
          &       & $13,11$ &$30$ &RTD$(5,24)-13$&$24$\\
          &       & $18,7$ &$31$ &RTD$(5,25)-18$&$25$\\
          &       & $23,3$ &$33$ &RTD$(5,26)-23$&$26$\\
          & $5,6$ & $1,17$ & $21$ &RTD$(6,18)+G(1)-1$&$19$\\
          &       & $7,12$ & $22$ &RTD$(6,19)-7$&$19$\\
          &       & $13,7$ & $24$ &RTD$(5,20)\_Q(19,1)+7$&$19$\\
          &       & $19,2$ & $25$ &RTD$(6,21)-19$&$21$\\
     \hline
    $108$ & $3$ & $36,0$ &$53$ &NKTS$(108)$&$53$\\
           & $4$ & $27,0$ &$35$ &RGDD$(108,4,3)$&$35$\\
           & $6$ & $18,0$ &$21$ &RTD$(6,18)+G(1)$&$19$\\
           & $9$ & $12,0$ &$13$ &MOLRs$(9,12)$&$6$\\
          & $4,5$ & $2,20$ &$27$ &RITD$(22,3;5)-2$&$19$\\
          &       & $7,16$  &$28$ &RTD$(5,23)-7$&$23$\\
          &       & $12,12$ &$30$ &RTD$(5,24)-12$&$24$\\
          &       & $17,8$ &$31$ &RTD$(5,25)-17$&$25$\\
          &       & $22,4$ &$33$ &RTD$(5,26)-22$&$26$\\
          & $5,6$ & $6,13$ & $22$ &RTD$(6,19)-6$&$19$\\
          &       & $12,8$ & $24$ &RTD$(5,20)\_Q(19,1)+8$&$19$\\
          &       & $18,3$ & $25$ &RTD$(6,21)-18$&$21$\\
    \hline
    $109$ 
          & $4,5$ & $1,21$ &$27$ &RITD$(22,3;5)-1$&$19$\\
          &       & $6,17$  &$28$ &RTD$(5,23)-6$&$23$\\
          &       & $11,13$ &$30$ &RTD$(5,24)-11$&$24$\\
          &       & $16,9$ &$31$ &RTD$(5,25)-16$&$25$\\
          &       & $21,5$ &$33$ &RTD$(5,26)-21$&$26$\\
          &       & $26,1$ &$35$ &RTD$(5,27)-26$&$27$\\
          & $5,6$ & $5,14$ & $22$ &RTD$(6,19)-5$&$19$\\
          &       & $11,9$ & $24$ &RTD$(5,20)\_Q(19,1)+9$&$19$\\
          &       & $17,4$ & $25$ &RTD$(6,21)-17$&$21$\\
     \hline
\end{tabular}
\caption{\label{table:parallel19} SGP and SGA solutions for  $107\leq v\leq 109$. \\Block sizes  $K=\{3\}, \{4\},\{6\},\{9\},\{4,5\}$, or $\{5,6\}$.}
\end{table}

\begin{table}[p]
\centering
\begin{tabular}{|c|c|c|c|c|c|}
\hline
    $v$ & $K$ & $m_{1},m_{2}$ & $MAX$ & solution& $NR$  \\
    \hline
    $110$ & $5$ & $22,0$ &$27$ & RITD$(22,3;5)$&$19$\\
          & $10$ & $11,0$ &$12$ &RTD$(10,11)$&$11$\\
          & $4,5$ & $5,18$ &$28$ &RTD$(5,23)-5$&$23$\\
          &       & $10,14$  &$29$ &RTD$(5,24)-10$&$24$\\
          &       & $15,10$ &$31$ &RTD$(5,25)-15$&$25$\\
          &       & $20,6$ &$33$ &RTD$(5,26)-20$&$26$\\
          &       & $25,2$ &$35$ &RTD$(5,27)-25$&$27$\\
          & $5,6$ & $4,15$ & $22$ &RTD$(6,19)-4$&$19$\\
          &       & $10,10$ & $23$ &RTD$(5,20)\_Q(19,1)+10$&$19$\\
          &       & $16,5$ & $25$ &RTD$(6,21)-16$&$21$\\
     \hline
    $111$ & $3$ & $37,0$ &$55$ &KTS$(111)$&$55$\\
          & $4,5$ & $4,19$ &$28$ &RTD$(5,23)-4$&$23$\\
          &       & $9,15$  &$29$ &RTD$(5,24)-9$&$24$\\
          &       & $14,11$ &$31$ &RTD$(5,25)-14$&$25$\\
          &       & $19,7$ &$33$ &RTD$(5,26)-19$&$26$\\
          &       & $24,3$ &$35$ &RTD$(5,27)-24$&$27$\\
          & $5,6$ & $3,16$ & $22$ &RTD$(6,19)-3$&$19$\\
          &       & $9,11$ & $23$ &RTD$(5,20)\_Q(19,1)+11$&$19$\\
          &       & $15,6$ & $25$ &RTD$(6,21)-15$&$21$\\
          &       & $21,1$ & $27$ &ownSG$(132,6)-21$ & $14$\\
    \hline
    $112$ & $4$ & $28,0$ &$37$ &RBIBD$(112,4)$&$37$\\
           & $7$ & $16,0$ &$18$ &RTD$(7,16)$&$16$\\
           & $8$ & $14,0$ &$15$ &ownSG$(112,8)$&$8$\\
          & $4,5$ & $3,20$ &$28$ &RTD$(5,23)-3$&$23$\\
          &       & $8,16$  &$29$ &RTD$(5,24)-8$&$24$\\
          &       & $13,12$ &$31$ &RTD$(5,25)-13$&$25$\\
          &       & $18,8$ &$33$ &RTD$(5,26)-18$&$26$\\
          &       & $23,4$ &$34$ &RTD$(5,27)-23$&$27$\\
          & $5,6$ & $2,17$ & $22$ &RTD$(6,19)-2$&$19$\\
          &       & $8,12$ & $23$ &RTD$(5,20)\_Q(19,1)+12$&$19$\\
          &       & $14,7$ & $25$ &RTD$(6,21)-14$&$21$\\
          &       & $20,2$ & $26$ &ownSG$(132,6)-20$ & $14$\\
     \hline
     
\end{tabular}
\caption{\label{table:parallel20} SGP and SGA solutions for  $110\leq v\leq 112$. \\Block sizes  $K=\{3\}, \{4\},\{5\},\{7\},\{8\},\{10\},\{4,5\}$, or $\{5,6\}$.}
\end{table}

\begin{table}[p]
\centering
\begin{tabular}{|c|c|c|c|c|c|}
\hline
    $v$ & $K$ & $m_{1},m_{2}$ & $MAX$ & solution& $NR$  \\
    \hline
    $113$& $4,5$ & $2,21$ &$28$ &RTD$(5,23)-2$&$23$\\
          &       & $7,17$  &$29$ &RTD$(5,24)-7$&$24$\\
          &       & $12,13$ &$31$ &RTD$(5,25)-12$&$25$\\
          &       & $17,9$ &$32$ &RTD$(5,26)-17$&$26$\\
          &       & $22,5$ &$34$ &RTD$(5,27)-22$&$27$\\
          &       & $27,1$ &$36$ &RTD$(5,28)-27$&$28$\\
          & $5,6$ & $1,18$ & $22$ &RTD$(6,19)-1$&$19$\\
          &       & $7,13$ & $23$ &RTD$(5,20)\_Q(19,1)+13$&$19$\\
          &       & $13,8$ & $25$ &RTD$(6,21)-13$&$21$\\
          &       & $19,3$ & $26$ &ownSG$(132,6)-19$ & $14$\\
     \hline
    $114$ & $3$ & $38,0$ &$56$ &NKTS$(114)$&$56$\\
           & $6$ & $19,0$ &$22$ &RTD$(6,19)$&$19$\\
          & $4,5$ & $1,22$ &$28$ &RTD$(5,23)-1$&$23$\\
          &       & $6,18$  &$29$ &RTD$(5,24)-6$&$24$\\
          &       & $11,14$ &$31$ &RTD$(5,25)-11$&$25$\\
          &       & $16,10$ &$32$ &RTD$(5,26)-16$&$26$\\
          &       & $21,6$ &$34$ &RTD$(5,27)-21$&$27$\\
          &       & $26,2$ &$36$ &RTD$(5,28)-26$&$28$\\
          & $5,6$ & $6,14$ & $23$ &RTD$(5,20)\_Q(19,1)+14$&$19$\\
          &       & $12,9$ & $25$ &RTD$(6,21)-12$&$21$\\
          &       & $18,4$ & $26$ &ownSG$(132,6)-18$ & $14$\\
    \hline
    $115$ & $5$ & $23,0$ &$28$ &RTD$(5,23)$&$23$\\
          & $4,5$ & $5,19$  &$29$ &RTD$(5,24)-5$&$24$\\
          &       & $10,15$ &$31$ &RTD$(5,25)-10$&$25$\\
          &       & $15,11$ &$32$ &RTD$(5,26)-15$&$26$\\
          &       & $20,7$ &$34$ &RTD$(5,27)-20$&$27$\\
          &       & $25,3$ &$36$ &RTD$(5,28)-25$&$28$\\
          & $5,6$ & $5,15$ & $23$ &RTD$(5,20)\_Q(19,1)+15$&$19$\\
          &       & $11,10$ & $25$ &RTD$(6,21)-11$&$21$\\
          &       & $17,5$ & $26$ &ownSG$(132,6)-17$ & $14$\\
     \hline
\end{tabular}
\caption{\label{table:parallel21} SGP and SGA solutions for  $113\leq v\leq 115$. \\Block sizes  $K=\{3\},\{5\},\{6\},\{4,5\}$, or $\{5,6\}$.}
\end{table}

\begin{table}
\centering
\begin{tabular}{|c|c|c|c|c|c|}
\hline
    $v$ & $K$ & $m_{1},m_{2}$ & $MAX$ & solution& $NR$  \\
    \hline
    $116$ & $4$ & $29,0$ &$38$ &RGDD$(116,4,2)$&$38$\\
          & $4,5$ & $4,20$ &$29$ &RGDD$(120,5,4)-4$&$29$\\
          &       & $9,16$  &$31$ &RTD$(5,25)-9$&$25$\\
          &       & $14,12$ &$32$ &RTD$(5,26)-14$&$26$\\
          &       & $19,8$ &$34$ &RTD$(5,27)-19$&$27$\\
          &       & $24,4$ &$36$ &RTD$(5,28)-24$&$28$\\
          & $5,6$ & $4,16$ & $23$ &ownSG$(120,6)-4$&$13$\\
          &       & $10,11$ & $25$ &RTD$(6,21)-10$&$21$\\
          &       & $16,6$ & $26$ &ownSG$(132,6)-16$ & $14$\\
          &       & $22,1$ & $28$ &RTD$(6,23)-22$&$23$\\
    \hline
    $117$ & $3$ & $39,0$ &$58$ &KTS$(117)$&$58$\\
           & $9$ & $13,0$ &$14$ &RTD$(9,13)$&$13$\\
          & $4,5$ & $3,21$ &$29$ &RGDD$(120,5,4)-3$&$29$\\
          &       & $8,17$  &$31$ &RTD$(5,25)-8$&$25$\\
          &       & $13,13$ &$32$ &RTD$(5,26)-13$&$26$\\
          &       & $18,9$ &$34$ &RTD$(5,27)-18$&$27$\\
          &       & $23,5$ &$36$ &RTD$(5,28)-23$&$28$\\
          &       & $28,1$ &$38$ &RTD$(5,29)-28$&$29$\\
          & $5,6$ & $3,17$ & $23$ &ownSG$(120,6)-3$&$13$\\
          &       & $9,12$ & $25$ &RTD$(6,21)-9$&$21$\\
          &       & $15,7$ & $26$ &ownSG$(132,6)-15$ & $14$\\
          &       & $21,2$ & $28$ &RTD$(6,23)-21$&$23$\\
    \hline
    $118$ 
          & $4,5$ & $2,22$ &$29$ &RGDD$(120,5,4)-2$&$29$\\
          &       & $7,18$  &$31$ &RTD$(5,25)-7$&$25$\\
          &       & $12,14$ &$32$ &RTD$(5,26)-12$&$26$\\
          &       & $17,10$ &$34$ &RTD$(5,27)-17$&$27$\\
          &       & $22,6$ &$35$ &RTD$(5,28)-22$&$28$\\
          &       & $27,2$ &$37$ &RTD$(5,29)-27$&$29$\\
          & $5,6$ & $2,18$ & $23$ &ownSG$(120,6)-2$&$13$\\
          &       & $8,13$ & $25$ &RTD$(6,21)-8$&$21$\\
          &       & $14,8$ & $26$ &ownSG$(132,6)-14$ & $14$\\
          &       & $20,3$ & $28$ &RTD$(6,23)-20$&$23$\\
     \hline
\end{tabular}
\caption{\label{table:parallel22} SGP and SGA solutions for  $116\leq v\leq 118$. \\Block sizes  $K=\{3\}, \{4\},\{9\},\{4,5\}$, or $\{5,6\}$.}
\end{table}

\begin{table}[p]
\centering
\begin{tabular}{|c|c|c|c|c|c|}
\hline
    $v$ & $K$ & $m_{1},m_{2}$ & $MAX$ & solution& $NR$  \\
    \hline
    $119$ &  $7$ & $17,0$ &$19$ &RTD$(7,17)$&$17$\\
          & $4,5$ & $1,23$ &$29$ &RGDD$(120,5,4)-1$&$29$\\
          &       & $6,19$  &$31$ &RTD$(5,25)-6$&$25$\\
          &       & $11,15$ &$32$ &RTD$(5,26)-11$&$26$\\
          &       & $16,11$ &$34$ &RTD$(5,27)-16$&$27$\\
          &       & $21,7$ &$35$ &RTD$(5,28)-21$&$28$\\
          &       & $26,3$ &$37$ &RTD$(5,29)-26$&$29$\\
          & $5,6$ & $1,19$ & $23$ &ownSG$(120,6)-1$&$13$\\
          &       & $7,14$ & $25$ &RTD$(6,21)-7$&$21$\\
          &       & $13,9$ & $26$ &ownSG$(132,6)-13$ & $14$\\
          &       & $19,4$ & $28$ &RTD$(6,23)-19$&$23$\\
     \hline
    $120$ & $3$ & $40,0$ &$59$ &NKTS$(120)$&$59$\\
           & $4$ & $30,0$ &$39$ &RGDD$(120,4,3)$&$39$\\
           & $5$ & $24,0$ &$29$ &RGDD$(120,5,4)$&$29$\\
           & $6$ & $20,0$ &$23$ &ownSG$(120,6)$&$13$\\
           & $8$ & $15,0$ &$17$ &RBIBD$(120,8)$&$17$\\
           & $10$ & $12,0$ &$13$ &MOLRs$(10,12)$&$6$\\
          & $4,5$ & $5,20$ &$31$ &RBIBD$(125,5)-5,B$&$30$\\
          &       & $10,16$  &$32$ &RTD$(5,26)-10$&$26$\\
          &       & $15,12$ &$34$ &RTD$(5,27)-15$&$27$\\
          &       & $20,8$ &$35$ &RTD$(5,28)-20$&$28$\\
          &       & $25,4$ &$37$ &RTD$(5,29)-25$&$29$\\
          & $5,6$ & $6,15$ & $25$ &RBIBD$(126,6)-6,B$&$24$\\
          &       & $12,10$ & $26$ &ownSG$(132,6)-12$ & $14$\\
          &       & $18,5$ & $28$ &RTD$(6,23)-18$&$23$\\
    \hline
\end{tabular}
\caption{\label{table:parallel23} SGP and SGA solutions for  $119\leq v\leq 120$. \\Block sizes  $K=\{3\}, \{4\},\{5\},\{6\},\{7\},\{8\},\{10\},\{4,5\}$, or $\{5,6\}$.}
\end{table}

\begin{table}
\centering
\begin{tabular}{|c|c|c|c|c|c|}
\hline
    $v$ & $K$ & $m_{1},m_{2}$ & $MAX$ & solution& $NR$  \\
    \hline
    $121$ & $11$ & $11,0$ &$12$ &RBIBD$(121,11)$&12\\
          & $4,5$ & $4,21$ &$31$ &RBIBD$(125,5)-4,B$&$30$\\
          &       & $9,17$  &$32$ &RTD$(5,26)-9$&$26$\\
          &       & $14,13$ &$33$ &RTD$(5,27)-14$&$27$\\
          &       & $19,9$ &$35$ &RTD$(5,28)-19$&$28$\\
          &       & $24,5$ &$37$ &RTD$(5,29)-24$&$29$\\
          &       & $29,1$ &$39$ &RTD$(5,30)-29$&$30$\\
          & $5,6$ & $5,16$ & $25$ &RBIBD$(126,6)-5,B$&$24$\\
          &       & $11,11$ & $26$ &ownSG$(132,6)-11$ & $14$\\
          &       & $17,6$ & $27$ &RTD$(6,23)-17$&$23$\\
          &       & $23,1$ & $29$ &RTD$(6,24)-23$&$24$\\
     \hline
     $122$ & $4,5$ & $3,22$ &$31$ &RBIBD$(125,5)-3,B$&$30$\\
          &       & $8,18$  &$32$ &RTD$(5,26)-8$&$26$\\
          &       & $13,14$ &$33$ &RTD$(5,27)-13$&$27$\\
          &       & $18,10$ &$35$ &RTD$(5,28)-18$&$28$\\
          &       & $23,6$ &$37$ &RTD$(5,29)-23$&$29$\\
          &       & $28,2$ &$39$ &RTD$(5,30)-28$&$30$\\
          & $5,6$ & $4,17$ & $25$ &RBIBD$(126,6)-4,B$&$24$\\
          &       & $10,12$ & $26$ &ownSG$(132,6)-10$ & $14$\\
          &       & $16,7$ & $27$ &RTD$(6,23)-16$&$23$\\
          &       & $22,2$ & $29$ &RTD$(6,24)-22$&$24$\\
     \hline
\end{tabular}
\caption{\label{table:parallel24} SGP and SGA solutions for  $121\leq v\leq 122$. \\Block sizes  $K=\{11\},\{4,5\}$, or $\{5,6\}$.}
\end{table}

\begin{table}
\centering
\begin{tabular}{|c|c|c|c|c|c|}
\hline
    $v$ & $K$ & $m_{1},m_{2}$ & $MAX$ & solution& $NR$  \\
    \hline
    $123$ & $3$ & $41,0$ &$61$ &KTS$(123)$&$61$\\
          & $4,5$ & $2,23$ &$31$ &RBIBD$(125,5)-2,B$&$30$\\
          &       & $7,19$  &$32$ &RTD$(5,26)-7$&$26$\\
          &       & $12,15$ &$33$ &RTD$(5,27)-12$&$27$\\
          &       & $17,11$ &$35$ &RTD$(5,28)-17$&$28$\\
          &       & $22,7$ &$37$ &RTD$(5,29)-22$&$29$\\
          &       & $27,3$ &$39$ &RTD$(5,30)-27$&$30$\\
          & $5,6$ & $3,18$ & $25$ &RBIBD$(126,6)-3,B$&$24$\\
          &       & $9,13$ & $26$ &ownSG$(132,6)-9$ & $14$\\
          &       & $15,8$ & $27$ &RTD$(6,23)-15$&$23$\\
          &       & $21,3$ & $29$ &RTD$(6,24)-21$&$24$\\
    \hline
     $124$ & $4$ & $31,0$ &$41$ &RBIBD$(124,4)$&$41$\\
          & $4,5$ & $1,24$ &$31$ &RBIBD$(125,5)-1$&$31$\\
          &       & $6,20$  &$32$ &RTD$(5,26)-6$&$26$\\
          &       & $11,16$ &$33$ &RTD$(5,27)-11$&$27$\\
          &       & $16,12$ &$35$ &RTD$(5,28)-16$&$28$\\
          &       & $21,8$ &$37$ &RTD$(5,29)-21$&$29$\\
          &       & $26,4$ &$38$ &RTD$(5,30)-26$&$30$\\
          & $5,6$ & $2,19$ & $25$ &RBIBD$(126,6)-2,B$&$24$\\
          &       & $8,14$ & $26$ &ownSG$(132,6)-8$ & $14$\\
          &       & $14,9$ & $27$ &RTD$(6,23)-14$&$23$\\
          &       & $20,4$ & $29$ &RTD$(6,24)-20$&$24$\\
     \hline
     
\end{tabular}
\caption{\label{table:paralle25} SGP and SGA solutions for  $123\leq v\leq 124$. \\Block sizes  $K=\{3\}, \{4\},\{4,5\}$, or $\{5,6\}$.}
\end{table}

\begin{table}
\centering
\begin{tabular}{|c|c|c|c|c|c|}
\hline
    $v$ & $K$ & $m_{1},m_{2}$ & $MAX$ & solution& $NR$  \\
    \hline
    $125$ & $5$ & $25,0$ &$31$ &RBIBD$(125,5)$&$31$\\
          & $4,5$ & $5,21$ &$32$ &RTD$(5,26)-5$&$26$\\
          &       & $10,17$  &$33$ &RTD$(5,27)-10$&$27$\\
          &       & $15,13$ &$35$ &RTD$(5,28)-15$&$28$\\
          &       & $20,9$ &$36$ &RTD$(5,29)-20$&$29$\\
          &       & $25,5$ &$38$ &RTD$(5,30)-25$&$30$\\
          &       & $30,1$ &$40$ &RTD$(5,31)-30$&$31$\\
          & $5,6$ & $1,20$ & $25$ &RBIBD$(126,6)-1$&$25$\\
          &       & $7,15$ & $26$ &ownSG$(132,6)-7$ & $14$\\
          &       & $13,10$ & $27$ &RTD$(6,23)-13$&$23$\\
          &       & $19,5$ & $29$ &RTD$(6,24)-19$&$24$\\
     \hline
     $126$ & $3$ & $42,0$ &$62$ &NKTS$(126)$&$62$\\
           & $6$ & $21,0$ &$25$ &RBIBD$(126,6)$&$25$\\
           & $7$ & $18,0$ &$20$ &ownSG$(126,7)$&$10$\\
           & $9$ & $14,0$ &$15$ &ownSG$(126,9)$&$7$\\
          & $4,5$ & $4,22$ &$32$ &RTD$(5,26)-4$&$26$\\
          &       & $9,18$  &$33$ &RTD$(5,27)-9$&$27$\\
          &       & $14,14$ &$35$ &RTD$(5,28)-14$&$28$\\
          &       & $19,10$ &$36$ &RTD$(5,29)-19$&$29$\\
          &       & $24,6$ &$38$ &RTD$(5,30)-24$&$30$\\
          &       & $29,2$ &$40$ &RTD$(5,31)-29$&$31$\\
          & $5,6$ & $6,16$ & $26$ &ownSG$(132,6)-6$ & $14$\\
          &       & $12,11$ & $27$ &RTD$(6,23)-12$&$23$\\
          &       & $18,6$ & $29$ &RTD$(6,24)-18$&$24$\\
          &       & $24,1$ & $30$ &RTD$(6,25)-24$&$25$\\
     \hline
     
\end{tabular}
\caption{\label{table:parallel26} SGP and SGA solutions for  $125\leq v\leq 126$. \\Block sizes  $K=\{3\},\{5\},\{6\},\{7\},\{9\},\{4,5\}$, or $\{5,6\}$.}
\end{table}

\begin{table}
\centering
\begin{tabular}{|c|c|c|c|c|c|}
\hline
    $v$ & $K$ & $m_{1},m_{2}$ & $MAX$ & solution& $NR$  \\
    \hline
    $127$ & $4,5$ & $3,23$ &$32$ &RTD$(5,26)-3$&$26$\\
          &       & $8,19$  &$33$ &RTD$(5,27)-8$&$27$\\
          &       & $13,15$ &$35$ &RTD$(5,28)-13$&$28$\\
          &       & $18,11$ &$36$ &RTD$(5,29)-18$&$29$\\
          &       & $23,7$ &$38$ &RTD$(5,30)-23$&$30$\\
          &       & $28,3$ &$40$ &RTD$(5,31)-28$&$31$\\
          & $5,6$ & $5,17$ & $26$ &ownSG$(132,6)-5$ & $14$\\
          &       & $11,12$ & $27$ &RTD$(6,23)-11$&$23$\\
          &       & $17,7$ & $29$ &RTD$(6,24)-17$&$24$\\
          &       & $23,2$ & $30$ &RTD$(6,25)-23$&$25$\\
     \hline
     $128$ & $4$ & $32,0$ &$42$ &RGDD$(128,4,2)$&$42$\\
           & $8$ & $16,0$ &$18$ &RTD$(8,16)+G(1)$&$17$\\
          & $4,5$ & $2,24$ &$32$ &RTD$(5,26)-2$&$26$\\
          &       & $7,20$  &$33$ &RTD$(5,27)-7$&$27$\\
          &       & $12,16$ &$35$ &RTD$(5,28)-12$&$28$\\
          &       & $17,12$ &$36$ &RTD$(5,29)-17$&$29$\\
          &       & $22,8$ &$38$ &RTD$(5,30)-22$&$30$\\
          &       & $27,4$ &$40$ &RTD$(5,31)-27$&$31$\\
          & $5,6$ & $4,18$ & $26$ &ownSG$(132,6)-4$ & $14$\\
          &       & $10,13$ & $27$ &RTD$(6,23)-10$&$23$\\
          &       & $16,8$ & $29$ &RTD$(6,24)-16$&$24$\\
          &       & $22,3$ & $30$ &RTD$(6,25)-22$&$25$\\
     \hline
     
\end{tabular}
\caption{\label{table:parallel27} SGP and SGA solutions for  $127\leq v\leq 128$. \\Block sizes  $K=\{4\},\{8\},\{4,5\}$, or $\{5,6\}$.}
\end{table}

\begin{table}
\centering
\begin{tabular}{|c|c|c|c|c|c|}
\hline
    $v$ & $K$ & $m_{1},m_{2}$ & $MAX$ & solution& $NR$  \\
    \hline
    $129$ & $3$ & $43,0$ &$64$ &KTS$(129)$&$64$\\
          &$4,5$ & $1,25$ &$32$ &RTD$(5,26)-1$&$26$\\
          &       & $6,21$  &$33$ &RTD$(5,27)-6$&$27$\\
          &       & $11,17$ &$34$ &RTD$(5,28)-11$&$28$\\
          &       & $16,13$ &$36$ &RTD$(5,29)-16$&$29$\\
          &       & $21,9$ &$38$ &RTD$(5,30)-21$&$30$\\
          &       & $26,5$ &$40$ &RTD$(5,31)-26$&$31$\\
          &       & $31,1$ &$42$ &RTD$(5,32)-31$&$32$\\
          & $5,6$ & $3,19$ & $26$ &ownSG$(132,6)-3$ & $14$\\
          &       & $9,14$ & $27$ &RTD$(6,23)-9$&$23$\\
          &       & $15,9$ & $28$ &RTD$(6,24)-15$&$24$\\
          &       & $21,4$ & $30$ &RTD$(6,25)-21$&$25$\\
     \hline
     $130$ & $5$ & $26,0$ &$32$ &RTD$(5,26)$&$26$\\
           & $10$ & $13,0$ &$14$ &RTD$(10,13)$&$13$\\
          & $4,5$ & $5,22$ &$33$ &RTD$(5,27)-5$&$27$\\
          &       & $10,18$  &$34$ &RTD$(5,28)-10$&$28$\\
          &       & $15,14$ &$36$ &RTD$(5,29)-15$&$29$\\
          &       & $20,10$ &$38$ &RTD$(5,30)-20$&$30$\\
          &       & $25,6$ &$39$ &RTD$(5,31)-25$&$31$\\
          &       & $30,2$ &$41$ &RTD$(5,32)-30$&$32$\\
          & $5,6$ & $2,20$ & $26$ &ownSG$(132,6)-2$ & $14$\\
          &       & $8,15$ & $27$ &RTD$(6,23)-8$&$23$\\
          &       & $14,10$ & $28$ &RTD$(6,24)-14$&$24$\\
          &       & $20,5$ & $30$ &RTD$(6,25)-20$&$25$\\
     \hline
\end{tabular}
\caption{\label{table:parallel28} SGP and SGA solutions for  $129\leq v\leq 130$. \\Block sizes  $K=\{3\}, \{5\},\{10\},\{4,5\}$, or $\{5,6\}$.}
\end{table}

\begin{table}
\centering
\begin{tabular}{|c|c|c|c|c|c|}
\hline
    $v$ & $K$ & $m_{1},m_{2}$ & $MAX$ & solution& $NR$  \\
    \hline
    $131$ & $4,5$ & $4,23$ &$33$ &RTD$(5,27)-4$&$27$\\
          &       & $9,19$  &$34$ &RTD$(5,28)-9$&$28$\\
          &       & $14,15$ &$36$ &RTD$(5,29)-14$&$29$\\
          &       & $19,11$ &$38$ &RTD$(5,30)-19$&$30$\\
          &       & $24,7$ &$39$ &RTD$(5,31)-24$&$31$\\
          &       & $29,3$ &$41$ &RTD$(5,32)-29$&$32$\\
          & $5,6$ & $1,21$ & $26$ &ownSG$(132,6)-1$ & $14$\\
          &       & $7,16$ & $27$ &RTD$(6,23)-7$&$23$\\
          &       & $13,11$ & $28$ &RTD$(6,24)-13$&$24$\\
          &       & $19,6$ & $30$ &RTD$(6,25)-19$&$25$\\
          &       & $25,1$ & $32$ &RTD$(5,26)\_Q(21,5)+1$&$22$\\
     \hline
    $132$ & $3$ & $44,0$ & $65$ & NKTS$(132)$ & $65$\\
          & $4$ & $33,0$ & $43$ & RGDD$(132,4,3)$&$43$\\
          & $6$ & $22,0$ & $22$ & ownSG$(132,6)$ & $14$\\
          & $11$ & $12,0$ & $13$ &MOLRs$(11,12)$&$6$\\
          & $4,5$ & $3,24$ & $33$ &RTD$(5,27)-3$&$27$\\
          &       & $8,20$ & $34$ &RTD$(5,28)-8$&$28$\\
          &       & $13,16$ & $36$ &RTD$(5,29)-13$&$29$\\
          &       & $18,12$ & $37$ &RTD$(5,30)-18$&$30$\\
          &       & $23,8$ & $39$ &RTD$(5,31)-23$&$31$\\
          &       & $28,4$ & $41$ &RTD$(5,32)-28$&$32$\\
          & $5,6$ & $6,17$ & $27$ & RTD$(6,23)-6$&$23$\\
          &       & $12,12$ &$28$ &RTD$(6,24)-12$&$24$\\
          &       & $18,7$ & $30$ &RTD$(6,25)-18$&$25$\\
          &       & $24,2$ &$32$ &ownSG$(156,6)-24$&$15$\\     
    \hline
\end{tabular}
\caption{\label{table:parallel29} SGP and SGA solutions for  $131\leq v\leq 132$. \\Block sizes $K=\{3\}, \{4\},\{6\},\{11\},\{4,5\}$, or $\{5,6\}$.}
\end{table}

\begin{table}
\centering
\begin{tabular}{|c|c|c|c|c|c|}
\hline
    $v$ & $K$ & $m_{1},m_{2}$ & $MAX$ & solution& $NR$  \\
    \hline
$133$&$7$ &$19,0$ &$22$ &RTD$(7,19)$ &$19$\\
     &$4,5$&$2,25$&$33$&RTD$(5,27)-2$&$27$\\
     &     &$7,21$&$34$&RTD$(5,28)-7$&$28$\\
     &     &$12,17$&$36$&RTD$(5,29)-12$&$29$\\
     &     &$17,13$&$37$&RTD$(5,30)-17$&$30$\\
     &     &$22,9$&$39$&RTD$(5,31)-22$&$31$\\
     &     &$27,5$&$41$&RTD$(5,32)-27$&$32$\\
     &     &$32,1$ &$43$ &RTD$(5,33)-32$&$33$\\
     &$5,6$&$5,18$&$27$&RTD$(6,23)-5$&$23$\\
     &     &$11,13$&$28$&RTD$(6,24)-11$&$24$\\
     &    &$17,8$ &$30$ &RTD$(6,25)-17$&$25$\\
     &     &$23,3$&$31$&ownSG$(156,6)-23$&$15$\\ 
\hline
$134$&$4,5$&$1,26$ &$33$ &RTD$(5,27)-1$&$27$\\
     &    &$6,22$ &$34$ &RTD$(5,28)-6$&$28$\\
     &    &$11,18$ &$36$ &RTD$(5,29)-11$&$29$\\
     &    &$16,14$ &$37$ &RTD$(5,30)-16$&$30$\\
     &    &$21,10$ &$39$ &RTD$(5,31)-21$&$31$\\
     &    &$26,6$ &$41$ & RTD$(5,32)-26$&$32$\\
     &    &$31,2$ &$43$ & RTD$(5,33)-31$&$33$\\
     &$5,6$&$4,19$&$27$&  RTD$(6,23)-4$&$23$\\ 
     &    &$10,14$ &$28$ &RTD$(6,24)-10$&$24$\\
     &    &$16,9$ &$30$ &RTD$(6,25)-16$&$25$\\
     &    &$22,4$ &$31$ &ownSG$(156,6)-22$&$15$\\
    \hline
\end{tabular}
\caption{\label{table:parallel30} SGP and SGA solutions for  $133\leq v\leq 134$. \\Block sizes $K=\{7\}, \{4,5\}$, or $\{5,6\}$.}
\end{table}

\begin{table}
\centering
\begin{tabular}{|c|c|c|c|c|c|}
\hline
    $v$ & $K$ & $m_{1},m_{2}$ & $MAX$ & solution& $NR$\\
    \hline 
    $135$ & $3$ & $45,0$ &$67$ &KTS$(135)$ &$67$\\
          & $5$ & $27,0$ &$33$ &RTD$(5,27)$&$27$\\
          & $9$ & $15,0$ &$16$ &ownSG$(135,9)$&$7$\\
          &$4,5$ &$5,23$ &$34$ & RTD$(5,28)-5$&$28$\\
          &      &$10,19$ &$36$ & RTD$(5,29)-10$&$29$\\
          &      &$15,15$ &$37$ & RTD$(5,30)-15$&$30$\\
          &      &$20,11$ &$39$ & RTD$(5,31)-20$&$31$\\
          &      &$25,7$ &$41$ & RTD$(5,32)-25$&$32$\\
          &      &$30,3$ &$43$ & RTD$(5,33)-30$&$33$\\
          &$5,6$ &$3,20$ &$27$ & RTD$(6,23)-3$&$23$\\
          &      &$9,15$ &$28$ & RTD$(6,24)-9$&$24$\\
          &      &$15,10$ &$30$ & RTD$(6,25)-15$&$25$\\
          &      &$21,5$ &$31$ & ownSG$(156,6)-21$&$15$\\
    \hline
    $136$ &$4$ & $34,0$ &$45$ & RBIBD$(136,4)$&$45$\\
          &$8$ & $17,0$ &$19$ &RTD$(8,17)$&$17$\\
          &$4,5$ &$4,24$ &$34$ & RTD$(5,28)-4$&$28$\\
          &       &$9,20$ &$36$ & RTD$(5,29)-9$&$29$\\
          &       &$14,16$ &$37$ & RTD$(5,30)-14$&$30$\\
          &       &$19,12$ &$39$ & RTD$(5,31)-19$&$31$\\
          &       &$24,8$ &$40$ & RTD$(5,32)-24$&$32$\\
          &       &$29,4$ &$42$ & RTD$(5,33)-29$&$33$\\
          &$5,6$ &$2,21$ &$27$ & RTD$(6,23)-2$&$23$\\
          &      &$8,16$ &$28$ & RTD$(6,24)-8$&$24$\\
          &      &$14,11$ &$30$ & RTD$(6,25)-14$&$25$\\
          &      &$20,6$ &$31$ & ownSG$(156,6)-20$&$15$\\
          &      &$26,1$ &$33$ & RTD$(6,27)-26$&$27$\\ 
    \hline
\end{tabular}
\caption{\label{table:parallel31} SGP and SGA solutions for  $135\leq v\leq 136$. \\Block sizes $K=\{3\}, \{4\},\{5\},\{8\},\{9\},\{4,5\}$, or $\{5,6\}$.}
\end{table}

\begin{table}
\centering
\begin{tabular}{|c|c|c|c|c|c|}
\hline
    $v$ & $K$ & $m_{1},m_{2}$ & $MAX$ & solution& $NR$\\
    \hline 
    $137$ & $4,5$ & $3,25$ & $34$ & RTD$(5,28)-3$&$28$\\
          &       & $8,21$ &$36$& RTD$(5,29)-8$&$29$\\
          &       & $13,17$ &$37$& RTD$(5,30)-13$&$30$\\
          &       & $18,13$ &$39$& RTD$(5,31)-18$&$31$\\
          &       & $23,9$ &$40$& RTD$(5,32)-23$&$32$\\
          &       & $28,5$ &$42$& RTD$(5,33)-28$&$33$\\
          &       & $33,1$ &$44$& RTD$(5,34)-33$&$34$\\
          & $5,6$ & $1,22$ &$27$& RTD$(6,23)-1$&$23$\\
          &       &$7,17$ &$28$ & RTD$(6,24)-7$&$24$\\
          &       &$13,12$ &$30$ & RTD$(6,25)-13$&$25$\\
          &       &$19,7$ &$31$ & ownSG$(156,6)-19$&$15$\\
          &       &$25,2$ &$33$ & RTD$(6,27)-25$&$27$\\ 
          \hline
    $138$ & $3$ & $46,0$ & $68$ & NKTS$(138)$ & $68$\\
          & $6$ & $23,0$ & $27$ & RTD$(6,23)$&$23$\\
          &$4,5$& $2,26$ & $34$ & RTD$(5,28)-2$&$28$\\
          &     & $7,22$ & $36$ & RTD$(5,29)-7$&$29$\\
          &     & $12,18$ & $37$ & RTD$(5,30)-12$&$30$\\
          &       & $17,14$ & $39$ & RTD$(5,31)-17$&$31$\\
          &       & $22,10$ & $40$ & RTD$(5,32)-22$&$32$\\
          &       & $27,6$ & $42$ & RTD$(5,33)-27$&$33$\\
          &       & $32,2$ & $44$ & RTD$(5,34)-32$&$34$\\
          & $5,6$ & $6,18$ & $28$ & RTD$(6,24)-6$&$24$\\
          &       &$12,13$ & $30$ & RTD$(6,25)-12$&$25$\\
          &       &$18,8$ & $31$ & ownSG$(156,6)-18$&$15$\\
          &       &$24,3$ & $33$ & RTD$(6,27)-24$&$27$\\   
    \hline
\end{tabular}
\caption{\label{table:parallel32} SGP and SGA solutions for  $137\leq v\leq 138$. \\Block sizes $K=\{3\}, \{6\},\{4,5\}$, or $\{5,6\}$.}
\end{table}

\begin{table}
\centering
%
\caption{\label{table:parallel33} SGP and SGA solutions for  $139\leq v\leq 140$. \\Block sizes $K=\{4\},\{5\},\{7\},\{10\},\{4,5\}$, or $\{5,6\}$.}
\end{table}

\begin{table}
\centering
%
\caption{\label{table:parallel34} SGP and SGA solutions for  $141\leq v\leq 142$. \\Block sizes $K=\{3\}, \{4,5\}$, or $\{5,6\}$.}
\end{table}

\begin{table}
\centering
%
\caption{\label{table:parallel35} SGP and SGA solutions for  $143\leq v\leq 144$. \\Block sizes $K=\{3\}, \{4\},\{6\},\{8\},\{9\},\{11\},\{12\},\{4,5\}$, or $\{5,6\}$.}
\end{table}

\begin{table}
\centering
%
\caption{\label{table:parallel36} SGP and SGA solutions for  $145\leq v\leq 146$. \\Block sizes $K=\{5\},\{4,5\}$, or $\{5,6\}$.}
\end{table}

\begin{table}
\centering
%
\caption{\label{table:parallel37} SGP and SGA solutions for  $147\leq v\leq 148$. \\Block sizes $K=\{3\}, \{4\},\{7\},\{4,5\}$, or $\{5,6\}$.}
\end{table}

\begin{table}
\centering
%
\caption{\label{table:parallel38} SGP and SGA solutions for  $149\leq v\leq 150$. \\Block sizes $K=\{3\}, \{5\},\{6\},\{10\},\{4,5\}$, or $\{5,6\}$.}
\end{table}

\begin{table}
\centering
%
\caption{\label{table:additional} Some example SGP and SGA solutions for  $v>150$ and a single block size.
}
\end{table}

\end{document}